\documentclass[12pt]{amsart} 
\usepackage[english]{babel}
\usepackage[latin1]{inputenc}
\usepackage{amssymb}
\usepackage{amsmath}
\usepackage{amscd}
\usepackage{latexsym}
\usepackage{amsthm}
\usepackage{mathrsfs}
\usepackage{mathabx}
\usepackage{color}
\theoremstyle{plain}
\newtheorem{thm}{Theorem}[section]
\newtheorem{cor}[thm]{Corollary}
\newtheorem{lem}[thm]{Lemma}
\newtheorem{prop}[thm]{Proposition}

\newtheorem{defn}[thm]{Definition}
\newtheorem{oss}[thm]{Remark}
\newtheorem{esem}[thm]{Example}
\newtheorem{notation}[thm]{Notation}

\begin{document}

\title[Moduli spaces of sheaves on K3 surfaces]{Irreducible symplectic varieties from moduli spaces of sheaves on K3 and Abelian surfaces}
\author{Arvid Perego}
\address{Arvid Perego: Dipartimento di Matematica dell'Universit\`a di Genova, 16146 Genova, Italy. E-mail: perego@dima.unige.it}
\author{Antonio Rapagnetta}
\address{Antonio Rapagnetta: Dipartimento di Matematica dell'Universit\`a di Roma II - Tor Vergata, 00133 Roma, Italy. E-mail: rapagnet@axp.mat.uniroma2.it}

\subjclass[2020]{14D20,14J60,14J28,14J42}

\keywords{moduli spaces of sheaves, K3 surfaces, irreducible symplectic varieties}

\begin{abstract}
We show that the moduli spaces of sheaves on a projective K3 surface are irreducible symplectic varieties, and that the same holds for the fibers of the Albanese map of moduli spaces of sheaves on an Abelian surface.
\end{abstract}

\maketitle

\section{Introduction and main results}

A \textit{holomorphic symplectic form} on a complex manifold $X$ is an everywhere nondegenerate, closed, holomorphic $2-$form on $X$. A complex manifold admitting a holomorphic symplectic form is called \textit{holomorphic symplectic manifold}. We let $h^{p,0}(X)$ be the dimension of the vector space $H^{0}(X,\Omega_{X}^{p})$.

A connected compact K\"ahler manifold $X$ is an \textit{irreducible symplectic manifold} if it is holomorphic symplectic, simply connected and $h^{2,0}(X)=1$. In particular, an irreducible symplectic manifold has even complex dimension and trivial canonical bundle. 

By the Bogomolov Decomposition Theorem, irreducible symplectic manifolds are one of the three types of manifolds which are building blocks for compact K\"ahler manifolds with numerically trivial canonical bundle. There are very few known deformation classes of irreducible symplectic manifolds, namely:

\begin{enumerate}
 \item a compact, connected smooth complex surface is an irreducible symplectic manifold if and only if it is a K3 surface;
 \item if $S$ is a K3 surface and $n\in\mathbb{N}$, $n\geq 2$, the Hilbert scheme $Hilb^{n}(S)$ of $n$ points on $S$ is an irreducible symplectic manifold of dimension $2n$ (see \cite{B}, Th\'eor\`eme 3 and Proposition 6);
 \item if $T$ is a $2-$dimensional complex torus and $n\in\mathbb{N}$, $n\geq 2$, the generalized Kummer variety $Kum^{n}(T)$ is an irreducible symplectic manifold of dimension $2n$ (see \cite{B}, Th\'eor\`eme 4 and Proposition 8);
 \item there are two more deformation classes: $OG_{6}$, in dimension 6, and $OG_{10}$, in dimension 10 (see \cite{OG2}, \cite{OG3}).
\end{enumerate}

A possible way to obtain new examples of varieties behaving like irreducible symplectic manifolds is to enlarge the family of varieties we are considering by including singular varieties. This is a very natural step, in particular in view of the Minimal Model Program.

Indeed, if $X$ is a connected complex projective manifold with $\kappa(X)=0$, if the MMP works for $X$ then it produces a birational map $X\dashrightarrow Y$, where $Y$ has terminal singularities and nef canonical divisor. Assuming the abundance conjecture, we get that a multiple of $K_{Y}$ is trivial. So, for the classification of projective varieties whose Kodaira dimension is 0 it is central to extend the Bogomolov decomposition to normal projective varieties having terminal singularities and torsion (i. e. numerically trivial by Theorem 8.2 of \cite{Kaw}) canonical divisor.

A singular version of the Bogomolov Decomposition Theorem has been recently obtained. For singular projective varieties with klt singularities this is Theorem 1.5 of \cite{HP} (whose proof is the combination of several results contained in \cite{GGK}, \cite{D}, \cite{GKP}, \cite{DG}), then extended to compact K\"ahler spaces with log-terminal singularities and numerically trivial canonical bundle by Theorem A of \cite{BGL}. 

The role played by irreducible symplectic manifolds in the Bogomolov Decomposition Theorem is played in these generalizations by \textit{irreducible symplectic varieties}, whose definition was first given in \cite{GKP}. We will present the definition of an irreducible symplectic variety only in the projective setting, since this is the one we will need in the present paper (for a more general definition, see \cite{BGL}). 

We need the following notation: if $X$ is a normal algebraic variety and $X_{reg}$ is the smooth locus of $X$ whose open embedding in $X$ is $j:X_{reg}\longrightarrow X$, for every $p\in\mathbb{N}$ such that $0\leq p\leq\dim(X)$ we let $$\Omega_{X}^{[p]}:=j_{*}\Omega^{p}_{X_{reg}}=\big(\wedge^{p}\Omega_{X}\big)^{**},$$whose global sections are called \textit{reflexive $p-$forms} on $X$. A reflexive $p-$form on $X$ is then a holomorphic $p-$form on $X_{reg}$. 

If $f:Y\longrightarrow X$ is a finite, dominant morphism between irreducible normal algebraic varieties, there is a morphism of reflexive sheaves $f^{*}\Omega_{X}^{[p]}\longrightarrow\Omega_{Y}^{[p]}$, induced by the usual pull-back morphism of forms on the smooth loci, giving a morphism $f^{[*]}:H^{0}(X,\Omega_{X}^{[p]})\longrightarrow H^{0}(Y,\Omega_{Y}^{[p]})$, called \textit{reflexive pull-back morphism}.

We first recall the definitions of a symplectic form and a symplectic variety (see \cite{B2}).

\begin{defn}
{\rm Let $X$ be a normal algebraic variety.
\begin{enumerate}
 \item A \textit{symplectic form} on $X$ is a closed reflexive $2-$form $\sigma$ on $X$ which is nondegenerate at each point of $X_{reg}$.
 \item If $\sigma$ is a symplectic form on $X$, the pair $(X,\sigma)$ is a \textit{symplectic variety} if for every resolution $f:\widetilde{X}\longrightarrow X$ of the singularities of $X$, the holomorphic symplectic form $\sigma_{reg}:=\sigma_{|X_{reg}}$ extends to a holomorphic $2-$form on $\widetilde{X}$.
 \item If $(X,\sigma)$ is a symplectic variety and $f:\widetilde{X}\longrightarrow X$ is a resolution of the singularities over which $\sigma_{reg}$ extends to a holomorphic symplectic form, we say that $f$ is a \textit{symplectic resolution}.
\end{enumerate}}
\end{defn}

We now define irreducible symplectic varieties following \cite{GKP}. Recall that if $X$ and $Y$ are two irreducible normal projective varieties, a \textit{finite quasi-\'etale morphism} $f:Y\longrightarrow X$ is a finite morphism which is \'etale in codimension one.

\begin{defn}
\label{defn:irrvar}
{\rm An \textit{irreducible symplectic variety} is a normal projective variety $X$ with canonical singularities that has a symplectic form $\sigma\in H^{0}(X,\Omega_{X}^{[2]})$ such that for every finite quasi-\'etale morphism $f:Y\longrightarrow X$ the exterior algebra of reflexive forms on $Y$ is spanned by $f^{[*]}\sigma$.}
\end{defn}

\begin{oss}
\label{oss:sicuro}
{\rm If $X$ is a normal projective variety, by definition of $\Omega_{X}^{[p]}$ we have $H^{0}(X,\Omega_{X}^{[p]})=H^{0}(X_{reg},\Omega_{X_{reg}}^{p})$. Theorem 1.4 of \cite{GKKP} implies that if $X$ is a quasi-projective with klt singularities and $\pi:\widetilde{X}\longrightarrow X$ is a log-resolution, then for every $p\in\mathbb{N}$ such that $0\leq p\leq\dim(X)$ the sheaf $\pi_{*}\Omega^{p}_{\widetilde{X}}$ is reflexive. This implies in particular (see Observation 1.3 therein) that $H^{0}(X,\Omega_{X}^{[p]})\simeq H^{0}(\widetilde{X},\Omega^{p}_{\widetilde{X}})$.}
\end{oss}

The definition of irreducible symplectic variety is motivated by the description of the algebra of holomorphic forms of an irreducible symplectic manifold, which is spanned by a holomorphic symplectic form (see Proposition 3 in \cite{B}).

By Corollary 13.3 of \cite{GGK} an irreducible symplectic variety $X$ is simply connected, so the $\mathbb{Z}-$module $H^{2}(X,\mathbb{Z})$ is free. Moreover, the fact that irreducible symplectic varieties are simply connected together with the Bogomolov Decomposition Theorem imply that smooth irreducible symplectic varieties are irreducible symplectic manifolds.

\begin{oss}
\label{oss:gkk}
{\rm A symplectic resolution $Y$ of an irreducible symplectic variety $X$ is an irreducible symplectic manifold. Indeed $X$ is klt and simply connected, so $Y$ is simply connected as well (see \cite{T}), and by Theorem 1.4 of \cite{GKKP} we have $h^{0}(Y,\Omega_{Y}^{2})=h^{0}(X,\Omega_{X}^{[2]})=1$, since $X$ is irreducible symplectic. Anyway there are symplectic varieties having an irreducible symplectic manifold as a symplectic resolution, but which are not irreducible symplectic varieties, as in the following:}
\end{oss}

\begin{esem}
\label{esem:sym}
{\rm If $X$ is a K3 surface and $m\geq 2$, then $Y=Sym^{m}(X)$ has $Hilb^{m}(X)$ as symplectic resolution, but it is not an irreducible symplectic variety since $X^{m}$ is a finite quasi-\'etale cover of $Y$.}
\end{esem}

Examples of irreducible symplectic varieties in dimension 4 are known (see \cite{Pe} for an overview). Among them we cite the partial resolution of the quotient of $Hilb^{2}(S)$ for $S$ a K3 surface (resp. $Kum^{2}(T)$ for $T$ an Abelian surface) by the action of a symplectic involution (\cite{Men1}, \cite{KaMe}), and the quotients of $Hilb^{2}(S)$ by the action of a symplectic automorphism of order 3, 5, 7 or 11 (\cite{Men2}, \cite{Men3}).

The aim of this work is to provide a wide family of examples of irreducible symplectic varieties in higher dimension.  

\subsection{Notation and main results of the paper} 

The aim of the present paper is to provide families of irreducible symplectic varieties using moduli spaces of sheaves on K3 or Abelian surfaces.

In the following, $S$ will be a projective K3 surface or an Abelian surface, and we let $\widetilde{H}(S,\mathbb{Z}):=H^{2*}(S,\mathbb{Z})$. An element $v\in \widetilde{H}(S,\mathbb{Z})$ will be written $v=(v_{0},v_{1},v_{2})$, where $v_{i}\in H^{2i}(S,\mathbb{Z})$, and $v_{0},v_{2}\in\mathbb{Z}$. It will be called \textit{Mukai vector} if $v_{0}\geq 0$, $v_{1}\in NS(S)$ and if $v_{0}=0$, then either $v_{1}$ is the first Chern class of an effective divisor, or $v_{1}=0$ and $v_{2}>0$. Moreover, we let $\rho(S)$ be the rank of the N\'eron-Severi group of $S$. 

Recall that $\widetilde{H}(S,\mathbb{Z})$ has a pure weight-two Hodge structure and a compatible lattice structure given by  the Mukai pairing $(.,.)$ (see \cite{HL}, Definitions 6.1.5 and 6.1.11). We let $v^{2}:=(v,v)$ for every Mukai vector $v$, and we call $\widetilde{H}(S,\mathbb{Z})$ the \textit{Mukai lattice} of $S$.

If $\mathscr{F}$ is a coherent sheaf on $S$, we define its \textit{Mukai vector} as $$v(\mathscr{F}):=ch(\mathscr{F})\sqrt{td(S)}=(rk(\mathscr{F}),c_{1}(\mathscr{F}),ch_{2}(\mathscr{F})+\epsilon(S)rk(\mathscr{F})),$$where $\epsilon(S):=1$ if $S$ is K3, and $0$ if $S$ is Abelian.

Let now $v$ be a Mukai vector on $S$ and suppose that $H$ is a polarization which is general with respect to $v$ (see Definition 2.9). We write $M_{v}(S,H)$ (resp. $M_{v}^{s}(S,H)$) for the moduli space of Gieseker $H-$semistable (resp. $H-$stable) sheaves on $S$ with Mukai vector $v$.
 
If $S$ is Abelian and $v^{2}>0$, we have a dominant isotrivial fibration $a_{v}:M_{v}(S,H)\longrightarrow S\times\widehat{S}$ (see section 4.1 of \cite{Y2}), where $\widehat{S}$ is the dual of $S$. We let $K_{v}(S,H):=a_{v}^{-1}(0_{S},\mathscr{O}_{S})$, and $K_{v}^{s}(S,H):=K_{v}(S,H)\cap M^{s}_{v}(S,H)$. The morphism $a_{v}$ is known to be the Albanese morphism of $M_{v}(S,H)$ if $v$ is primitive (see Theorem 0.1 of \cite{Y2}): we will show that this holds in all the cases we will consider (see Corollary \ref{cor:albanese}).

If no confusion on $S$ and $H$ is possible, we drop them from the notation. Moreover, we will always write $v=mw$, where $m\in\mathbb{N}$ and $w$ is a primitive Mukai vector on $S$. 

If $M_{v}^{s}\neq\emptyset$, then it is a holomorphic symplectic quasi-projective manifold of dimension $v^{2}+2$ (see \cite{M1}).

If $m=1$ and $S$ is K3, then $M_{v}^{s}\neq\emptyset$ if and only if $v^{2}\geq-2$ (see Theorem 0.1 of \cite{Y1}). If $S$ is Abelian, then $M_{v}^{s}\neq\emptyset$ if and only if $v^{2}\geq 0$ (see Theorem 0.1 of \cite{Y2}, and compare with section 2.4 of \cite{KLS}). If $w^{2}>0$, then $M_{v}$ and $K_{v}$ are normal, irreducible projective varieties (see Theorem 4.4 of \cite{KLS} and Remark A.1 of \cite{PR2}).

If $v^{2}\leq 0$ we have a precise description of $M_{v}$ and $K_{v}$.
\begin{enumerate}
 \item If $v^{2}<0$ and $S$ is K3, then $M_{v}$ is either empty or a point (see \cite{M2}). If $S$ is Abelian, then $M_{v}=\emptyset$ (see \cite{Y2}). 
 \item If $v^{2}=0$ and $S$ is K3, then either $M_{v}$ is a K3 surface (if $m=1$, see \cite{M2}) or a symmetric product of a K3 surface (see section 1 of \cite{KLS}), in which case $M_{v}$ is not irreducible symplectic (see Example \ref{esem:sym}). 
 \item If $v^{2}=0$ and $S$ is Abelian, then either $M_{v}$ is an Abelian surface $A$ (if $m=1$, see \cite{M2}) or a symmetric product of an Abelian surface $A$ (see section 1 of \cite{KLS}). Then $M_{v}$ is not simply connected, and the fiber of the sum morphism $M_{v}\longrightarrow A$ is either a point (if $m=1$) or a symplectic variety which is not irreducible symplectic (the proof of this is similar to Example \ref{esem:sym}).
\end{enumerate}

Because of this, we will only be interested in Mukai vectors $v$ such that $v^{2}>0$. In what follows we will need the following:

\begin{defn}
{\rm Let $m,k\in\mathbb{N}$, $m,k>0$. A Mukai vector $v$ on $S$ will be said of \textit{type} $(m,k)$ if $v=mw$ for a primitive Mukai vector $w\in\widetilde{H}(S,\mathbb{Z})$ such that $w^{2}=2k$.}
\end{defn}

If $S$ a projective K3 surface or an Abelian surface, $v$ is a Mukai vector of type $(m,k)$ on $S$ and $H$ is a polarization on $S$ that is general with respect to $v$, then $M_{v}$ is a nonempty, irreducible, normal projective variety of dimension $2m^{2}k+2$ (see Theorem 4.4 of \cite{KLS}), which is symplectic and whose smooth locus is $M^{s}_{v}$. If $S$ is Abelian and $(m,k)\neq(1,1)$, then $K_{v}$ is a nonempty, irreducible, normal projective variety of dimension $2m^{2}k-2$, which is symplectic and whose regular locus is $K^{s}_{v}$. If $(m,k)=(1,1)$, then $M_{v}$ is isomorphic to $S\times\widehat{S}$ and $K_{v}$ is a point.

The first result we will prove is the following:

\begin{thm}
\label{thm:main}Let $m,k\in\mathbb{N}$, $m,k>0$, and for $i=1,2$ let $S_{i}$ be a projective K3 or Abelian surface, $v_{i}$ a Mukai vector on $S_{i}$ of type $(m,k)$ and $H_{i}$ a polarization on $S_{i}$ which is general with respect to $v_{i}$. 
\begin{enumerate}
 \item If $S_{1}$ and $S_{2}$ are both K3 surfaces or both Abelian surfaces, then $M_{v_{1}}(S_{1},H_{1})$ and $M_{v_{2}}(S_{2},H_{2})$ are deformation equivalent, and the deformation is locally trivial.
 \item If $S_{1}$ and $S_{2}$ are two Abelian surfaces, then $K_{v_{1}}(S_{1},H_{1})$ and $K_{v_{2}}(S_{2},H_{2})$ are deformation equivalent, and the deformation is locally trivial.
\end{enumerate}
\end{thm}

\begin{oss}
\label{oss:loctriv}
{\rm In the statement, and all along the paper, we will say that a morphism $f:\mathscr{X}\longrightarrow T$ of complex varieties is a \textit{locally trivial deformation} of a complex variety $X$ if $T$ is connected and $f$ is a proper flat morphism verifying the two following conditions (see \cite{FK}):
\begin{enumerate}
	\item there is a point $t_{0}\in T$ such that $X=f^{-1}(t_{0})$,
\item for every $x\in \mathscr{X}$ the deformation germ $(\mathscr{X},x)\longrightarrow(T,f(x))$ is isomorphic to the trivial deformation $(f^{-1}(f(x)),x)\times(T,f(x))$ of the germ $(f^{-1}(f(x)),x)$.
\end{enumerate}}
\end{oss}

\begin{oss}
\label{oss:defotype}
{\rm The deformation relating $M_{v_{1}}(S_{1},H_{1})$ and $M_{v_{2}}(S_{2},H_{2})$ in Theorem \ref{thm:main} is obtained using only deformations of the moduli spaces induced by deformations of the base surfaces $S_{i}$ together with the Mukai vectors and the polarizations (see section 2.3 for the definition) and isomorphisms between moduli spaces induced by Fourier-Mukai transforms. In partcular $M_{v_{1}}(S_{1},H_{1})$ and $M_{v_{2}}(S_{2},H_{2})$ are locally trivially deformation equivalent in the algebraic category.}
\end{oss}

As a consequence of Theorem \ref{thm:main}, starting from K3 surfaces (resp. Abelian surfaces) we get a single locally trivial deformation class for every pair $(m,k)$ of strictly positive integers. 

The proof of Theorem \ref{thm:main} is the content of section 2 of the present paper. For $m=1$ it is due to several authors (see \cite{M1}, \cite{B}, \cite{OG1}, \cite{Y1}, \cite{Y2}). For $(m,k)=(2,1)$, the proof of Theorem \ref{thm:main} is in \cite{PR}. The deformation equivalence in point (1) of Theorem \ref{thm:main} is basically due to Yoshioka: for Mukai vectors with positive rank it is Proposition 3.6 of \cite{Y6}; the rank 0 case is not explicitly stated, but can be obtained as in Corollary 3.5 of \cite{Y7}. As it is an important result which plays a central role in our paper, we decided to include a complete proof. The local triviality of the deformation follows from the Main Theorem of \cite{N2}. 

Yoshioka's original proof of the deformation equivalence involves two main technical tools: deformations of K3 or Abelian surfaces inducing deformations of the moduli spaces, and Fourier-Mukai transforms. As a by-product of it one gets nonemptyness and normality of the moduli spaces. Based on his proof of the deformation equivalence of the moduli spaces, and using a different argument to deal with particular cases, Yoshioka proves also their irreducibility (see Theorem 3.18 of \cite{Y6}).

The proof of point (1) of Theorem \ref{thm:main} we propose is a re-elaborated version of Yoshioka's proof that we tried to keep as elementary as possible. It uses the same tools together with Theorem 4.4 of \cite{KLS}, which proves the irreducibility and the normality of the moduli spaces independently of \cite{Y6} and \cite{Y7} (and which implies the irreducibility and the normality of the $K_{v}$'s, as shown in Remark A.1 of \cite{PR2}). 

We remark that the proof of Theorem 4.4 of \cite{KLS} uses the fact that if $v$ is primitive and $v^{2}\geq 0$, then $M_{v}\neq\emptyset$, which was proved by Yoshioka in \cite{Y1} and \cite{Y2}, independently of \cite{Y6} and \cite{Y7} (compare with section 2.4 of \cite{KLS}).

The use of Theorem 4.4 of \cite{KLS} allows us to minimize the technicalities  involved in the proof of the deformation equivalence of moduli spaces. The only Fourier-Mukai transforms we will use are those corresponding to tensorization with line bundles: the one whose kernel is the ideal of the diagonal (for K3 surfaces), and the one whose kernel is the Poincar\'e bundle (for Abelian surfaces), and we only need to check the preservation of the semistability under these functors in the most natural direction (see section 2.4).

The aim of section 3 is to show the following, which is the main result of this paper:

\begin{thm}
\label{thm:mainb1}
Let $m,k\in\mathbb{N}$, $m,k>0$, $S$ a projective K3 or Abelian surface, $v$ a Mukai vector on $S$ of type $(m,k)$ and $H$ a polarization on $S$ which is general with respect to $v$. 
\begin{enumerate}
 \item If $S$ is K3, then $M_{v}(S,H)$ is an irreducible symplectic variety.
 \item If $S$ is Abelian and $(m,k)\neq(1,1)$, then $K_{v}(S,H)$ is an irreducible symplectic variety.  
\end{enumerate}
\end{thm}

Theorem \ref{thm:mainb1} provides an answer to Question 14.10 of \cite{GGK}: it implies that, for general polarizations, all moduli spaces of sheaves on a projective K3 surface (and all the fibers of the Albanese morphism of moduli spaces of sheaves on Abelian surfaces) are irreducible symplectic varieties, with the only exception of symmetric products. 

The examples of irreducible symplectic varieties given in Theorem  \ref{thm:mainb1}  
naturally fall  into three  main cases:
\begin{enumerate}
 \item if $S$ is K3, then $M_{v}$ is smooth if and only if $m=1$ (and deformation equivalent to $Hilb^{k+1}(S)$). If $S$ is Abelian, then $K_{v}$ is smooth if and only if $m=1$ (it is a point if $k=1$, a K3 surface if $k=2$, and deformation equivalent to $Kum^{k-1}(S)$ if $k\geq 3$); 
 \item if $S$ is K3, then $M_{v}$ has a symplectic resolution if and only if $(m,k)=(2,1)$ (which is in the deformation class $OG_{10}$). If $S$ is Abelian, then $K_{v}$ has a symplectic resolution if and only if $(m,k)=(2,1)$ (which is in the deformation class $OG_{6}$);
  \item in all other cases $M_{v}$ and $K_{v}$ have terminal singularities and no symplectic resolutions. Indeed, both $M_{v}$ and $K_{v}$ are normal, irreducible projective varieties having a symplectic form on their smooth loci, and their singular loci have codimension at least equal to 4 (see for example \cite{KLS}). The main result of \cite{Fle} implies that they are symplectic varieties, and by Corollary 1 of \cite{N1} they have terminal singularities. As shown in  \cite{KLS}, $M_{v}$ and $K_{v}$ are locally factorial and do not admit any symplectic resolution.
\end{enumerate}

The proof of Theorem \ref{thm:mainb1} uses Theorem \ref{thm:main} to reduce to a surface $S$ such that $NS(S)=\mathbb{Z}\cdot h$, where $h$ is the first Chern class of an ample divisor $H$ with $H^{2}=2k$. Taking $v=m(0,h,0)$, if $S$ is a K3 surface we show that $M_{v}$ and $M^{s}_{v}$ are simply connected; if $S$ is Abelian, we show that $K_{v}$ and $K^{s}_{v}$ are simply connected (if $(m,k)\neq(2,1)$). We then calculate the numbers $h^{0}(M_{v},\Omega_{M_{v}}^{[p]})$ and $h^{0}(K_{v},\Omega_{K_{v}}^{[p]})$ by comparing them with $h^{0}(M_{v'},\Omega_{M_{v'}}^{p})$ and $h^{0}(K_{v'},\Omega_{K_{v'}}^{p})$ where $v'$ is the primitive Mukai vector $v'=(0,mh,1-m^{2}k)$. 

\begin{oss}
{\rm The irreducible symplectic varieties we get with Theorem \ref{thm:mainb1} have all simply connected smooth locus, up to one exception, namely when $S$ is an Abelian surface and $(m,k)=(2,1)$: in this case the fundamental group of $K^{s}_{v}$ is $\mathbb{Z}/2\mathbb{Z}$ (see section 4 of \cite{MRS}, or Theorem \ref{thm:simconab}). In any case, all the irreducible symplectic varieties we obtain have smooth locus with finite fundamental group.}
\end{oss}

\begin{oss}
{\rm A natural open question concerns the computation of the fundamental invariants, i.e. the Beauville-Bogomolov form and the Fujiki constant for $M_{v}$ and $K_{v}$ in the case where they do not admit a symplectic resolution.  
We treat this problem in \cite{PR3}.}
\end{oss}

\section{Deformations of moduli spaces}

In this section we study how moduli spaces vary under deformations and isomorphisms. In Section 2.1 we recall the notion of polarization which is general with respect to a Mukai vector $v$ and some notions related to that: $v-$generic polarizations, $v-$walls and $v-$chambers. Section 2.2 is devoted to the morphism $a_{v}$ in the case of Abelian surfaces. In Section 2.3 we introduce deformations of moduli spaces induced by deformations of the datum of a surface $S$, a Mukai vector $v$ on $S$ and a polarization $H$ along smooth, connected varieties. In section 2.4 we study isomorphisms between moduli spaces coming from Fourier-Mukai transforms. 

Section 2.5 is devoted to the proof of Theorem \ref{thm:main}, which is the main result of this chapter. Our proof of Theorem \ref{thm:main} is heavily based on Theorem 4.4 of \cite{KLS} which asserts that if $v$ is a Mukai vector on a projective K3 or Abelian surface $S$ and $H$ is a $v-$generic polarization, then $M_{v}(S,H)$ is a normal, irreducible projective variety of the expected dimension. Theorem 4.4 of \cite{KLS} is based on the nonemptyness of moduli spaces of sheaves for generic polarizations and primitive Mukai vectors of positive square (Theorems 0.1 and 8.1 of \cite{Y2}, compare with section 2.4 of \cite{KLS}).

These assumptions could be avoided (following Yoshioka) using Theorem 1.17 of \cite{M1} and stronger versions of Propositions \ref{prop:yoshi2} and \ref{prop:yoshi1}.

\subsection{General and generic polarizations}

We recall the definition of $v-$generic polarization introduced in \cite{HL} and  \cite{Y2}, that we will use all along the paper, and the notion of polarization which is general with respect to $v$ introduced in \cite{Y7}.

\subsubsection{Generic polarizations}

In what follows, $S$ will always denote a projective K3 or Abelian surface, and $v=(v_{0},v_{1},v_{2})$ a Mukai vector on $S$. We will furthermore suppose that when $\rho(S)>1$, if $v_{0}=0$ then $v_{2}\neq 0$ (the case of $v=(0,v_{1},0)$ will be briefly discussed at the beginning of section 2.1.2, see Example \ref{esem:noqvg}).

We associate to each Mukai vector $v$ of this form a set $W_{v}$ of divisors on $S$, whose definition depends on $v_{0}$: the case $v_{0}>0$ will be different from the case $v_{0}=0$.  

If $v_{0}>0$, we let $$|v|=\frac{v_{0}^{2}}{4}(v,v)+\frac{v_{0}^{2\epsilon(S)+2}}{2},$$where recall that $\epsilon(S)=1$ if $S$ is K3 and 0 if $S$ is Abelian. The rational number $|v|$ only depends on $v_{0}$ and $v^{2}$. If $|v|>0$, we define $$W_{v}:=\{D\in NS(S)\,|\,-|v|\leq D^{2}<0\},$$and we let $W_{v}:=\emptyset$ if $|v|=0$. We notice that if $m,k\in\mathbb{N}$, $m,k>0$, and $v$ is a Mukai vector of type $(m,k)$, then $|v|>0$.

If $v_{0}=0$, for every pure sheaf $E$ of Mukai vector $v$, and $0\neq F\subseteq E$ of Mukai vector $u=(0,u_{1},u_{2})$, the \textit{divisor associated to the pair} $(E,F)$ is defined as $D:=u_{2}v_{1}-v_{2}u_{1}$. The set $W_{v}$ will then be the set of the non-numerically trivial divisors associated to all the possible pairs of this type.

A primitive ample divisor $H$ on $S$ will be called \textit{polarization}\footnote{By a slight abuse of notation, the line bundle $\mathcal{O}_{S}(H)$ will be usually denoted $H$, and will still be called polarization.}. The set $W_{v}$ is used to define the notion of $v-$generic polarization as follows:

\begin{defn}
{\rm A polarization $H$ is $v-$\textit{generic} if $H\cdot D\neq 0$ for every $D\in W_{v}$.}
\end{defn}

If $\rho(S)=1$, then the ample generator of the Picard group of $S$ is $v-$generic for every $v$. If $\rho(S)\geq 2$, there can be polarizations which are not $v-$generic, and to characterize them we introduce $v-$walls and $v-$chambers. We let $Amp(S)$ be the ample cone of $S$.

\begin{defn}
{\rm If $D\in W_{v}$, the $v-$\textit{wall associated to} $D$ is $$D^{\perp}:=\{\alpha\in Amp(S)\,|\,D\cdot\alpha=0\}.$$}
\end{defn}

The $v-$wall associated to $D\in W_{v}$ is a hyperplane in $Amp(S)$. If $v_{0}>0$, by Theorem 4.C.2 of \cite{HL} the set of $v-$walls is locally finite in $Amp(S)$. If $v_{0}=0$, it is even finite (see section 1.4 of \cite{Y2}).

\begin{defn}
{\rm Suppose that $\rho(S)\geq 2$. A connected component of $Amp(S)\setminus\bigcup_{D\in W_{v}}D^{\perp}$ is called $v-$\textit{chamber}.}
\end{defn}

By definition, a polarization is $v-$generic if and only if it lies in a $v-$chamber. Since the family of $v-$walls is locally finite in $Amp(S)$, it follows that a $v-$generic polarization exists for every Mukai vector $v$ verifying the conditions above.

\begin{oss}
\label{oss:quasivgen}
{\rm If $H$ is $v-$generic and $E$ is a $H-$semistable sheaf with Mukai vector $v$, then any $H-$destabilizing subsheaf of $E$ has Mukai vector of the form $pv$ for some $p\in\mathbb{Q}$. In particular, if $v$ is primitive and $H$ is $v-$generic, any $H-$semistable sheaf of Mukai vector $v$ is $H-$stable (compare with section 2.4 of \cite{KLS}).}
\end{oss}

An important property of generic polarizations is that changing polarization inside a $v-$chamber does not affect the moduli space. More precisely, we have the following (see \cite{Q}, or section 4.C in \cite{HL}):

\begin{prop}
\label{prop:chamber}Suppose that $\rho(S)\geq 2$ and that $v$ is a Mukai vector on $S$. Let $\mathcal{C}$ be a $v-$chamber, and suppose that $H,H'\in\mathcal{C}$. A sheaf of Mukai vector $v$ is $H-$(semi)stable if and only if it is $H'-$(semi)stable. As a consequence, we have natural identifications $M_{v}(S,H)=M_{v}(S,H')$ and $M^{s}_{v}(S,H)=M^{s}_{v}(S,H')$.
\end{prop}

We conclude this section with the behaviour of $v-$genericity with respect to the tensorization with a line bundle. If $v$ is a Mukai vector on $S$ and $L\in Pic(S)$, we let $v_{L}:=v\cdot ch(L)$. If $L=\mathcal{O}_{S}(D)$ for some divisor $D$, then we let $v_{D}:=v_{L}$. Notice that if $E$ is a sheaf such that $v(E)=v$, then $v(E\otimes L)=v_{L}$. 

\begin{lem}
\label{lem:tensogen}
Let $v$ be a Mukai vector and $H$ a polarization on $S$.
\begin{enumerate}
 \item If $v=(r,\xi,a)$ with $r>0$ and $L\in Pic(S)$, we have that $H$ is $v-$generic if and only if it is $v_{L}-$generic.
 \item If $v=(0,\xi,a)$ where $a\neq 0$ and $d\in\mathbb{Z}$ is such that $d\neq-\frac{a}{\xi\cdot H}$, we have that $H$ is $v-$generic if and only if it is $v_{dH}-$generic.
\end{enumerate}
\end{lem}

\proof If $v=(r,\xi,a)$ with $r>0$, notice that $$v_{L}=(r,\xi+rc_{1}(L),a+\xi\cdot L+rL^{2}/2).$$Then $v$ and $v_{L}$ have the same rank, and as easily seen $v^{2}=v_{L}^{2}$. Hence $|v|=|v_{L}|$, so $W_{v}=W_{v_{L}}$, and we are done.

If $v=(0,\xi,a)$, notice that $v_{dH}=(0,\xi,a+d\xi\cdot H)$. There is a bijection between $W_{v}\cup\{0\}$ and $W_{v_{dH}}\cup\{0\}$ mapping the divisor associated to a pair $(E,F)$ to the one associated to $(E\otimes\mathcal{O}_{S}(dH),F\otimes\mathcal{O}_{S}(dH))$. Indeed, if $D\in W_{v}$ is associated to $(E,F)$, and if $v(F)=(0,\zeta,b)$, we get $D=b\xi-a\zeta$. The divisor associated to $(E\otimes\mathcal{O}_{S}(dH),F\otimes\mathcal{O}_{S}(dH))$ is $$D'=D+d(\xi\cdot H)\zeta-d(\zeta\cdot H)\xi,$$and the bijection maps $D$ to $D'$ (and conversely). Notice that $D\cdot H=D'\cdot H$, hence $H$ is $v-$generic if and only if it is $v_{dH}-$generic.\endproof

\subsubsection{General polarizations}

The definition of $v-$generic polarization makes perfect sense even for Mukai vectors of type $v=(0,v_{1},0)$, but it is not well-adapted to our goals. Indeed, if $v=(0,v_{1},0)$, by defining $W_{v}$ as before we see that $D\in W_{v}$ is of the form $D=bv_{1}$ for some $b\neq 0$. As $v_{1}$ is the first Chern class of an effective divisor, we get $H\cdot D\neq 0$ for all $D\in W_{v}$, and hence every polarization would be $v-$generic.
 
Now, the definition of $v-$genericity is motivated by the fact that if $v$ is a primitive Mukai vector and $H$ is $v-$generic, then a $H-$semistable sheaf with Mukai vector $v$ is $H-$stable: this holds if $v=(0,v_{1},v_{2})$ where $v_{2}\neq 0$, or if $v=(0,v_{1},0)$ and $\rho(S)=1$ (as a consequence of Remark \ref{oss:quasivgen}), but it is no longer true if $v=(0,v_{1},0)$ and $\rho(S)\geq 2$ as the following example shows (see Lemma 1.2 of \cite{Y2}).

\begin{esem}
\label{esem:noqvg}
{\rm Let $S$ be a K3 surface with $NS(S)=\mathbb{Z}\cdot c'\oplus\mathbb{Z}\cdot c''$, where $c'$ and $c''$ are the classes of two irreducible effective curves $C'$ and $C''$, e.g. $S$ is a generic quartic surface containing a line.  We let $j'$ and $j''$ be the inclusions of $C'$ and $C''$ respectively into $S$. Denote by $v$ the primitive Mukai vector $(0,c'+c'',0)$ and let $M'$ and $M''$ be line bundles on $C'$ and $C''$ with Euler characteristic zero. The sheaves  $j'_{*}M'$ and $j''_{*}M''$ are $H-$stable with respect to any polarization $H$, and we have $v(j'_{*}M')=(0,c',0)$ and $v(j''_{*}M'')=(0,c'',0)$. The sheaf $j'_{*}M'\oplus j''_{*}M''$ is then $H-$semistable of Mukai vector $v=(0,c'+c'',0)$, but it is not $H-$stable, since $j'_{*}M'$ and $j''_{*}M''$ are both $H-$destabilizing.}
\end{esem}

The definition of $v-$generic polarization we have given in the previous section is then not well adapted to Mukai vectors of the form $(0,v_{1},0)$. Because of this, we introduce a different definition of genericity for polarizations than can be found in \cite{Y7} (see Definition 1.4 therein for Mukai vectors $(v_{0},v_{1},v_{2})$ with $v_{0}>0$, and Definition 3.1 if $v_{0}=0$) and that is more suitable for Mukai vectors of the form $(0,v_{1},0)$. 

\begin{defn}
{\rm Let $S$ be a projective K3 surface or an Abelian surface, and $v=(v_{0},v_{1},v_{2})$ a Mukai vector on $S$. A polarization $H$ will be called \textit{general with respect to $v$} in the following cases.
\begin{enumerate}
 \item Case 1: when $v_{0}>0$. In this case $H$ is general with respect to $v$ if for every $\mu_{H}-$semistable sheaf $E$ such that $v(E)=v$ and every $0\neq F\subseteq E$, we have that if $\mu_{H}(E)=\mu_{H}(F)$ then $c_{1}(F)/rk(F)=c_{1}(E)/rk(E)$. 
 \item Case 2: when $v_{0}=0$. In this case $H$ is general with respect to $v$ if for every $H-$semistable sheaf $E$ such that $v(E)=v$ and every $0\neq F\subseteq E$, if $\chi(E)/(c_{1}(E)\cdot H)=\chi(F)/(c_{1}(F)\cdot H)$ then $v(F)\in\mathbb{Q}v$. 
\end{enumerate}}
\end{defn}
 
We first prove that this notion is more general than that of $v-$genericity if $v$ is not of the form $(0,v_{1},0)$.

\begin{lem}
\label{lem:vgengen}
Let $S$ be a projective K3 surface or an Abelian surface and $v$ a Mukai vector on $S$ such that if $\rho(S)\geq 2$ and $v=(0,v_{1},v_{2})$, then $v_{2}\neq 0$.
\begin{enumerate}
 \item If $\rho(S)=1$, the ample generator of $Pic(S)$ is both $v-$generic and general with respect to $v$.
 \item If $\rho(S)\geq 2$, if a polarization is $v-$generic, then it is general with respect to $v$. 
\end{enumerate} 
\end{lem}

\proof This is immediate if $\rho(S)=1$, so we suppose $\rho(S)\geq 2$. If $v=(v_{0},v_{1},v_{2})$ and $v_{0}>0$, this is a consequence of Theorem 4.C.3 of \cite{HL}. So we suppose $v_{0}=0$, and hence $v_{2}\neq 0$, and consider a $v-$generic polarization $H$. 

To show that $H$ is general with respect to $v$, let $E$ be a $H-$semistable sheaf with Mukai vector $v$ and $0\neq F\subseteq E$ a subsheaf with Mukai vector $u=(0,u_{1},u_{2})$. Notice that $\chi(E)=v_{2}$ and $\chi(F)=u_{2}$.

If $\chi(E)/(c_{1}(E)\cdot H)=\chi(F)/(c_{1}(F)\cdot H)$, then $(u_{2}v_{1}-v_{2}u_{1})\cdot H=0$. Since $u_{2}v_{1}-v_{2}u_{1}$ is the divisor associated with $(E,F)$, and since $H$ is $v-$generic, it follows that $u_{2}v_{1}-v_{2}u_{1}=0$, so that $$v(F)=u=(0,u_{1},u_{2})=\frac{u_{2}}{v_{2}}(0,v_{1},v_{2})\in\mathbb{Q}v,$$and hence $H$ is general with respect to $v$.\endproof

If $\rho(S)\geq 2$, a polarization $H$ which is general with respect to $v$ is not necessarily $v-$generic, so it may lie on a $v-$wall. Anyway the moduli space $M_{v}(S,H)$ is equal to a moduli space $M_{v}(S,H')$ where $H'$ is $v-$generic. This is the content of the following:

\begin{lem}
\label{lem:relgenvgen}
Let $S$ be a projective K3 surface or an Abelian surface and $v$ a Mukai vector on $S$. Suppose that $\rho(S)\geq 2$, and that if $v=(0,v_{1},v_{2})$ then $v_{2}\neq 0$. Suppose moreover that $H$ is a polarization which lies on a $v-$wall and is general with respect to $v$, and let $\mathcal{C}$ be a $v-$chamber such that $H\in\overline{\mathcal{C}}$ (the closure of $\mathcal{C}$ in $Amp(S)$). Then a sheaf $E$ of Mukai vector $v$ is $H-$(semi)stable if and only if it is $H'-$(semi)stable for every $H'\in\mathcal{C}$. In particular, we have a identifications $M_{v}(S,H)=M_{v}(S,H')$ and $M^{s}_{v}(S,H)=M^{s}_{v}(S,H')$.
\end{lem}

\proof We present a proof of this for a Mukai vector $v=(v_{0},v_{1},v_{2})$ with $v_{0}>0$, the case $v_{0}=0$ being easier. We will only consider semistable sheaves, the case of stable sheaves being similar.

To do so, suppose that $E$ is $H'-$semistable but not $H-$semistable. Then it is $\mu_{H}-$semistable (since $\mu-$semistability is preserved passing to limits on a wall) and it has a proper subsheaf $F$ such that $p_{H}(F)>p_{H}(E)$: this implies that $\mu_{H}(E)=\mu_{H}(F)$, hence $c_{1}(F)/rk(F)=c_{1}(E)/rk(E)$ (as $H$ is general with respect to $v$). But since $p_{H}(F)>p_{H}(E)$, we get that $p_{H'}(F)>p_{H'}(E)$, which is not possible. As a consequence if $E$ is $H'-$semistable, then it is $H-$semistable.

Conversely, suppose that $E$ is $H-$semistable but not $H'-$semistable. Then $E$ has a proper subsheaf $F$ such that $p_{H'}(F)>p_{H'}(E)$. If $\mu_{H'}(F)>\mu_{H'}(E)$, then we must have $\mu_{H}(F)=\mu_{H}(E)$ (otherwise the segment $[H,H']$ would meet a $v-$wall inside $\mathcal{C}$), so $c_{1}(F)/rk(F)=c_{1}(E)/rk(E)$. Since $p_{H}(F)\leq p_{H}(E)$, we get $p_{H'}(F)\leq p_{H'}(E)$, which is not possible. If $\mu_{H'}(F)=\mu_{H'}(E)$, by $v-$genericity of $H'$ this implies again $c_{1}(F)/rk(F)=c_{1}(E)/rk(E)$, hence we would get $p_{H}(F)>p_{H}(E)$, which is again not possible. As a consequence, if $E$ is $H-$semistable, then it is $H'-$semistable.\endproof

As a consequence of Lemmata \ref{lem:vgengen} and \ref{lem:relgenvgen}, it will be enough to consider $v-$generic polarizations in order to get results for polarizations which are general with respect to $v$, at least in the case where $v=(v_{0},v_{1},v_{2})$ with either $v_{0}>0$ or $v_{0}=0$ and $v_{2}\neq 0$.

The case of $v=(0,v_{1},0)$ was not yet considered, and it is indeed a very special case: while all polarizations are $v-$generic in this case, it may happen that no general polarizations with respect to $v$ exist at all. An example of this is the following:

\begin{esem}
\label{esem:noqvg2}
{\rm We use the same notation of Example \ref{esem:noqvg}, and we show that if $v=(0,c'+c'',0)$, then there is no general polarization with respect to $v$. Indeed, the sheaf $E:=j'_{*}M'\oplus j''_{*}M''$ is $H-$semistable for every polarization $H$, and the Mukai vector of the destabilizing subsheaf  $F:=j'_{*}M'$ is $v(F)=(0,c',0)\notin\mathbb{Q}v$, but $$\chi(F)/c_{1}(F)\cdot H=0=\chi(E)/c_{1}(E)\cdot H.$$This shows that $H$ is not general with respect to $v$.}
\end{esem}
 
If $v=(0,v_{1},0)$ and a general polarization $H$ with respect to $v_1$ exists, the following Lemma implies that $H$ is still general with respect to $v_{H}:=(0,v_{1},v_{1}\cdot H)$ and that $M_{v}(S,H)$ is isomorphic to $M_{v_H}(S,H)$. By Lemma \ref{lem:relgenvgen} it follows that $M_{v}(S,H)$ is also isomorphic to $M_{v_H}(S,H')$ for a $v_H-$generic polarization $H'$.

\begin{lem}
\label{lem:tensorh0}
Let $S$ be a projective K3 surface or an Abelian surface and $v=(0,v_{1},v_{2})$ a Mukai vector on $S$. A polarization $H$ is general with respect to $v$ if and only if it is general with respect to $v_{H}$. Moreover, the tensorization with $H$ induces an isomorphism between $M_{v}(S,H)$ and $M_{v_{H}}(S,H)$, and between $M_{v}^{s}(S,H)$ and $M_{v_{H}}^{s}(S,H)$.
\end{lem}

\proof Since the tensorization with a multiple of $H$ preserves the $H-$(semi)stability the Lemma follows. \endproof

\begin{oss} {\rm As a consequence of Lemmata \ref{lem:vgengen}, \ref{lem:relgenvgen} and \ref{lem:tensorh0}, if $v$ is a Mukai vector and $H$ is a polarization with is general with respect to $v$, then $M_{v}(S,H)$ is either identified with or isomorphic to a moduli space $M_{v'}(S,H')$ where $H'$ is $v'-$generic. This will allow us to restrict ourselves to $v-$generic polarizations in the proof of the main results of the present paper, and to prove all the results needed for the proofs of Theorems \ref{thm:main} and  \ref{thm:mainb1} only for $v-$generic polarizations.}
\end{oss}

A further reason to consider polarizations which are general with respect to a Mukai vector $v$ instead of $v-$generic polarizations is the openness, in the algebraic category, of the locus of the base of a deformation where a polarization stays general with respect to a Mukai vector. This is the content of the following Proposition.

\begin{prop}
	\label{prop:generalopen}
	Let $f:\mathscr{X}\longrightarrow T$ be a smooth projective family of K3 surfaces or Abelian surfaces over a connected algebraic variety $T$. Let $\mathscr{L},\mathscr{H}$ be two line bundles on $\mathscr{X}$, and for every $t\in T$ let $L_{t},H_{t}$ be the restrictions of $\mathscr{L},\mathscr{H}$ to the fiber $X_{t}$ of $f$ over $t$. Let $v_{t}:=(r,c_{1}(L_{t}),a)$ be a Mukai vector on $X_{t}$ and suppose that $H_t$ is ample on  $X_{t}$ for every $t\in T$. Then the locus  $$T_{ng}:=\left\{t\in T\,|\,H_{t} \;\mbox{is not general with respect to}\; v_{t}\right\}\subseteq T$$ is Zariski closed in $T$.
\end{prop}

\proof Let us first consider the case of $r>0$. As an immediate consequence of the definition, the polarization $H_t$ is not general with respect to $v_t$ if and only if there exists a  $\mu_{H_t}-$semistable sheaf $E$ of  Mukai vector $v_t$ on $X_t$ admitting a surjective map to a torsion free  sheaf $G$ having the same $
H_t-$slope of $E$ and such that $c_1(G)/rk(G)\neq c_{1}(E)/rk(E)$.

Since the family of $\mu_{H_t}-$semistable sheaves with Mukai vectors $v_t$ on the fibers of $f$ is  bounded (see as instance Theorem 3.3.7 of \cite{HL}),  there is a variety $Y$ with a morphism  $\pi_{Y}:Y\longrightarrow T$ and a $Y-$flat sheaf $\mathcal{E}_Y$ on $\mathscr{X}_Y:=Y\times_T \mathscr{X}$ such that, for every $t\in T$, every $\mu_{H_t}-$semistable sheaf with Mukai vector $v_t$ appears as the restriction $E_y$ of $\mathcal{E}_{Y}$ to $\{y\}\times X_t$ for some $y\in \pi_{Y}^{-1}(t)\subseteq Y$.

By a result of Grothendieck (see Lemma 1.7.9 of \cite{HL}), there exists a  finite  set $I\subset \mathbb{Q}[x]$ such that  $P\in I$ if and only if  there are $t\in T$, $y\in \pi_{Y}^{-1}(t)$ and a torsion free quotient of $E_y$ whose Hilbert polynomial with respect to $H_{t}$ is $P$ and whose $H_{t}-$slope is $\mu_{H_{t}}(E_y)$. 

For $P\in I$ we let $q_{P,Y}:Q_{P,Y}\rightarrow Y$ be the relative Quot scheme whose fiber over $y\in Y $ parametrizes the quotients of $E_y$ whose Hilbert polynomial with respect to $H_{\pi_{Y}(y)}$ is $P$ and let $Q^0_{P,Y}\subset Q_{P.Y}$ be the locus parametrizing quotients $E_y\rightarrow  G_y$ with $c_1(G_y)/rk(G_y)\neq c_{1}(E_y)/rk(E_y)$. Since  $L_t$ is the restriction of the global line bundle $\mathcal{L}$, the locus $Q^0_{P,Y}$ is a union of connected components of $Q_{P,Y}$ and by construction
$$T_{ng}=\bigcup_{P\in I} \pi_{Y}(q_{P,Y}(Q^0_{P,Y})).$$

As a finite union of images of regular morphisms, we have that $T_{ng}$ is constructible, and we can check its closure in $T$ by the evalutative criterion of properness for the inclusion $T_{ng}\subset T$. 

It then remains to prove that, if $C$ is a smooth curve and $C^{0}:=C\setminus \overline{c}$ is the complement of a point $\overline{c}$ and $\iota:C\rightarrow T$ is a  morphism such that $\iota (C^0)\subset T_{ng}$, then  $\iota (\overline{c})\in T_{ng}$ too.
Using the   universal family of the relative Quot scheme $Q_{P,Y}$ and replacing $C^{0}$ by a quasi-finite cover if necessary, we may assume that there exist two $C^0-$flat families $\mathcal{E}_{C^0}$ and $\mathcal{G}_{C^0}$ of torsion free sheaves on $\mathscr{X}_{C^0}:=C^0\times_T \mathscr{X}$, together with a surjective morphism $\mathcal{E}_{C^0}\rightarrow \mathcal{G}_{C^0}$, such that for every $c\in C^0$ the two following conditions are fulfilled: 
\begin{enumerate} 
	\item the restriction $E_{c}$ of  $\mathcal{E}_{C^0}$ to $X_{\iota(c)}$ has Mukai vector $v_{\iota(c)}$ and is $\mu_{H_{\iota(c)}}-$semistable, 
	\item the restriction $G_c$ of $\mathcal{G}_{C^0}$ to $X_{\iota(c)}$ has the same $\mu_{H_{\iota(c)}}-$slope of $E_c$ but  $c_1(G_c)/rk(G_c)\neq c_{1}(E_c)/rk(E_c)$.
\end{enumerate}

As for every coherent sheaf on $\mathscr{X}_{C^0}$, we can extend $\mathcal{E}_{C^0}$ to  a coherent sheaf $\mathcal{E}_C$ on $\mathscr{X}_{C}:=C\times_T \mathscr{X}$ and, by the relative version of Langton's theorem \cite{langton} (for a proof working in the relative case see Theorem 2.B.1 of \cite{HL}), we may also assume that its restriction $E_{\overline{c}}$ to $X_{\iota(\overline{c})}$ is $\mu_{H_{\iota(\overline{c})}}-$semistable. 

We need to show that there exists a torsion free  quotient   $E_{\overline{c}}\rightarrow G_{\overline{c}}$ such that
$\mu_{H_{\iota(\overline{c})}}(G_{\overline{c}})=\mu_{H_{\iota(\overline{c})}}(E_{\overline{c}})$ and
$c_1(G_{\overline{c}})/rk(G_{\overline{c}})\neq c_{1}(E_{\overline{c}})/rk(E_{\overline{c}})$.

By flatness, the Hilbert polynomial with respect to $H_{\iota(c)}$ does not depend on $c$ and has to be a fixed element $P\in I$ for every $c\in C^0$. 
The relative  Quot scheme 
$Q_{P,C}\rightarrow C$, whose fiber over $c\in C$ parametrizes quotients of $E_c$ having Hilbert polynomial $P$, has a connected component $Q_{P,C}^0$  containing the sheaves $G_c$ verifying the above condition (2) and, hence, dominating $C$:  by projectivity of the relative Quot scheme, the connected component $Q_{P,C}^0$ also surjects on $C$. 

This implies that  there exists a quotient $G_{\overline{c}}$ of $E_{\overline{c}}$ having Hilbert polynomial $P$ and, since $G_{\overline{c}}$ is in the same connected component $Q_{P,C}^0$ of the $G_c$ for $c\neq\overline{c}$, we see that $G_{\overline{c}}$ verifies the above condition (2) as well, so we have  $c_1(G_{\overline{c}})/rk(G_{\overline{c}})\neq c_{1}(E_{\overline{c}})/rk(E_{\overline{c}})$. 

If $G_{\overline{c}}$ is torsion-free, we see that $\iota(\overline{c})\in T_{ng}$, and we are done. If $G_{\overline{c}}$ is not 
torsion-free, since $E_{\overline{c}}$ is slope-semistable, by considering its quotient by the torsion subsheaf we get a torsion-free sheaf with the same rank and first Chern class and we are done.

We are left with the case $r=0$. Let $\phi:\mathscr{M}\longrightarrow T$ be the relative moduli space associated to $f:\mathscr{X}\longrightarrow T$, $\mathscr{H}$ and $\mathscr{L}$, and let $\Sigma\subseteq\mathscr{M}$ be the closed subset parameterizing strictly semistable sheaves. Let $\Sigma'\subseteq\Sigma$ be the closed subset parameterizing polystable sheaves of the form $F_{1}\oplus...\oplus F_{s}$ such that $v(F_{i})\not\in\mathbb{Q} v$ for some $i\le s$. Since $\phi$ is projective, the image $\phi(\Sigma')$ is closed in $T$ and, since $T_{ng}=\phi(\Sigma')$, we are done.\endproof

Proposition \ref{prop:generalopen} shows that being general with respect to a Mukai vector is an open property in the Zariski topology: this is remarkable in particular in view of the fact that the $v-$genericity is only open in the analytic topology (see \cite{PR}, Corollary 4.2).

\subsection{Yoshioka's fibration}

Here we recall the definition and the main properties of the morphism $a_{v}:M_{v}(S,H)\longrightarrow S\times\widehat{S}$ introduced by Yoshioka in \cite{Y8}, and relate it to another morphism used in \cite{OG3}.

We let $S$ be an Abelian surface, $\widehat{S}$ its dual and $\mathcal{P}$ the Poincar\'e line bundle on $S\times\widehat{S}$. Fix a Mukai vector $v$ and a polarization $H$ on $S$ which is general with respect to $v$, and let $M_{v}(S,H)$ be the moduli space of $H-$semistable sheaves on $S$ with Mukai vector $v$.

We recall that for a smooth projective variety $X$, the Grothendieck group $K(X)$ has a ring structure: if $F$ and $G$ are two locally free sheaves, we let $[F]+[G]=[F\oplus G]$ and $[F]\cdot[G]=[F\otimes G]$; if $F$ and $G$ are coherent but not locally free, replace them by a finite locally free resolution of both. 

If $f:X\longrightarrow Y$ is a morphism of smooth projective varieties, then we have the ring morphisms of pull-back $f^{*}:K(Y)\longrightarrow K(X)$ and the group morphism of push-forward $f_{!}:K(X)\longrightarrow K(Y)$. Moreover, the determinant map $\det:K(X)\longrightarrow Pic(X)$ is well-defined. 

\subsubsection{Yoshioka's fibration}

We now define the morphism $a_{v}$ following Yoshioka (see \cite{Y8} and \cite{Y2}). To do so, fix a coherent sheaf $\mathcal{F}_{0}$ with Mukai vector $v$. For every coherent sheaf $\mathcal{F}$ on $S$ of Mukai vector $v$, we set $$\delta_{v}(\mathcal{F}):=\det(p_{\widehat{S}!}(p_{S}^{*}([\mathcal{F}]-[\mathcal{F}_{0}])\cdot([\mathcal{P}]-[\mathcal{O}_{S\times\widehat{S}}]))\in Pic^{0}(\widehat{S}),$$where $p_{S}$ and $p_{\widehat{S}}$ are the two projections of $S\times\widehat{S}$ onto $S$ and $\widehat{S}$ respectively, and $Pic^{0}(\widehat{S})$ is the group of topologically trivial line bundles on $\widehat{S}$. Letting $$F:D^{b}(S)\longrightarrow D^{b}(\widehat{S}),\,\,\,\,\,\,\,\,F(E^{\bullet}):=Rp_{\widehat{S}*}(p_{S}^{*}E^{\bullet}\otimes\mathcal{P})$$be the Fourier-Mukai functor with kernel $\mathcal{P}$, then we have $$\delta_{v}(\mathcal{F})=\det(F(\mathcal{F}))\otimes\det(F(\mathcal{F}_{0}))^{\vee}\in Pic^{0}(\widehat{S}).$$Notice that we have an isomorphism $Pic^{0}(\widehat{S})\simeq S$, hence we have a morphism $\delta_{v}:M_{v}(S,H)\longrightarrow S$.

We then let $$a_{v}:M_{v}(S,H)\longrightarrow S\times\widehat{S},\,\,\,\,\,\,\,\,a_{v}(\mathcal{F}):=(\delta_{v}(\mathcal{F}),\det(\mathcal{F})\otimes\det(\mathcal{F}_{0})^{\vee}).$$Let now $K_{v}(S,H):=a_{v}^{-1}(0_{S},\mathcal{O}_{S})$, where $0_{S}$ is the zero of the Abelian group $S$: if $v^{2}>0$, the morphism $$\tau_{v}:K_{v}(S,H)\times S\times\widehat{S}\longrightarrow M_{v}(S,H),\,\,\,\,\,\,\tau(\mathcal{E},p,L):=\tau_{p}^{*}(\mathcal{E})\otimes L,$$is a finite \'etale cover (for a proof of this, see section 4.2 of \cite{Y2}). We will moreover let $K_{v}^{s}(S,H):=K_{v}(S,H)\cap M^{s}_{v}(S,H)$.

\subsubsection{O'Grady's fibration}

Another morphism $b_{v}:M_{v}(S,H)\longrightarrow S\times\widehat{S}$ was used by O'Grady in \cite{OG3}. For $\gamma\in CH_{0}(S)$, we let $\Sigma(\gamma)\in S$ be the sum of the points of the support of a representative of $\gamma$, counted with multiplicities (i. e. the Albanese image of $\gamma$). For a coherent sheaf $\mathcal{F}$ on $S$ we let $\textbf{c}_{2}(\mathcal{F})\in CH_{0}(S)$ be the second Chern class of $\mathcal{F}$, and we set $$\beta(\mathcal{F}):=\Sigma(\textbf{c}_{2}(\mathcal{F})).$$

The morphism $b_{v}:M_{v}(S,H)\longrightarrow S\times\widehat{S}$ is defined by 
$$b_{v}(\mathcal{F}):=(\beta(\mathcal{F}),\det(\mathcal{F})\otimes\det(\mathcal{F}_{0})^{\vee}).$$

The relation between $a_{v}$ and $b_{v}$ is explained in the following

\begin{lem}
\label{lem:avbv}There is an automorphism $g:S\longrightarrow S$ such that $$b_{v}=(g\times id_{\widehat{S}})\circ a_{v}.$$
\end{lem}

\proof To prove this, we just need to show that for every $\mathcal{F}_{1},\mathcal{F}_{2}\in M_{v}(S,H)$, we have $a_{v}(\mathcal{F}_{1})=a_{v}(\mathcal{F}_{2})$ if and only if $b_{v}(\mathcal{F}_{1})=b_{v}(\mathcal{F}_{2})$. Equivalently, we just need to show that for every $\mathcal{F}_{1},\mathcal{F}_{2}\in M_{v}(S,H)$ such that $\det(\mathcal{F}_{1})\simeq\det(\mathcal{F}_{2})$, we have $\delta_{v}(\mathcal{F}_{1})=\delta_{v}(\mathcal{F}_{2})$ if and only if $\beta(\mathcal{F}_{1})=\beta(\mathcal{F}_{2})$.

Suppose first that $\det(\mathcal{F}_{1})\simeq\det(\mathcal{F}_{2})$ and $\delta_{v}(\mathcal{F}_{1})=\delta_{v}(\mathcal{F}_{2})$, and let $\Gamma:=[\mathcal{F}_{1}]-[\mathcal{F}_{2}]\in K(S)$. As $v(\mathcal{F}_{1})=v(\mathcal{F}_{2})=v$ and $\det(\mathcal{F}_{1})\simeq\det(\mathcal{F}_{2})$, the only nontrivial Chern class of $\Gamma$ (in the Chow ring of $S$) is $\textbf{c}_{2}(\Gamma)$.

Moreover, there is a representative of $\textbf{c}_{2}(\Gamma)$ of the form $$\overline{\Gamma}:=\sum_{i=1}^{n}p_{i}-\sum_{i=1}^{n}q_{i},$$where $p_{1},\cdots,p_{n}$ and $q_{1},\cdots,q_{n}$ are points of $S$. We then notice that $\Gamma\in K(S)$ has the same rank and Chern classes of the class $$\Gamma':=\big[\oplus_{i=1}^{n}\mathbb{C}_{p_{i}}-\oplus_{i=1}^{n}\mathbb{C}_{q_{i}}\big]\in K(S).$$Notice that if we let $\widetilde{F}:K(S)\longrightarrow K(\widehat{S})$ be the morphism induced by $F$ on the level of the Grothendieck groups, we have $$\det[\widetilde{F}(\Gamma')]=\otimes_{i=1}^{n}\mathcal{P}_{p_{i}}\otimes\mathcal{P}_{q_{i}}^{\vee}.$$

As $\det(\widetilde{F}(\Gamma))$ depends only on the rank and the Chern classes of $\Gamma$ in the Chow ring of $\widehat{S}$, we get $$\det(\widetilde{F}([\mathcal{F}_{1}]))\otimes\det(\widetilde{F}([\mathcal{F}_{2}]))^{\vee}=\det(\widetilde{F}(\Gamma))=\det(\widetilde{F}(\Gamma'))=$$ $$=\otimes_{i=1}^{n}\mathcal{P}_{p_{i}}\otimes\mathcal{P}_{q_{i}}^{\vee}=\mathcal{P}_{\sum_{i=1}^{n}(p_{i}-q_{i})}=\mathcal{P}_{\overline{\Gamma}},$$where the equality $\otimes_{i=1}^{n}\mathcal{P}_{p_{i}}\otimes\mathcal{P}_{q_{i}}^{\vee}=\mathcal{P}_{\sum_{i=1}^{n}(p_{i}-q_{i})}$ follows from the fact that the map $$S\longrightarrow\widehat{S},\,\,\,\,\,\,\,\,p\mapsto\mathcal{P}_{p}$$is a group isomorphism.

Now, notice that as $\det(\mathcal{F}_{1})=\det(\mathcal{F}_{2})$, we have that $\delta_{v}(\mathcal{F}_{1})=\delta_{v}(\mathcal{F}_{2})$ if and only if $\det(\widetilde{F}([\mathcal{F}_{1}]))=\det(\widetilde{F}([\mathcal{F}_{2}]))$. The previous equalities give that this holds if and only if $\mathcal{P}_{\overline{\Gamma}}=\mathcal{O}_{S}$. But is equivalent to $\Sigma([\overline{\Gamma}])=0_{S}$, where $[\overline{\Gamma}]$ is the class of $\overline{\Gamma}$ in $CH_{0}(S)$. 

As this class is $\textbf{c}_{2}(\Gamma)$, we finally get that $\delta_{v}(\mathcal{F}_{1})=\delta_{v}(\mathcal{F}_{2})$ if and only if $\Sigma(\textbf{c}_{2}(\Gamma))=0_{S}$. This last condition is equivalent to $\Sigma(\textbf{c}_{2}(\mathcal{F}_{1}))=\Sigma(\textbf{c}_{2}(\mathcal{F}_{2}))$, i. e. to $\beta(\mathcal{F}_{1})=\beta(\mathcal{F}_{2})$, concluding the proof.\endproof

As a consequence, we see that $b_{v}$ is an isotrivial fibration.

\subsubsection{Fibers of the Yoshioka fibration}

The behaviour of moduli spaces of sheaves on an Abelian surface $S$ under changing of polarization in the ample cone of $S$ may now be generalized to the fibers of their Yoshioka (or O'Grady) fibration. 

The main result is the following:

\begin{lem}
\label{lem:polchangeab}
Suppose that $S$ is an Abelian surface with $\rho(S)\geq 2$ and that $v$ is a Mukai vector on $S$. \begin{enumerate}
	\item If $v=(v_{0},v_{1},v_{2})$ is such that $v_{0}>0$ or $v_{0}=0$ and $v_{2}\neq 0$, let $\mathcal{C}$ be a $v-$chamber. 
	\begin{enumerate}
		\item If $H,H'\in\mathcal{C}$, then we have natural identifications $K_{v}(S,H)=K_{v}(S,H')$ and $K^{s}_{v}(S,H)=K^{s}_{v}(S,H')$.
		\item If $H\in\overline{\mathcal{C}}$ is general with respect to $v$ and $H'\in\mathcal{C}$, then we have natural identifications $K_{v}(S,H)=K_{v}(S,H')$ and $K^{s}_{v}(S,H)=K^{s}_{v}(S,H')$.
	\end{enumerate}
	\item If $v=(0,v_{1},0)$ and $H$ is general with respect to $v$, then the tensorization with $H$ induces an isomorphism between $K_{v}(S,H)$ and $K_{v_{H}}(S,H)$, and an isomorphism between $K^{s}_{v}(S,H)$ and $K^{s}_{v_{H}}(S,H)$. 
\end{enumerate}
\end{lem} 

\proof The first item of the statement is an immediate consequence of Lemma \ref{lem:relgenvgen}. For the second, use Lemma \ref{lem:tensorh0} by noticing that by the very definition of $b_{v}$ the tensorization with $H$ maps fibers of $b_{v}$ to fibers of $b_{v_{H}}$ and hence by Lemma \ref{lem:avbv} fibers of $a_{v}$ to fibers of $a_{v_{H}}$.\endproof

\subsection{Deformations of a surface with a Mukai vector and a polarization}

We introduce the main construction we use in what follows. Let $T$ be a smooth, connected algebraic variety, and use the following notation: if $f:Y\longrightarrow T$ is a morphism and $\mathscr{L}\in Pic(Y)$, for every $t\in T$ we let $Y_{t}:=f^{-1}(t)$ and $\mathscr{L}_{t}:=\mathscr{L}_{|Y_{t}}$.

\begin{defn}
{\rm Let $S$ be a projective K3 or Abelian surface, $v$ a Mukai vector on $S$ and $H$ a polarization on $S$. Write $v=m(r,\xi,a)$, where $\xi=c_{1}(L)$. If $T$ is a smooth, connected  algebraic variety, a \textit{deformation of $(S,v,H)$ along $T$} is a triple $(\mathscr{X},\mathscr{L},\mathscr{H})$, where:
\begin{enumerate}
 \item $\mathscr{X}$ is a projective, smooth deformation of $S$ along $T$, i. e. there is a smooth, projective, surjective map $f:\mathscr{X}\longrightarrow T$ such that $\mathscr{X}_{t}$ is a projective surface for every $t\in T$, and there is $0\in T$ such that $\mathscr{X}_{0}\simeq S$;
 \item $\mathscr{L}$ is a line bundle on $\mathscr{X}$ such that $\mathscr{L}_{0}\simeq L$.
 \item $\mathscr{H}$ is a line bundle on $\mathscr{X}$ such that $\mathscr{H}_{t}$ is ample for every $t\in T$ and such that $\mathscr{H}_{0}\simeq H$;
\end{enumerate}
For every $t\in T$ we will write $v_{t}:=m(r,c_{1}(\mathscr{L}_{t}),a)$.}
\end{defn}

\begin{oss}
\label{oss:ap}{\rm Let $S$ be a projective K3 (resp. Abelian) surface, $v=(v_{0},v_{1},v_{2})$ a Mukai vector on $S$. Let $T$ be a smooth, connected variety, and $f:\mathscr{X}\longrightarrow T$ a smooth, projective deformation of $S$ such that $\mathscr{X}_{0}\simeq S$ for some $0\in T$. Suppose that on $\mathscr{X}$ there are two line bundles $\mathscr{H}$ and $\mathscr{L}$, and let $H:=\mathscr{H}_{0}$, $L:=\mathscr{L}_{0}$ and suppose that $H$ is ample and that $c_{1}(L)=v_{1}$. Then $(\mathscr{X},\mathscr{L},\mathscr{H})$ is a deformation of $(S,v,H)$ along $T$ if and only if $\mathscr{H}_{t}$ is ample for every $t\in T$. As the set of $t\in T$ such that $\mathscr{H}_{t}$ is ample is Zariski open in $T$, by restricting to a nonempty Zariski open subset of $T$ we may assume that $(\mathscr{X},\mathscr{L},\mathscr{H})$ is a deformation of $(S,v,H)$ along $T$. Moreover, if we assume that $H$ is general with respect to $v$, thanks to Proposition \ref{prop:generalopen}, by restricting to a smaller nonempty Zariski open subset of $T$, we may assume that $\mathscr{H}_{t}$ is general with respect to $v_{t}$ for every $t\in T$.}
\end{oss}

Let $S$ be a projective K3 (resp. Abelian) surface, $v$ a Mukai vector on $S$ and $H$ a polarization on $S$ that is general with respect to $v$. If $(\mathscr{X},\mathscr{L},\mathscr{H})$ is a deformation of $(S,v,H)$ along a smooth, connected algebraic variety $T$, we let $\phi:\mathscr{M}\longrightarrow T$ be the relative moduli space of semistable sheaves and $\phi^{s}:\mathscr{M}^{s}\longrightarrow T$ the relative moduli space of stable sheaves. This means that for every $t\in T$ we have $\mathscr{M}_{t}=M_{v_{t}}(\mathscr{X}_{t},\mathscr{H}_{t})$ and $\mathscr{M}^{s}_{t}=M^{s}_{v_{t}}(\mathscr{X}_{t},\mathscr{H}_{t})$.

If $S$ is Abelian, let $\widehat{\mathscr{X}}\longrightarrow T$ be the dual family, i. e. the connected component of the relative Picard variety $g:Pic_{\mathscr{X}/T}\longrightarrow T$ containing the section of $g$ corresponding to the family $\mathscr{O}_{\mathscr{X}}$. The dual family is then the smooth projective family whose fiber over $t\in T$ is the dual of $\mathscr{X}_{t}$. Consider the following condition:
\begin{center}
($\star$) the morphism $\phi:\mathscr{M}\longrightarrow T$ has a section, and $\mathscr{X}\longrightarrow T$ is a $T-$group scheme.
\end{center}
If the condition ($\star$) holds, we have a $T-$morphism $a_{v}:\mathscr{M}\longrightarrow\mathscr{X}\times_{T}\widehat{\mathscr{X}}$ such that for every $t\in T$ the restriction morphism $a_{v|\mathscr{M}_{t}}$ is the Yoshioka fibration defined in section 2.2. If $$Z:=\{(0_{\mathscr{X}_{t}},\mathcal{O}_{\mathscr{X}_{t}})\in\mathscr{X}_{t}\times\widehat{\mathscr{X}}_{t}\,|\,t\in T\}\subseteq \mathscr{X}\times_{T}\widehat{\mathscr{X}},$$we will let $\mathscr{K}:=a_{v}^{-1}(Z)$: restricting the morphism $\phi$ to $\mathscr{K}$ we get a morphism $\phi_{0}:\mathscr{K}\longrightarrow T$, whose fiber over $t\in T$ is $\mathscr{K}_{t}=K_{v_{t}}(\mathscr{X}_{t},\mathscr{H}_{t})$. A similar definition, but using $\mathscr{M}^{s}$ instead of $\mathscr{M}$, gives the family $\phi_{0}^{s}:\mathscr{K}^{s}\longrightarrow T$.

\begin{oss}
{\rm Condition ($\star$) is always verified up to shrinking $T$ and taking a finite \'etale cover of $T$ thanks to the smoothness of $\mathscr{M}^{s}$ over $T$.}
\end{oss}

The first result we need is that the families $\mathscr{M}$ and $\mathscr{K}$ are $T-$flat over a Zariski open neighborhood of any $t\in T$ such that $\mathscr{H}_{t}$ is general with respect to $v_{t}$. This is the content of the following Lemma:

\begin{lem}
\label{lem:redflat}
Let $S$ be a projective K3 (resp. Abelian) surface, $v$ a Mukai vector on $S$ and $H$ a polarization on $S$ that is general with respect to $v$. Let $T$ be a smooth, connected algebraic  variety, $(\mathscr{X},\mathscr{L},\mathscr{H})$ a deformation of $(S,v,H)$ along $T$, and assume that condition ($\star$) holds if $S$ is Abelian. Suppose that $t\in T$ is such that $\mathscr{H}_{t}$ is general with respect to $v_{t}$. 
\begin{enumerate}
 \item The morphisms $\phi:\mathscr{M}\longrightarrow T$ and $\phi_{0}:\mathscr{K}\longrightarrow T$ are flat at $t$.
 \item The morphisms $\phi^{s}:\mathscr{M}^{s}\longrightarrow T$ and $\phi_{0}^{s}:\mathscr{K}^{s}\longrightarrow T$ are smooth at $t$.
\end{enumerate}
\end{lem}

\proof By Remark \ref{oss:ap} we may suppose that for every $t\in T$ the polarization $\mathscr{H}_{t}$ is general with respect to $v_{t}$. 

(1) Notice that $\mathscr{M}$ (resp. $\mathscr{K}$) is connected (since $T$ and the fibers are connected). Moreover, by Theorem 4.4 of \cite{KLS} (resp. by Remark A.1 of \cite{PR2} for $\mathscr{K}$), and using Lemmata \ref{lem:relgenvgen} and \ref{lem:tensorh0} (risp. Lemma \ref{lem:polchangeab}), we see that the fibers of $\phi$ (resp. $\phi_{0}$) are reduced, irreducible and equidimensional. Now, by Theorem 14.4.4 of \cite{EGA} it follows that $\phi$ (resp. $\phi_{0}$) is universally open. Using Corollary 15.2.3 of \cite{EGA} we get that $\phi$ (resp. $\phi_{0}$) is flat. 

(2) follows from (1) since   $\phi^{s}$ and $\phi_{0}^{s}$ have smooth fibers.
\endproof

Let $S$ be a projective K3 (resp. Abelian) surface, $v$ a Mukai vector on $S$ and $H$ a polarization on $S$ that is general with respect to $v$. By choosing a nontrivial deformation of $(S,v,H)$ along a smooth, connected variety $T$ we get a flat, projective deformation $\phi:\mathscr{M}\longrightarrow T$ of $M_{v}$, and a smooth quasi-projective deformation $\phi^{s}:\mathscr{M}^{s}\longrightarrow T$ of $M^{s}_{v}$. Moreover, if $S$ is Abelian we get a flat, projective deformation $\phi_{0}:\mathscr{K}\longrightarrow T$ of $K_{v}$, and a smooth quasi-projective deformation $\phi_{0}^{s}:\mathscr{K}^{s}\longrightarrow T$ of $K^{s}_{v}$. We now prove that this deformation is locally trivial:

\begin{lem}
\label{lem:loctriv}
Let $S$ be a projective K3 (resp. Abelian) surface, $v$ a Mukai vector on $S$ and $H$ a polarization on $S$ that is general with respect to $v$. Let $T$ be a smooth connected variety, $(\mathscr{X},\mathscr{L},\mathscr{H})$ a deformation of $(S,v,H)$ along $T$, and assume that condition ($\star$) holds if $S$ is Abelian.
\begin{enumerate}
 \item If $p\in\mathscr{M}$ and $t:=\phi(p)$ is such that $\mathscr{H}_{t}$ is general with respect to $v_{t}$, then $(\mathscr{M},p)\simeq(\mathscr{M}_{t},p)\times(T,t)$ as germs of analytic spaces.
 \item If $p\in\mathscr{K}$ and $t:=\phi_{0}(p)$ is such that $\mathscr{H}_{t}$ is general with respect to $v_{t}$, then $(\mathscr{K},p)\simeq(\mathscr{K}_{t},p)\times(T,t)$ as germs of analytic spaces.
\end{enumerate}
\end{lem}

\proof If $v$ is primitive, then $\phi$ is a smooth, projective morphism, and there is nothing to prove. If $v=2w$ where $w$ is primitive and $w^{2}=2$, this is Proposition 2.16 of \cite{PR} if $T$ is a curve (the proof given in \cite{PR} works for polarizations which are general with respect to the Mukai vector). If $T$ is a smooth connected variety of dimension $d\geq 2$, this implies that the statement holds along any smooth connected curve through $t$. By Corollary 0.2 of \cite{FK} the general statement follows. 

For the remaining cases, by \cite{KLS}, $\mathscr{M}_{t}=M_{v_{t}}(\mathscr{X}_{t},\mathscr{H}_{t})$ and $\mathscr{K}_{t}=K_{v_{t}}(\mathscr{X}_{t},\mathscr{H}_{t})$ are symplectic varieties which are locally factorial, and by Corollary 1 of \cite{N1} they have terminal singularities. The Main Theorem of \cite{N2} tells us that for every $p\in \mathscr{M}_{t}$ (resp. $p\in \mathscr{K}_{t}$) and for every $n\in\mathbb{N}$, the infinitesimal $n-$th order deformation of $\mathscr{M}_{t}$ (resp. of $\mathscr{K}_{t}$) induced by $\phi$ (resp. by $\phi_{0}$), which is flat by Lemma \ref{lem:redflat}, is locally trivial at $p$: the statement follows again by Corollary 0.2 of \cite{FK}.\endproof

As a corollary of this, using the Thom First Isotopy Lemma (see Theorem 3.5 of \cite{Dim}) we have the following:

\begin{lem}
\label{lem:toploctriv}
Let $S$ be a projective K3 (resp. Abelian) surface, $v$ a Mukai vector on $S$ and $H$ a polarization on $S$ that is general with respect to $v$. Let $T$ be a smooth connected algebraic variety, let $(\mathscr{X},\mathscr{L},\mathscr{H})$ be  a deformation of $(S,v,H)$ along $T$, and assume that condition ($\star$) holds if $S$ is Abelian.
\begin{enumerate}
 \item If $p\in\mathscr{M}$ and $t:=\phi(p)$ is such that $\mathscr{H}_{t}$ is general with respect to $v_{t}$, there is an analytic open neighborhood $U\subseteq T$ of $t$ such that $\phi^{-1}(U)$ is homeomorphic over $U$ to $\mathscr{M}_{t}\times U$ and $(\phi^{s})^{-1}(U)$ is homeomorphic over $U$ to $\mathscr{M}^{s}_{t}\times U$.
 \item If $p\in\mathscr{K}$ and $t:=\phi_{0}(p)$ is such that $\mathscr{H}_{t}$ is general with respect to $v_{t}$, there is an analytic open neighborhood $U\subseteq T$ of $t$ such that $\phi_{0}^{-1}(U)$ is homeomorphic over $U$ to $\mathscr{K}_{t}\times U$ and  $(\phi_{0}^{s})^{-1}(U)$ is homeomorphic over $U$ to $\mathscr{K}^{s}_{t}\times U$.
\end{enumerate}
\end{lem}

\subsection{Isomorphisms between moduli spaces}

We now describe several isomorphisms that will be frequently used in the proof of Theorem \ref{thm:main}. All of them are induced by Fourier-Mukai transforms, either the tensorization with a line bundle or the one whose kernel is the ideal sheaf of the diagonal (for K3 surfaces) or the Poincar\'e bundle (for Abelian surfaces).

\subsubsection{Isomorphisms from tensorization with line bundles.}

Let $S$ be a projective K3 or Abelian surface, and let $v=m(r,\xi,a)$ be a Mukai vector. Recall that if $L\in Pic(S)$, we defined $v_{L}:=v\cdot ch(L)$, and that if $D$ is a divisor on $S$, we let $v_{D}:=v_{\mathcal{O}_{S}(D)}$ (see section 2.1, Lemma \ref{lem:tensogen}).

\begin{defn}
{\rm Let $v,v'\in\widetilde{H}(S,\mathbb{Z})$ be two Mukai vectors, and let $v=(r,\xi,a)$, $v'=(r',\xi',a')$.
\begin{enumerate}
 \item If $H$ is a polarization on $S$, we say that $v$ and $v'$ are $H-$\textit{equivalent} if there is $s\in\mathbb{Z}$ such that $v'=v_{sH}$.
 \item If $r,r'>0$, we say that $v$ and $v'$ are \textit{equivalent} if there is $L\in Pic(S)$ such that $v'=v_{L}$.
\end{enumerate}}
\end{defn}

The following is the main result about isomorphisms induced by tensorization with a line bundle\footnote{In the following Lemma, by a slight abuse of notation we let $K_{v}$ denote not only the fiber of $a_{v}:M_{v}\longrightarrow S\times\widehat{S}$ over $(0_{S},\mathcal{O}_{S})$, but also any other fiber. This is justified since $a_{v}$ is an isotrivial fibration, so all its fibers are isomorphic.}, which shows that moduli spaces of sheaves corresponding to equivalent (or $H-$equivalent) Mukai vectors are isomorphic (and the isomorphism is induced by tensorization with a suitable line bundle).

\begin{lem}
\label{lem:tensorlb}
Let $S$ be a projective K3 or Abelian surface, $v$ a Mukai vector and $H$ an ample line bundle on $S$.
\begin{enumerate}
 \item For every $d\in\mathbb{Z}$ the morphism $$M_{v}(S,H)\longrightarrow M_{v_{dH}}(S,H),\,\,\,\,\,\,\,\,\mathcal{F}\mapsto\mathcal{F}\otimes\mathcal{O}_{S}(dH)$$is an isomorphism, which induces isomorphisms $M^{s}_{v}(S,H)\simeq M^{s}_{v_{dH}}(S,H)$. If $S$ is an Abelian surface it also induces isomorphisms $K_{v}(S,H)\simeq K_{v_{dH}}(S,H)$, and $K^{s}_{v}(S,H)\simeq K^{s}_{v_{dH}}(S,H)$.
 \item If $v=(v_{0},v_{1},v_{2})$ and $v_{0}>0$, $L\in Pic(S)$ and $H$ is $v-$generic, the morphism $$M_{v}(S,H)\longrightarrow M_{v_{L}}(S,H),\,\,\,\,\,\,\,\mathcal{F}\mapsto\mathcal{F}\otimes L$$is an isomorphism, which induces isomorphisms $M^{s}_{v}(S,H)\simeq M^{s}_{v_{L}}(S,H)$. If $S$ is an  Abelian surface it also induces isomorphisms $K_{v}(S,H)\simeq K_{v_{L}}(S,H)$ and $K^{s}_{v}(S,H)\simeq K^{s}_{v_{L}}(S,H)$. 
\end{enumerate}
\end{lem}

\proof First, notice that $v(\mathcal{F}\otimes L)=v(\mathcal{F})\cdot ch(L)$. To prove the first point of the statement, it is enough to remark that a sheaf $\mathcal{F}$ of Mukai vector $v$ is $H-$(semi)stable if and only if $\mathcal{F}\otimes\mathcal{O}_{S}(dH)$ is $H-$(semi)stable.

For the second point, we need to show that if $\mathcal{F}$ is $H-$(semi)stable, then $\mathcal{F}\otimes L$ is $H-$(semi)stable. This is proved for stable sheaves by Yoshioka (see Lemma 1.1 of \cite{Y2}), and the proof goes through for semistable sheaves.

If $S$ is Abelian, by definition of the morphism $a_{v}$ (see section 2.2) we have that if $\mathcal{F}_{1}$ and $\mathcal{F}_{2}$ are in the same fiber of $a_{v}$, then $\mathcal{F}_{1}\otimes L$ and $\mathcal{F}_{2}\otimes L$ are in the same fiber of $a_{v_{L}}$. As $a_{v}$ and $a_{v_{L}}$ are both isotrivial fibrations, the isomorphism between $M_{v}$ and $M_{v_{L}}$ obtained by tensorization with $L$ induces an isomorphism between $K_{v}$ and $K_{v_{L}}$.\endproof

\subsubsection{Isomorphisms from Fourier-Mukai transforms.}

We now recall two basic results, originally due to Yoshioka, about isomorphisms between moduli spaces of sheaves over K3 or Abelian surfaces coming from Fourier-Mukai transforms. Yoshioka's theorems are stated in a more general setting: here we present simplified adapted proofs for the convenience of the reader. We will only consider the Fourier-Mukai transform whose kernel is the ideal of the diagonal (for K3 surfaces) or the Poincar\'e bundle (for Abelian surfaces).

We need the following notation: if $S$ a projective K3, we let $\Delta\subseteq S\times S$ be the diagonal and $\mathcal{I}$ the ideal of $\Delta$. We have an exact sequence of coherent sheaves on $S\times S$:
\begin{equation}
\label{eq:exideal}
0\longrightarrow\mathcal{I}\longrightarrow\mathcal{O}_{S\times S}\longrightarrow\mathcal{O}_{\Delta}\longrightarrow 0
\end{equation}
We moreover let $\pi_{1},\pi_{2}:S\times S\longrightarrow S$ be the two projections.

If $S$ is an Abelian surface and $\widehat{S}$ is its dual, we let $\mathcal{P}$ be the Poincar\'e line bundle on $S\times\widehat{S}$, $\pi_{1}:S\times\widehat{S}\longrightarrow S$ and $\pi_{2}:S\times\widehat{S}\longrightarrow\widehat{S}$ the two projections, and $\iota:S\longrightarrow S$ the involution acting as $-1$.

We will moreover consider the functors  $$F_{K3}:D^{b}(S)\longrightarrow D^{b}(S),\,\,\,\,\,\,\,\,F_{K3}(E^{\bullet}):=R\pi_{2*}(\pi_{1}^{*}E^{\bullet}\otimes^{L}\mathcal{I}),$$ $$\widehat{F}_{K3}:D^{b}(S)\longrightarrow D^{b}(S),\,\,\,\,\,\,\,\,\widehat{F}_{K3}(E^{\bullet}):=R\mathcal{H}om_{\pi_{1}}(\mathcal{I},\pi_{2}^{*}E^{\bullet})$$if $S$ is a K3 surface, and  $$F_{Ab}:D^{b}(S)\longrightarrow D^{b}(\widehat{S}),\,\,\,\,\,\,\,\,F_{Ab}(E^{\bullet}):=R\pi_{2*}(\pi_{1}^{*}E^{\bullet}\otimes\mathcal{P}),$$ $$\widehat{F}_{Ab}:D^{b}(\widehat{S})\longrightarrow D^{b}(S),\,\,\,\,\,\,\,\,\widehat{F}_{Ab}(E^{\bullet}):=\iota^{*}R\pi_{1*}(\pi_{2}^{*}E^{\bullet}\otimes\mathcal{P})$$if $S$ is an Abelian surface. 

By \cite{Br} we know that $F_{K3}$ and $F_{Ab}$ are equivalences of triangulated categories. Moreover, the functor $\widehat{F}_{K3}[2]$ is the right and left adjoint to $F_{K3}$, so that $F_{K3}\circ\widehat{F}_{K3}=[-2]$ (see Proposition 1.26 of \cite{Hub}), and $\widehat{F}_{Ab}[2]$ is the right and left adjoint to $F_{Ab}$, so that $F_{Ab}\circ\widehat{F}_{Ab}=[-2]$ (see Theorem 2.2 of \cite{M3}).

We will make use of the following definition due to Mukai:

\begin{defn}
{\rm Let $S$ be a projective K3 or Abelian surface, $\mathcal{G}$ a coherent sheaf on $S$ and $F$ a Fourier-Mukai functor on $D^{b}(S)$. For $i\in\{0,1,2\}$, we say that $\mathcal{G}$ \textit{verifies WIT($i$) with respect to $F$} if $F(\mathcal{G})=F^{i}(\mathcal{G})[-i]$.}
\end{defn}

If $S$ is a projective $K3$ surface and $\mathcal{G}$ is a coherent sheaf on $S$, then the functor $R\pi_{2*}(\pi_{1}^{*}\mathcal{G}\otimes^{L}\cdot)$ applied to the exact sequence (\ref{eq:exideal}) gives the long exact sequence of coherent sheaves on $S$
\begin{equation}
\label{eq:exf}
\begin{split}
0\longrightarrow F_{K3}^{0}(\mathcal{G})\longrightarrow\mathcal{O}_{S}\otimes H^{0}(\mathcal{G})\stackrel{ev}\longrightarrow\mathcal{G}\longrightarrow\\
\longrightarrow F_{K3}^{1}(\mathcal{G})\longrightarrow\mathcal{O}_{S}\otimes H^{1}(\mathcal{G})\longrightarrow 0\longrightarrow\\
\longrightarrow F_{K3}^{2}(\mathcal{G})\longrightarrow\mathcal{O}_{S}\otimes H^{2}(\mathcal{G})\longrightarrow 0
\end{split}
\end{equation}
and if $\mathcal{G}$ is torsion free the functor $R\mathcal{H}om_{\pi_{1}}(\cdot,\pi_{2}^{*}\mathcal{G})$ applied to the exact sequence (\ref{eq:exideal}) gives the long exact sequence of coherent sheaves on $S$
\begin{equation}
\label{eq:exfhat}
\begin{split}
0\longrightarrow\mathcal{O}_{S}\otimes H^{0}(\mathcal{G})\longrightarrow\widehat{F}_{K3}^{0}(\mathcal{G})\longrightarrow\\
\longrightarrow 0\longrightarrow\mathcal{O}_{S}\otimes H^{1}(\mathcal{G})\longrightarrow\widehat{F}_{K3}^{1}(\mathcal{G})\longrightarrow\\
\longrightarrow \mathcal{G}\longrightarrow\mathcal{O}_{S}\otimes H^{2}(\mathcal{G})\longrightarrow\widehat{F}_{K3}^{2}(\mathcal{G})\longrightarrow 0
\end{split}
\end{equation}

From now on we will make use of the following notation:

\begin{notation}
\label{not:vettori}
{\rm If $S$ is a projective K3 surface and $v=(r,\xi,a)$ is a Mukai vector on $S$, we let $\widetilde{v}:=(a,-\xi,r)$.}
	
{\rm If $S$ is an Abelian surface and $L\in Pic(S)$, we let $\widehat{L}:=\det(F(L))^{-1}\in Pic(\widehat{S})$; moreover, if $\xi=c_{1}(L)$, we let $\widehat{\xi}:=c_{1}(\widehat{L})$. Finally if $v=(r,\xi,a)$ is a Mukai vector on the Abelian surface $S$, we let $\widetilde{v}:=(a,-\widehat{\xi},r)$ (which belongs to the  Mukai lattice of $\widehat{S}$).}
\end{notation}

\begin{oss}
\label{oss:genfourier}
{\rm If $H$ is an ample line bundle on an Abelian surface $S$, then $\widehat{H}$ is ample (see Proposition 3.11 of \cite{M3}).}
\end{oss}

The first result we need is the following (see Theorem 3.18 of \cite{Y3}).

\begin{lem}
\label{lem:step1}Let $S$ be a projective K3 (resp. Abelian) surface, $v=(v_{0},v_{1},v_{2})$ a Mukai vector on $S$, $H$ a polarization on $S$ and let $h:=c_{1}(H)$. Let $r,k\in\mathbb{N}^{*}$ and $\xi\in NS(S)$ effective, and suppose furthermore that we are in one of the following situations:
\begin{enumerate}
	\item we have $NS(S)=\mathbb{Z}h$, $v^{2}=2k$ and $v_{0}=r$, i. e. there are $p,a\in\mathbb{Z}$ such that $v=(r,ph,a)$ and $v^{2}=2k$, or
	\item we have $v_{0}=0$, $v_{1}=\xi$ and $v^{2}=2k$, i. e. there is $p\in\mathbb{Z}$ such that $v=(0,\xi,p)$ and $v^{2}=2k$.
\end{enumerate}
Then there is $p_{0}\in\mathbb{N}$ such that if $p>p_{0}$ every $H-$semistable sheaf $\mathcal{E}$ with Mukai vector $v$ on $S$ verifies WIT(0) with respect to $F_{K3}$ (resp. $F_{Ab}$), and $F^{0}_{K3}(\mathcal{E})$ (resp. $F^{0}_{Ab}(\mathcal{E})$) is locally free of Mukai vector $\widetilde{v}$.
\end{lem}

\proof We let $$V_{r,k}:=\{v'=(r,p'h,a')\in\widetilde{H}(S,\mathbb{Z})\,|\,(v')^{2}=2k,\,\,0\leq p'\leq r\}$$and $$V_{0,\xi,k}:=\{v'=(0,\xi,p')\,|\,(v')^{2}=2k,\,\,0\leq p'\leq\xi\cdot H\}.$$

Notice that $V_{r,k}$ and $V_{0,\xi,k}$ are finite sets. As the family of semistable sheaves with fixed Mukai vector $v$ is bounded, by Serre Theorem there is $T\in\mathbb{N}$, 
such that for every $s>T$ and for every $H-$semistable sheaf $\mathcal{E}'$ with Mukai vector in $V_{r,k}$ or $V_{0,\xi,k}$ we have $H^{1}(\mathcal{E}'\otimes\mathcal{O}_{S}(sH))=H^{2}(\mathcal{E}'\otimes \mathcal{O}_{S}(sH))=0$ and the evaluation morphism $$H^{0}(\mathcal{E}'\otimes\mathcal{O}_{S}(sH))\otimes\mathcal{O}_{S}\longrightarrow\mathcal{E}'\otimes\mathcal{O}_{S}(sH)$$ is surjective. 

Under the hypotheses of case (1), there are $s\in\mathbb{N}$ and $v'\in V_{r,k}$ such that $v'_{sH}=v=(r,ph,a)$, and notice that if $p>p_0:=r+rT$ then we have $s>T$. Since tensorization by $sH$ induces an isomorphism 
from $M_{v'}(S,H)$ to $M_v(S,H)$ by Lemma \ref{lem:tensorlb}, we conclude that for every $H-$semistable sheaf $\mathcal{E}$ with $v(\mathcal{E})=v$ there is $\mathcal{E}'\in M_{v'}(S,H)$ such that $\mathcal{E}\simeq\mathcal{E}'\otimes\mathscr{O}_{S}(sH)$, and hence $H^{1}(\mathcal{E})=H^{2}(\mathcal{E})=0$ and the evaluation morphism $$H^{0}(\mathcal{E})\otimes\mathcal{O}_{S}\longrightarrow\mathcal{E}$$ is surjective.

Similarily, under the hypotheses of case (2) there are $s\in\mathbb{N}$ and $v'\in V_{0,\xi,k}$ such that $v'_{sH}=v=(0,\xi,p)$, and if $p>p_0:=\xi\cdot H +(\xi\cdot H)T$, by the same argument, the same conclusion holds
for every $H-$semistable sheaf $\mathcal{E}$ with $v(\mathcal{E})=v$.

Now, if $S$ is K3, as $H^{1}(\mathcal{E})=H^{2}(\mathcal{E})=0$ and the evaluation morphism $H^{0}(\mathcal{E})\otimes\mathcal{O}_{S}\longrightarrow\mathcal{E}$ is surjective, the exact sequence (\ref{eq:exf}) for $\mathcal{E}$ implies that $F_{K3}^{1}(\mathcal{E})=F_{K3}^{2}(\mathcal{E})=0$, so that $F_{K3}(\mathcal{E})=F_{K3}^{0}(\mathcal{E})$, which proves that $\mathcal{E}$ verifies WIT(0) with respect to $F_{K3}$.

If $S$ is Abelian, we do not only have that $H^{1}(\mathcal{E})=H^{2}(\mathcal{E})=0$ and the evaluation morphism $H^{0}(\mathcal{E})\otimes\mathcal{O}_{S}\longrightarrow\mathcal{E}$ is surjective, but the same holds for $\mathcal{E}\otimes L$ for every $L\in\widehat{S}$ (since $v(\mathcal{E}\otimes L)=v$). By cohomology and base change it follows that $F_{Ab}^{1}(\mathcal{E})=F_{Ab}^{2}(\mathcal{E})=0$, so that $F_{Ab}(\mathcal{E})=F_{Ab}^{0}(\mathcal{E})$, which proves that $\mathcal{E}$ verifies WIT(0) with respect to $F_{Ab}$.

We are left with showing that $F_{K3}^{0}(\mathcal{E})$ (resp. $F^{0}_{Ab}(\mathcal{E})$) is locally free and its Mukai vector is $\widetilde{v}$. 

To show this, let us first consider $S$ to be a K3 surface. As $\mathcal{E}$ verifies WIT(0) with respect to $F_{K3}$, the exact sequence (\ref{eq:exf}) applied to $\mathcal{E}$ gives the exact sequence 
\begin{equation}
	\label{eq:exf0}
	0\longrightarrow F_{K3}^{0}(\mathcal{E})\longrightarrow H^{0}(\mathcal{E})\otimes\mathcal{O}_{S}\longrightarrow\mathcal{E}\longrightarrow 0.
\end{equation}
We then see that $v(F_{K3}^{0}(\mathcal{E}))=\widetilde{v}$.

As $\mathcal{E}$ is a coherent sheaf of pure dimension $2$ or $1$ on the smooth surface  $S$, the projective dimension of $\mathcal{E}$ is at most 1:
since the sequence (\ref{eq:exf0}) is exact and  $H^{0}(\mathcal{E})\otimes\mathcal{O}_{S}$ is locally free, this implies that  $F_{K3}^{0}(\mathcal{E})$ is locally free too.

This completes the proof when $S$ is a K3 surface. The case of Abelian surfaces is easier: the fact that $v(F^{0}_{Ab}(\mathcal{E}))=\widetilde{v}$ is well-known. Moreover, as $H^{i}(\mathcal{E}\otimes L)=0$ for every $L\in\widehat{S}$ and for $i=1,2$, we see that $\mathcal{E}$ is an IT-sheaf of index 0, and hence $F^{0}_{Ab}(\mathcal{E})$ is locally free (see as instance Lemma 14.2.1 of \cite{BL}).\endproof

We are now in the position to prove the main results of this section. The first is the following, showing that if $S$ is a K3 (resp. Abelian) surface with Picard number 1 and $v$ is a Mukai vector on $S$ with positive rank, then $F_{K3}$ (resp. $F_{Ab}$) induces an isomorphism between $M_{v}$ and $M_{\widetilde{v}}$ (see Notation \ref{not:vettori} for the definition of $\widetilde{v}$).

\begin{prop}
\label{prop:yoshi2}Let $S$ be a K3 or Abelian surface such that $NS(S)=\mathbb{Z}\cdot h$, and let $h=c_{1}(H)$ for an an ample line bundle $H$. Let $r,k\in\mathbb{N}^{*}$ and let $v=(v_0,v_1,v_2)$ be a Mukai vector on $S$ such that $v_0=r$ and $v^{2}=2k$, i.e. there are $n,a\in\mathbb{Z}$ such that $v=(r,nh,a)$ and $v^2=2k$.
\begin{enumerate}
	\item If $S$ is K3, there is $n_{0}\in\mathbb{N}$ such that for every $n>n_{0}$ the functor $F_{K3}$ induces isomorphisms $M_{v}(S,H)\simeq M_{\widetilde{v}}(S,H)$ and $M^{s}_{v}(S,H)\simeq M^{s}_{\widetilde{v}}(S,H)$.
	\item If $S$ is Abelian, there is $n_{0}\in\mathbb{N}$ such that for every $n>n_{0}$ the functor $F_{Ab}$ induces isomorphisms $M_{v}(S,H)\simeq M_{\widetilde{v}}(\widehat{S},\widehat{H})$, $M^{s}_{v}(S,H)\simeq M^{s}_{\widetilde{v}}(\widehat{S},\widehat{H})$, $K_{v}(S,H)\simeq K_{\widetilde{v}}(\widehat{S},\widehat{H})$ and $K^{s}_{v}(S,H)\simeq K^{s}_{\widetilde{v}}(\widehat{S},\widehat{H})$.
\end{enumerate} 
\end{prop}

\proof The surfaces involved in the statements and in the proofs of this Proposition and of the related Lemmas (i. e. a K3 surface $S$ or an Abelian surface $S$ and its dual $\widehat{S}$) have all cyclic N\'eron-Severi group. In particular, if $\Sigma$ is such a surface, then $NS(\Sigma)=\mathbb{Z}\cdot\ell$ where $\ell=c_{1}(L)$ and $L$ is an ample generator.

This allows us to unify and simplify the notation for the Mukai vectors on such a surface $\Sigma$: the Mukai vector $v=(r,n\ell,a)$ will always be written under the simpler form $(r,n,a)$. 

Similarly, while discussing (semi)stability, we systematically avoid any explicit reference to the (essentially unique) polarization. In particular, if $\mathcal{E}$ is a coherent sheaf on $\Sigma$, we will simply write $p(\mathcal{E})$ for its reduced Hilbert polinomial with respect to the primitive ample divisor.

Finally, we will use the notation $F$ for both $F_{K3}$ and $F_{Ab}$, and $\widehat{F}$ for both $\widehat{F}_{K3}$ and $\widehat{F}_{Ab}$. Making use of this simplified notation, we start the proof of the Proposition.

We first notice that by Lemma \ref{lem:step1} there is $n_{0}\in\mathbb{N}$ such that for every $n>n_{0}$ and for every semistable sheaf $\mathcal{E}$ of Mukai vector $v=(r,n,a)$ on $S$ we have that $F(\mathcal{E})=F^{0}(\mathcal{E})$ is a locally free sheaf with Mukai vector $\widetilde{v}$.

Our aim is to prove that the locally free sheaf $F^{0}(\mathcal{E})$ is semistable. Once this will be done, it will imply that $F$ induces an injective morphism $f_{K3}:M_{v}(S,H)\longrightarrow M_{\widetilde{v}}(S,H)$ if $S$ is $K3$ and $f_{Ab}:M_{v}(S,H)\longrightarrow M_{\widetilde{v}}(\widehat{S},\widehat{H})$ if $S$ is Abelian. By Theorem 4.4 of \cite{KLS} these moduli spaces are irreducible of the same dimension, so $f_{K3}$ and $f_{Ab}$ are isomorphisms and induce  isomorphisms between the smooth loci of the moduli spaces.

Moreover, if $S$ is Abelian and $\mathcal{E},\mathcal{E}_{0}\in M_{v}(S,H)$, as $\iota^{*}\circ\widehat{F}:D^{b}(\widehat{S})\longrightarrow D^{b}(S)$ is the Fourier-Mukai transform with kernel $\mathcal{P}$, by definition of $a_{\widetilde{v}}$ we have $$a_{\widetilde{v}}(F(\mathcal{E}))=$$ $$=(\det(\iota^{*}\widehat{F}(F(\mathcal{E})))\otimes\det(\iota^{*}\widehat{F}(F(\mathcal{E}_{0})))^{\vee},\det(F(\mathcal{E}))\otimes\det(F(\mathcal{E}_{0}))^{\vee})=$$ $$=(\iota^{*}(\det(\mathcal{E})\otimes\det(\mathcal{E}_{0})^{\vee}),\det(F(\mathcal{E}))\otimes\det(F(\mathcal{E}_{0}))^{\vee})=\varepsilon(a_{v}(\mathcal{E})),$$where $$\varepsilon:S\times\widehat{S}\longrightarrow\widehat{S}\times S,\,\,\,\,\,\,\,\,\,\varepsilon(p,q):=(\widehat{\iota}(q),p),$$and $\widehat{\iota}:\widehat{S}\longrightarrow\widehat{S}$ is the involution acting as $-1$.

It follows that $\mathcal{E}_{1}$ and $\mathcal{E}_{2}$ lie in the same fiber of $a_{v}$ if and only if $F(\mathcal{E}_{1})$ and $F(\mathcal{E}_{2})$ lie in the same fiber of $a_{\widetilde{v}}$. As $a_{v}$ and $a_{\widetilde{v}}$ are isotrivial fibrations, it follows that the functor $F$ induces an injection $K_{v}(S,H)\longrightarrow K_{\widetilde{v}}(\widehat{S},\widehat{H})$. Since by Remark A.1 of \cite{PR2} we know that $K_{v}(S,H)$ and $K_{\widetilde{v}}(\widehat{S},\widehat{H})$ are irreducible and of the same dimension, the previous injection is an isomorphism and induces an isomorphism between the smooth loci.

Hence, it only remains to prove that the sheaf $F^{0}(\mathcal{E})$ is semistable. The proof will be by contradiction, supposing that $F^{0}(\mathcal{E})$ is not semistable. 

This implies that there is a desemistabilizing subsheaf $\mathcal{G}_{1}\subseteq F^{0}(\mathcal{E})$, i. e. a coherent subsheaf such that $p(\mathcal{G}_{1})>p(F^{0}(\mathcal{E}))$. We may and will choose it to be stable with maximal reduced Hilbert polynomial: such a $\mathcal{G}_{1}$ is the first term of a Jordan-H\"older filtration of the first term of a Harder-Narasimhan filtration of $F^{0}(\mathcal{E})$. 

We will moreover let $\mathcal{G}_{2}$ be the quotient of $F^{0}(\mathcal{E})$ by $\mathcal{G}_{1}$, so that we have an exact sequence
\begin{equation}
\label{eq:succesattag1}
0\longrightarrow\mathcal{G}_{1}\longrightarrow F^{0}(\mathcal{E})\longrightarrow\mathcal{G}_{2}\longrightarrow 0
\end{equation} 
Applying the functor $\widehat{F}$ to it, and using the fact that $\widehat{F}\circ F=[-2]$, we get the two exact sequences
\begin{equation}
	\label{eq:exfhat2a}
	0\longrightarrow\widehat{F}^{0}(\mathcal{G}_{2})\longrightarrow\widehat{F}^{1}(\mathcal{G}_{1})\longrightarrow 0
\end{equation}
\begin{equation}
	\label{eq:exfhat2b}
	0\longrightarrow\widehat{F}^{1}(\mathcal{G}_{2})\longrightarrow\widehat{F}^{2}(\mathcal{G}_{1})\longrightarrow\mathcal{E}\longrightarrow\widehat{F}^{2}(\mathcal{G}_{2})\longrightarrow 0
\end{equation}
The sheaf $\mathcal{G}:=\widehat{F}^{2}(\mathcal{G}_{1})/\widehat{F}^{1}(\mathcal{G}_{2})$ is then a subsheaf of $\mathcal{E}$. Our aim is to show that $p(\mathcal{G})>p(\mathcal{E})$: this would contradict the semistability of $\mathcal{E}$, hence completing the argument.

In order to prove that $p(\mathcal{G})>p(\mathcal{E})$, we start by collecting in the following Lemma some properties of the sheaf $\mathcal{G}_{1}$.

\begin{lem}
\label{lem:propg1}
The sheaf $\mathcal{G}_{1}$ is locally free, verifies WIT(2) with respect $\widehat{F}$, and if $v(\mathcal{G}_{1})=(a_{1},-n_{1},r_{1})$, then $a_{1},r_{1},n_{1}>0$, $a_{1}<a$, $n_{1}<n$ and either $n_{1}/a_{1}<n/a$, or $n_{1}/a_{1}=n/a$ and $r_{1}/a_{1}>r/a$.
\end{lem}

\proof The proof of this will be divided in six different steps.

\textit{Step 1: the sheaf $\mathcal{G}_{1}$ is locally free.} Indeed, otherwise $\mathcal{G}^{**}_{1}$ would be a locally free subsheaf of the locally free sheaf $F^{0}(\mathcal{E})$ with $p(\mathcal{G}^{**}_{1})>p(\mathcal{G}_{1})$, contradicting the maximality of $p(\mathcal{G}_{1})$.

\textit{Step 2: if $S$ is K3, then $-c_1(\mathcal{G}_1)$  is effective and nonzero, i.e. $n_{1}>0$}. As $S$ is K3, by the exact sequence (\ref{eq:exf}) there exists  an injection  $\mathcal{G}_{1}\subseteq F^{0}(\mathcal{E})\subseteq H^{0}(\mathcal{E})\otimes\mathcal{O}_{S}$. As a consequence, letting $V$ the generic quotient of $H^{0}(\mathcal{E})$ having rank $a_1$, we also have an injective morphism $\mathcal{G}_{1} \subseteq V\otimes\mathcal{O}_{S}$. 

Since $V\otimes\mathcal{O}_{S}$ is a trivial vector bundle having the same rank of $\mathcal{G}_{1}$, we deduce that $-c_1(\mathcal{G}_{1})$ is effective and is zero only if $\mathcal{G}_{1}$ is trivial. Finally $\mathcal{G}_{1}$ cannot be  trivial since it is contained $F^0(\mathcal{E})$, which is the kernel of the evaluation morphism $H^{0}(\mathcal{E})\otimes\mathcal{O}_{S}\rightarrow \mathcal{E}$ and cannot have nonzero global sections.

\textit{Step 3: if $S$ is Abelian, then $-c_1(\mathcal{G}_1)$  is effective and nonzero, i.e. $n_{1}>0$.} Let $Z\subseteq S$ be a reduced $0-$dimensional subscheme of degree $d\gg 0$, that we may and will choose so that for every $L\in Pic^{0}(S)$ no section of $\mathcal{E}\otimes L$ vanishes along $Z$. We let $\mathcal{E}_{Z}$ be the restriction of $\mathcal{E}$ to $Z$, and consider the exact sequence $$0\longrightarrow\mathcal{K}\longrightarrow\mathcal{E}\longrightarrow\mathcal{E}_{Z}\longrightarrow 0.$$By construction, for every $L\in Pic^{0}(S)$ we have $H^{0}(\mathcal{K}\otimes L)=0$. 

As a consequence we have $F^{0}(K)=0$ so, applying $F$ to the previous exact sequence, we get an inclusion $F^{0}(\mathcal{E})\subseteq F^{0}(\mathcal{E}_{Z})$. Now, notice that, as $Z$ is zero dimensional, $F^{0}(\mathcal{E}_{Z})$ is a direct sum of line bundles of degree 0 on $\widehat{S}$: as $\mathcal{G}_{1}\subseteq F^{0}(\mathcal{E})$, we then get an inclusion of $\mathcal{G}_{1}$ in a direct sum of line bundles of degree 0 on $\widehat{S}$. As a consequence, the rank $a_1$ vector bundle $\mathcal{G}_1$ also admits an injective morphism to direct sum $\oplus_{i=1}^{a_1}L_i$ of $a_1$ degree $0$ line bundles on $\widehat{S}$. This implies that 
$-c_1(\mathcal{G}_1)$ is represented by an effective divisor and is zero only if $\mathcal{G}_1\simeq \oplus_{i=1}^{a_1}L_i$.

Finally the last isomorphism does not hold since it would imply
that $F^0(\mathcal{E})\otimes L_i^{*}$ admits a nonzero section.
On the other hand, by projection furmula, $$H^0(F^0(\mathcal{E})\otimes L_i^{*})\simeq H^0(p_{S}^{*}(\mathcal{E})\otimes\mathcal{P}\otimes p_{\widehat{S}}^{*}(L^{*}))\simeq H^0(\mathcal{E}\otimes \iota^{*}\widehat{F}^0(L_i^{*}))$$
and the latter is zero since  $\widehat{F}^0(L_i^{*})=0$.

\textit{Step 4: the sheaf $\mathcal{G}_{1}$ verifies WIT(2) with respect to $\widehat{F}$.} First notice that, by definition of $\mathcal{G}_1$, its slope $\mu(\mathcal{G}_1)$  is maximal among the slopes of the subsheaves of $F^{0}(\mathcal{E})$, this implies that $\mu(\mathcal{G}_1)$ is an upper bound for the slopes of the subsheaves of $F^{0}(\mathcal{E})$ and $\mathcal{G}_2$. Since, by Steps 2 and 3, $n_1>0$, the class  $-c_1(\mathcal{G})$ is effective and nonzero,
all subsheaves of $F^{0}(\mathcal{E})$ and of $\mathcal{G}_{2}$ have strictly negative slope.

In particular, we see that $H^{0}(\mathcal{G}_{i})=0$ for $i=1,2$. This implies that $\widehat{F}^{0}(\mathcal{G}_{i})=0$ for $i=1,2$: if $S$ is K3 this is a consequence of the exact sequence (\ref{eq:exfhat}) applied to $\mathcal{G}_{i}$; if $S$ is Abelian, the same argument as before shows that $H^{0}(\mathcal{G}_{i}\otimes L)=0$ for every $L\in\widehat{S}$, hence $\widehat{F}^{0}(\mathcal{G}_{i})=0$ by cohomology and base change. 

Now, the exact sequence (\ref{eq:exfhat2a}) implies $\widehat{F}^{1}(\mathcal{G}_{1})=0$, and we conclude that  $\mathcal{G}_{1}$ verifies WIT(2) with respect to $\widehat{F}$.

\textit{Step 5: we have $r_{1}>0$.} As $\mathcal{G}_{1}$ verifies WIT(2) with respect to $\widehat{F}$, we see that $v(\widehat{F}^{2}(\mathcal{G}_{1}))=(r_{1},n_{1},a_{1})$ (this follows from the exact sequence (\ref{eq:exfhat}) applied to $\mathcal{G}_{1}$ if $S$ is K3, and it is well-known if $S$ is Abelian), hence $r_{1}\geq 0$. If $r_{1}=0$, the morphism $\widehat{F}^{2}(\mathcal{G}_{1})\longrightarrow\mathcal{E}$ in the exact sequence (\ref{eq:exfhat2b}) would be trivial (since $\mathcal{E}$ is torsion-free and $\widehat{F}^{2}(\mathcal{G}_{1})$ is torsion). As $\widehat{F}$ is fully faithful, this would imply that the inclusion morphism $\mathcal{G}_{1}\longrightarrow F^{0}(\mathcal{E})$ is trivial, getting a contradiction: it follows that $r_{1}>0$.

\textit{Step 6: conclusion of the proof}. We notice that $a_{1}$ is the rank of $\mathcal{G}_{1}$, which is a locally free subsheaf of $F^{0}(\mathcal{E})$, and $F^{0}(\mathcal{E})$ is a locally free sheaf of rank $a$. It follows that $0<a_{1}<a$. 
As $$v(F^{0}(\mathcal{E}))=\widetilde{v}=(a,-n,r),\,\,\,\,\,\,\, v(\mathcal{G}_{1})=(a_1,-n_1,r_1)$$
and as $p(\mathcal{G}_{1})>p(F^{0}(\mathcal{E}))$, we have either $n_{1}/a_{1}<n/a$ or $n_{1}/a_{1}=n/a$ and $r_{1}/a_{1}>r/a$. Finally $n_{1}<n$ follows from $n_{1}/a_{1}\leq n/a$ and $a_{1}<a$.\endproof

We will moreover need the following property of $\widehat{F}^{1}(\mathcal{G}_{2})$:

\begin{lem}
\label{lem:c1g2} If $\widehat{F}^{1}(\mathcal{G}_{2})\ne 0$ its first Chern class is strictly negative, i.e. 
$c_{1}(\widehat{F}^{1}(\mathcal{G}_{2}))=mh$ with $m<0$.
\end{lem}

\proof If $S$ is K3 the exact sequence (\ref{eq:exfhat}) applied to $\mathcal{G}_{2}$ shows that $\widehat{F}^{1}(\mathcal{G}_{2})$ is an extension of a subsheaf of $\mathcal{G}_{2}$ by $\mathcal{O}_{S}\otimes H^{1}(\mathcal{G}_{2})$ so, if $\widehat{F}^{1}(\mathcal{G}_{2})\ne 0$, its first Chern class is strictly negative. 

We now suppose that $S$ is Abelian. First of all, recall that $\widehat{F}^{0}(\mathcal{G}_{2})=0$ and that all the subsheaves of $\mathcal{G}_{2}$ have strictly negative first Chern class (see Step 4 of the proof of Lemma \ref{lem:propg1}). 

Since $\widehat{F}^{0}(\mathcal{G}_{2})=0$, it follows that $F^{i}(\widehat{F}^{0}(\mathcal{G}_{2}))=0$ for all $i$: the spectral sequence $$E_{2}^{p,q}:=F^{p}(\widehat{F}^{q}(\mathcal{G}_{2}))\longrightarrow E^{p+q}=(F\circ\widehat{F})^{p+q}(\mathcal{G}_{2})=\mathcal{G}_{2}[-2]^{p+q}$$then provides an inclusion $i:F^{1}(\widehat{F}^{1}(\mathcal{G}_{2}))\longrightarrow\mathcal{G}_{2}$ and the equality $F^{0}(\widehat{F}^{1}(\mathcal{G}_{2}))=0$.

In order to show that, if $\widehat{F}^{1}(\mathcal{G}_{2})\ne0$, its first Chern class is strictly negative, we first prove that $\widehat{F}^{1}(\mathcal{G}_{2})$ is torsion free. 

The sheaf  $\widehat{F}^{1}(\mathcal{G}_{2})$ cannot contain a zero dimensional subsheaf $T$ since $F^{0}(T)$ would be a nonzero semistable sheaf whose associated polystable sheaf is a direct sum of degree $0$ line bundles and $F^{0}(T)\subseteq F^{0}(\widehat{F}^{1}(\mathcal{G}_{2}))=0$. 

It follows that if $\widehat{F}^{1}(\mathcal{G}_{2})$ is not torsion free, there are a pure one dimensional sheaf
 $K$,  such that $c_1(K)=ph$ for a $p>0$, a torsion free sheaf $Q$ and an exact sequence 
 \begin{equation} 
\label{eq:exseqk}
0\longrightarrow K\longrightarrow\widehat{F}^{1}(\mathcal{G}_{2})\longrightarrow Q\longrightarrow 0.
\end{equation}
The sheaf $K$ verifies WIT(1) with respect to $F$: indeed, applying  $F$ to the sequence (\ref{eq:exseqk}) we obtain  $F^{0}(K)\subseteq F^{0}(\widehat{F}^{1}(\mathcal{G}_{2}))=0$ and, since $K$ is one dimensional, $F^2(K)=0$. As a consequence we get $c_1(F^1(K))=p\widehat{h}$. 

Moreover, 
since $Q$ is torsion free, the argument used for $\mathcal{E}$, in Step 3 of Lemma \ref{lem:propg1}, shows that $F^0(Q)$ is a subsheaf of a direct sum of line bundles belonging to  $Pic^0(\widehat{S})$: hence $c_1(F^0(Q))\le 0$.
 
 Applying the functor $F$ to the exact sequence (\ref{eq:exseqk}) we then get the exact sequence $$0\longrightarrow F^{0}(Q)\longrightarrow F^{1}(K)\longrightarrow F^{1}(\widehat{F}^{1}(\mathcal{G}_{2})),$$which yields an inclusion $$F^{1}(K)/F^{0}(Q)\longrightarrow F^{1}(\widehat{F}^{1}(\mathcal{G}_{2}))\stackrel{i}\longrightarrow\mathcal{G}_{2}.$$

Since  $c_{1}(F^{0}(Q))\leq 0$, we get $$c_{1}(F^{1}(K)/F^{0}(Q))\geq c_{1}(F^{1}(K))\geq 0,$$where the last inequality comes from the fact that $c_{1}(F^{1}(K))=p\widehat{h}$. It follows that $\mathcal{G}_{2}$  has a nonzero subsheaf with nonnegative first Chern class: but this is not possible (see Step 4 of the proof of Lemma \ref{lem:propg1}) and implies that $\widehat{F}^{1}(\mathcal{G}_{2})$  is torsion free.

It remains to prove that the first Chern class of the torsion free sheaf $\widehat{F}^{1}(\mathcal{G}_{2})$ is strictly negative: we will show it by distinguishing the  case where this sheaf is not $\mu-$semistable from the case  where it is $\mu-$semistable.

If $\widehat{F}^{1}(\mathcal{G}_{2})$  is not $\mu-$semistable and $c_1(\widehat{F}^{1}(\mathcal{G}_{2}))\ge 0$ there exist a $\mu-$stable sheaf $K$ with $c_1(K)=ph$ for $p>0$, a torsion free sheaf $Q$ and an exact sequence as in (\ref{eq:exseqk}). We can copy  the argument used above to show that the torsion of $\widehat{F}^{1}(\mathcal{G}_{2})$ cannot be pure of dimension $1$, the only difference  is that 
$F^{2}(K)=0$ follows from $\mu(K)>0$ and $\mu-$stability of $K$. As above we deduce that $\mathcal{G}_{2}$  has a nonzero subsheaf with nonnegative first Chern class. Since this is absurd, if  $\widehat{F}^{1}(\mathcal{G}_{2})$ is not   
$\mu-$semistable, its first Chern class is strictly negative.

Assume finally that  $\widehat{F}^{1}(\mathcal{G}_{2})$  is $\mu-$semistable and $c_1(\widehat{F}^{1}(\mathcal{G}_{2}))=mh$ for $m\ge0$.  By $\mu-$stability of $\widehat{F}^{1}(\mathcal{G}_{2})$, if $m>0$, there are no nonzero morphism from   $\widehat{F}^{1}(\mathcal{G}_{2})$ to any $L\in Pic^0(S)$ and, if $m=0$, the locus of $Pic^0(S)$ consisting of line bundles admitting nontrivial morphisms from $\widehat{F}^{1}(\mathcal{G}_{2})$ is, at most, finite.
This implies that the support of $F^2(\widehat{F}^{1}(\mathcal{G}_{2}))$ is empty or finite and, since $F^0(\widehat{F}^{1}(\mathcal{G}_{2}))=0$, we obtain $c_1(F^1(\widehat{F}^{1}(\mathcal{G}_{2})))=m\widehat{h}$. Again this is absurd since $i:F^1(\widehat{F}^{1}(\mathcal{G}_{2}))\rightarrow \mathcal{G}_2$
is an injection and $\mathcal{G}_2$ cannot contain subsheaves with nonnegative first Chern class. Hence, also in the case where  $\widehat{F}^{1}(\mathcal{G}_{2})$  is $\mu-$semistable, its first Chern class is srtictly negative. \endproof

We are now ready to conclude the proof of Proposition \ref{prop:yoshi2}. We first notice that as $\mathcal{G}=\widehat{F}^{2}(\mathcal{G}_{1})/\widehat{F}^{1}(\mathcal{G}_{2})$, we have $$c_{1}(\mathcal{G})=c_{1}(\widehat{F}^{2}(\mathcal{G}_{1}))-c_{1}(\widehat{F}^{1}(\mathcal{G}_{2})),$$hence by Lemma \ref{lem:c1g2} it follows that $p(\mathcal{G})\ge p(\widehat{F}^{2}(\mathcal{G}_{1}))$. 

Let us now prove that $p(\widehat{F}^{2}(\mathcal{G}_{1}))>p(\mathcal{E})$: it will then follow that $p(\mathcal{G})>p(\mathcal{E})$, concluding the contradiction argument.

To show that $p(\widehat{F}^{2}(\mathcal{G}_{1}))>p(\mathcal{E})$, recall first that $v(\widehat{F}^{2}(\mathcal{G}_{1}))=(r_{1},n_{1},a_{1})$ and $v(\mathcal{E})=(r,n,a)$.
Moreover, setting  $l=h^2/2$ (hence also $l=\widehat{h}^2/2$ if $S$ is abelian), by hypothesis and Lemma \ref{lem:propg1}:
\begin{enumerate}
 \item we have $r,n,a>0$ and $ln^{2}-ra=v^{2}/2=k>0$;
 \item we have $r_{1},n_{1},a_{1}>0$;
 \item we have $a_{1}<a$, and either $n_{1}/a_{1}<n/a$, or $n_{1}/a_{1}=n/a$ and $r_{1}/a_{1}>r/a$;
 \item we have $ln_{1}^{2}-r_{1}a_{1}=v(\mathcal{G}_{1})^{2}/2\geq -1$, since $\mathcal{G}_{1}$ is $H-$stable (if $S$ is Abelian, we even have $ln_{1}^{2}-r_{1}a_{1}\geq 0$).
\end{enumerate}
By Lemma \ref{lem:numeri} below, it follows that there is $n_{0}\in\mathbb{N}$ such that for every $n>n_{0}$, we have $n_{1}/r_{1}\geq n/r$, and if $n_{1}/r_{1}=n/r$, then $a_{1}/r_{1}>a/r$. This  exactly  means that $p(\widehat{F}^{2}(\mathcal{G}_{1}))>p(\mathcal{E})$, completing the proof.\endproof

We now prove the following, which is used to conclude the proof of the previous Lemma:

\begin{lem}
\label{lem:numeri}Fix $k,l,r\in\mathbb{N}$, $k,l,r>0$, and let $n,a,r_{1},a_{1},n_{1}\in\mathbb{N}$, $n,a,r_{1},a_{1},n_{1}>0$, be such that the following conditions are fulfilled:
\begin{enumerate}
 \item $ln^{2}-ra=k$;
 \item $ln_{1}^{2}-r_{1}a_{1}\geq -1$;
 \item $a_{1}<a$;
 \item $n_{1}/a_{1}<n/a$ or $n_{1}/a_{1}=n/a$ and $r_{1}/a_{1}>r/a$. 
\end{enumerate}
If $n>32r^{3}k$, then either $n_{1}/r_{1}>n/r$, or $n_{1}/r_{1}=n/r$ and $a_{1}/r_{1}>a/r$.
\end{lem}

\proof We let $k_{1}:=ln_{1}^{2}-r_{1}a_{1}$, so that $k_{1}\geq -1$. As $n_{1}/a_{1}\leq n/a$, it follows that $n_{1}/n\leq a_{1}/a$. Moreover, as $a_{1}=\frac{ln^{2}_{1}-k_{1}}{r_{1}}$ and $a=\frac{ln^{2}-k}{r}$, we get the inequality $$\frac{n_{1}}{n}\leq\frac{r}{r_{1}}\cdot\frac{n_{1}}{n}\cdot\frac{n_{1}-\frac{k_{1}}{ln_{1}}}{n-\frac{k}{ln}}.$$  This implies that 
\begin{equation}
\label{eq:in0}
1\leq\frac{r}{r_{1}}\cdot\frac{n_{1}-\frac{k_{1}}{ln_{1}}}{n-\frac{k}{ln}}.\end{equation}
We claim that as $n>32r^{3}k$ we have $r>r_{1}$. Indeed, as $n>32r^{3}k\geq 3k$ and $k_{1}\geq -1$, we get 
\begin{equation}
\label{eq:in2}
\frac{n_{1}-\frac{k_{1}}{ln_{1}}}{n-\frac{k}{ln}}\leq\frac{n_{1}+\frac{1}{ln_{1}}}{n-\frac{1}{3l}}\leq\frac{n_{1}+\frac{1}{n_{1}}}{n_{1}+\frac{2}{3}}
\end{equation}
where the last inequality follows from the fact that $n>n_{1}$ (so that $n-n_{1}\geq 1$). 

As $n>n_{1}$ and $n>3k$, the first term of the inequality (\ref{eq:in2}) is strictly smaller than 1: this can be checked by analyzing the third term if $n_{1}\geq 2$ and the second term
if $n_{1}=1$ (recall that  under our assumption $n\ge 3$). In any case we get $$\frac{n_{1}-\frac{k_{1}}{ln_{1}}}{n-\frac{k}{ln}}<1,$$hence the inequality (\ref{eq:in0}) gives $r>r_{1}$. 

We now write the inequality (\ref{eq:in0}) in a different form. More precisely, we have
$$1\leq\frac{\frac{n_{1}}{r_{1}}}{\frac{n}{r}}\cdot\frac{1-\frac{k_{1}}{ln_{1}^{2}}}{1-\frac{k}{ln^{2}}}\leq\frac{\frac{n_{1}}{r_{1}}}{\frac{n}{r}}\cdot\frac{1+\frac{1}{ln_{1}^{2}}}{1-\frac{k}{ln^{2}}}$$where the last equality follows from $k_{1}\geq-1$. As $n>32r^{3}k$ we see that $1-\frac{k}{ln^{2}}>0$, hence the previous inequality becomes 
\begin{equation}
\label{eq:in0rew}
\frac{\frac{n_{1}}{r_{1}}}{\frac{n}{r}}\geq\frac{1-\frac{k}{ln^{2}}}{1+\frac{1}{ln_{1}^{2}}}=1-\frac{\frac{k}{ln^{2}}+\frac{1}{ln_{1}^{2}}}{1+\frac{1}{ln_{1}^{2}}}.
\end{equation}

We first want to show that $n_{1}/r_{1}\geq n/r$. As the first term of the inequality (\ref{eq:in0rew}) is an integral multiple of $\frac{1}{r_{1}n}$, by the inequality (\ref{eq:in0rew}) it is enough to show that $$\frac{\frac{k}{ln^{2}}+\frac{1}{ln_{1}^{2}}}{1+\frac{1}{ln_{1}^{2}}}\leq\frac{1}{r_{1}n}.$$

To do so, notice first that $1+\frac{1}{ln_{1}^{2}}>1$, and that as $1/ln_{1}^{2}\leq 1$ we have $$\frac{\frac{1}{ln_{1}^{2}}}{1+\frac{1}{ln_{1}^{2}}}\leq\frac{1}{2}.$$Moreover, as $n>32r^{3}k$, we get $\frac{k}{ln^{2}}<1/4$. We then find that $$1-\frac{\frac{k}{ln^{2}}+\frac{1}{ln_{1}^{2}}}{1+\frac{1}{ln_{1}^{2}}}>\frac{1}{4},$$so that , by inequality (\ref{eq:in0rew}), we finally get $$\frac{n_{1}}{r_{1}}>\frac{n}{4r}.$$This implies that $n_{1}>\frac{n}{4r}$, hence we get $$\frac{1}{ln_{1}^{2}}<\frac{1}{l\frac{n^{2}}{16r^{2}}}<\frac{1}{ln\frac{32r^{3}k}{16r^{2}}}=\frac{1}{2lrnk}\leq\frac{1}{2rn}.$$

Now, using again  $n>32r^{3}k$, we even get that $k/ln^{2}<\frac{1}{32r^{3}n}$, hence  $$\frac{\frac{k}{ln^{2}}+\frac{1}{ln_{1}^{2}}}{1+\frac{1}{ln_{1}^{2}}}<\frac{k}{ln^{2}}+\frac{1}{ln_{1}^{2}}<\frac{1}{32r^{3}n}+\frac{1}{2rn}<\frac{1}{rn}<\frac{1}{r_{1}n},$$where the last inequality comes from $r>r_{1}$.

We then have $n_{1}/r_{1}\geq n/r$ if $n>32r^{3}k$. To complete the proof, notice that if $n_{1}/r_{1}=n/r$, then $r_{1}/r=n_{1}/n<a_{1}/a$. But this means that $a_{1}/r_{1}>a/r$, and we are done.\endproof

We conclude this section with a proposition and a corollary, which allows us to pass from a Mukai vector of rank 0 to a Mukai vector of strictly positive rank (see Proposition 3.14 of \cite{Y3} for a proof for stable sheaves).

We will have to deal also with surfaces whose N\'eron-Severi group has rank bigger than one and, in this case, we need an explicit lower bound on the Euler characteristic of the rank zero Mukai vector which makes possible  checking   genericity of the polarizations  for both Mukai vectors.

In order to give this bound we recall that if $S$ is smooth projective surface, $H$ is an ample divisor on $S$ and $\xi\in NS(S)$ is the class of an effective curve, the set of numerical equivalence classes of effective curves $C$ on $S$ such that $C\cdot H\leq \xi\cdot H$ is finite.

It follows that the set $$\bigg\{ \frac{(C^{2}+2)(\xi\cdot H)}{2C\cdot H}\,\bigg|\,C\subseteq S\,\,{\rm is}\,\,{\rm an}\,\,{\rm effective}\,\,{\rm curve},\,\,C\cdot H\le\xi\cdot H\bigg\}$$ is bounded  and we  denote by $N_{S,H,\xi}$ its maximum.

\begin{prop}
	\label{prop:yoshi1}
	Let $S$ be a projective K3 or Abelian surface and  set $F:=F_{K3}$ if $S$ is a K3 and $F:=F_{Ab}$ if $S$ is Abelian. Let $H$ a primitive ample line bundle on $S$, and set $h:=c_{1}(H)$. Let $k\in\mathbb{N}^{*}$ and let $\xi\in NS(S)$ be the first Chern class of an effective divisor on $S$ such that $\xi^{2}=2k$. Finally, let $v=(0,\xi,a)$ be a Mukai vector on $S$ and assume that WIT(0) holds with respect to $F$  for every $H-$semistable sheaf $\mathcal{E}$ with Mukai vector $v$, and that $F^{0}(\mathcal{E})$ is locally free (e.g. assume $a>p_{0}$ for  $p_{0}$ as in Lemma \ref{lem:step1}(2)).
	\begin{enumerate}
		\item If $S$ is K3, $a>N_{S,H,\xi}$ and $H$ is both $v-$generic and $\widetilde{v}-$generic, the functor $F$ induces isomorphisms $M_{v}(S,H)\simeq M_{\widetilde{v}}(S,H)$ and $M^{s}_{v}(S,H)\simeq M^{s}_{\widetilde{v}}(S,H)$.
		\item If $S$ is Abelian, $a>N_{S,H,\xi}$, $H$ is $v-$generic and $\widehat{H}$ is $\widetilde{v}-$generic, the functor $F$ induces isomorphisms $M_{v}(S,H)\simeq M_{\widetilde{v}}(\widehat{S},\widehat{H})$, $M^{s}_{v}(S,H)\simeq M^{s}_{\widetilde{v}}(\widehat{S},\widehat{H})$, $K_{v}(S,H)\simeq K_{\widetilde{v}}(\widehat{S},\widehat{H})$ and $K^{s}_{v}(S,H)\simeq K^{s}_{\widetilde{v}}(\widehat{S},\widehat{H})$.
	\end{enumerate} 
\end{prop}

\proof As in the proof of Proposition \ref{prop:yoshi2}, in order to get the statement we just need to prove that for $a>N_{S,H,\xi}$ the sheaf $F^{0}(\mathcal{E})$ is $H-$semistable in the K3 case and $\widehat{H}-$semistable in the Abelian case.

The proof is by contradiction: we suppose that $F^{0}(\mathcal{E})$ is not semistable and contradict in several steps the $H-$semistability of $\mathcal{E}$. 

So, suppose that $\mathcal{E}$ is $H-$semistable with Mukai vector $v=(0,\xi,a)$ where $a>N_{S,H,\xi}$, and suppose that $F^{0}(\mathcal{E})$ is not $H-$semistable if $S$ is K3, and that it is not $\widehat{H}-$semistable if $S$ is Abelian. 

Then $F^{0}(\mathcal{E})$ has a desemistabilizing  stable locally free subsheaf $\mathcal{G}_1$ with maximal reduced Hilbert polynomial and we set 
$\mathcal{G}_2:=F^{0}(\mathcal{E})/\mathcal{G}_1$. If $S$ is a K3 we write $v(\mathcal{G}_1)=(a_1,-\xi_1,r_1)$, if $S$ is Abelian
we write $v(\mathcal{G}_1)=(a_1,-\widehat{\xi}_1,r_1)$. 

\textit{Step 1: if $S$ is K3, the class  $\xi_1$ is effective and $\xi_1\neq 0$.} This follows by repeating the argument of Step 2 of Lemma \ref{lem:propg1}.

\textit{Step 2: if $S$ is Abelian, the class  $\widehat{\xi}_1$ is effective and $\widehat{\xi}_1\neq 0$.} This follows as in  Step 3 of Lemma \ref{lem:propg1}.

\textit{Step 3:  $r_1\le 0$ for $a>N_{S,H,\xi}$.} Suppose first that $S$ is a K3 surface. Since  $\mathcal{G}_1$ is $H-$stable  we have $v(\mathcal{G}_1)^{2}\geq -2$. Hence, if $r_1\geq 1$ we get 
$$-2\leq v(\mathcal{G}_1)^{2}=(\xi_1)^{2}-2a_1r_1\leq(\xi_1)^{2}-2a_1.$$
Setting $d:=\xi\cdot H$ and $d_1=\xi_1\cdot H$, the inequality  $p_{H}(\mathcal{G}_1)>p_{H}(F^{0}(\mathcal{E}))$ implies $-d_1/a_1\geq -d/a$ and hence $a_1\geq a\frac{d_1}{d}$
, so that $$-2\leq(\xi_1)^{2}-2a\frac{d_1}{d}.$$
Moreover as $a_1<a$, the inequality $a_1\geq a\frac{d_1}{d}$ implies $d_1<d$

As $\xi_1$ is an effective divisor such that $\xi_1\cdot H=d_1<d$, we have $N_{S,H,\xi}\ge \frac{\xi_1^{2}+2}{2d_1}\cdot d$ and since $a>N_{S,H,\xi}$ we obtain
$$-2\le (\xi_1)^{2}-2a\frac{d_1}{d}<(\xi_1)^{2}-2\frac{((\xi_1)^{2}+2)d}{2d_1}\cdot\frac{d_1}{d}=-2,$$getting a contradiction.

The same argument works also if $S$ is Abelian, simply by replacing any occurrence of $\xi ,\xi_1,H$ and $S$ by $\widehat{\xi}, \widehat{\xi}_1,\widehat{H}$ and $\widehat{S}$.

\textit{Step 4:  the sheaf $\mathcal{G}_1$ verifies WIT(2) with respect to $\widehat{F}$, and in particular $r_1=0$.} The first part follows as in Step 4 
of Lemma \ref{lem:propg1}. Since $r_1=rk(\widehat{F}(\mathcal{G}_1))$, it has to be $0$ by Step 3.

\textit{Step 5: the sheaf $\widehat{F}^2(\mathcal{G}_1)$ is a subsheaf of $\mathcal{E}$.} By maximality of the reduced Hilbert polynomial of $\mathcal{G}_1$, the sheaf $\mathcal{G}_2:=F^{0}(\mathcal{E})/\mathcal{G}_1$ is torsion-free. By applying the functor $\widehat{F}$ to the exact sequence $$0\longrightarrow\mathcal{G}_1\longrightarrow F^{0}(\mathcal{E})\longrightarrow\mathcal{G}_2\longrightarrow 0$$ we get the exact sequence
\begin{equation}
	\label{eq:exactrank0}
	0\longrightarrow\widehat{F}^{1}(\mathcal{G}_2)\longrightarrow\widehat{F}^{2}(\mathcal{G}_1)\longrightarrow\mathcal{E}\longrightarrow\widehat{F}^{2}(\mathcal{G}_2)\longrightarrow 0.
\end{equation}
so it suffices to show that  $\widehat{F}^{1}(\mathcal{G}_2)=0$. 

First of all, notice that by Step 4 the rank of $\widehat{F}^{2}(\mathcal{G}_1)$ is $r_1=0$. As $\widehat{F}^{1}(\mathcal{G}_2)\subseteq\widehat{F}^{2}(\mathcal{G}_1)$, we know that $\widehat{F}^{1}(\mathcal{G}_2)$ is a torsion sheaf.

If $S$ is a K3 surface,  the exact sequence (\ref{eq:exfhat}) applied to $\mathcal{G}_2$ shows that $\widehat{F}^{1}(\mathcal{G}_2)$ is an extension of a subsheaf of $\mathcal{G}_2$ by a locally free sheaf. As $\mathcal{G}_2$ is torsion-free, it follows that $\widehat{F}^{1}(\mathcal{G}_2)$ is torsion-free and, since it is also a torsion sheaf, we obtain $\widehat{F}^{1}(\mathcal{G}_2)=0$.

If $S$ is Abelian, the maximality of the reduced Hilbert polynomial of $\mathcal{G}_1$ among the Hilbert polynomials of the subsheaves of $F(\mathcal{E})$ implies that 
the slope of any subsheaf of $\mathcal{G}_2$ cannot be bigger than the slope of $\mathcal{G}_1$, which is strictly negative by Step 2. Hence any subsheaf of $\mathcal{G}_2$ has strictly negative slope and $H^0(\mathcal{G}_2\otimes L)=0$ for every 
$L\in Pic^0(\widehat{S})$: it follows that $\widehat{F}^{0}(\mathcal{G}_2)=0$. 

As in the proof Lemma \ref{lem:c1g2}, the spectral sequence $$E_{2}^{p,q}:=F^{p}(\widehat{F}^{q}(\mathcal{G}_2))\longrightarrow E^{p+q}=(F\circ\widehat{F})^{p+q}(\mathcal{G}_2)=\mathcal{G}_2[-2]^{p+q}$$then provides an inclusion $$i:F^{1}(\widehat{F}^{1}(\mathcal{G}_2))\longrightarrow\mathcal{G}_2,$$and the equality $F^{0}(\widehat{F}^{1}(\mathcal{G}_2))=0$. Since the support of the torsion sheaf $\widehat{F}^{1}(\mathcal{G}_2)$ has dimension at most $1$ we also have 
$F^{2}(\widehat{F}^{1}(\mathcal{G}_2))=0$ and 
$\widehat{F}^{1}(\mathcal{G}_2)$ satisfies WIT(1) with respect to  $F$. The equality $F^{0}(\widehat{F}^{1}(\mathcal{G}_2))=0$ implies that $\widehat{F}^{1}(\mathcal{G}_2)$ has no zero-dimensional torsion: it follows that if $\widehat{F}^{1}(\mathcal{G}_2)\neq 0$ then $c_1(\widehat{F}^{1}(\mathcal{G}_2))$ would be  the class of an effective curve and the same would hold
for $c_1(F^{1}(\widehat{F}^{1}(\mathcal{G}_2)))$. But this is impossible since $F^{1}(\widehat{F}^{1}(\mathcal{G}_2))\subseteq \mathcal{G}_2$ and the slope of every subsheaf of 
$\mathcal{G}_2$ is negative. It follows $\widehat{F}^{1}(\mathcal{G}_2)$=0.

\textit{Step 6: conclusion of the proof.} Suppose that $S$ is K3. Recall that $v(F(\mathcal{E}))=\widetilde{v}=(a,-\xi,0)$ and, by Step 4, $v(\mathcal{G}_1)=(a_1,-\xi_1,0)$ (with $\xi_1$ effective and nonzero by Step 1). Since $\mathcal{G}_1$ is a desemistabilizing subsheaf of  $F(\mathcal{E})$ we have  $\frac{-d_1}{a_1}> \frac{-d}{a}$ or equivalently 
\begin{equation}
	\label{eu}\frac{a_1}{d_1}>\frac{a}{d}.
\end{equation} 
By Step 4 the sheaf $\mathcal{G}_1$  verifies WIT(2) with respect to $\widehat{F}$, hence $v(\widehat{F}^2(\mathcal{G}_1))= (0,\xi_1,a_1)$ and, by Step 5, the sheaf $\widehat{F}^2(\mathcal{G}_1)$ is a subsheaf of $\mathcal{E}$: the inequality (\ref{eu}) then implies that $\widehat{F}(\mathcal{G}_1)$ is a desemistabilizing subsheaf for $\mathcal{E}$, which is not possible since $\mathcal{E}$ is semistable. 

The proof in the Abelian case is similiar, simply replace $\xi$ and $\xi_1$ with $\widehat{\xi}$ and $\widehat{\xi}_1$ respectively where necessary.\endproof

\begin{oss}
	\label{oss:+} {\rm In the important case where  
		$\rho(S)=1$, Proposition \ref{prop:yoshi1} simply says that if $v=(0,\xi,a)$ and $a\gg0$, the functor $F$ induces  isomorphism $M_v(S,H)\simeq M_{\widetilde{v}}(S,H)$ if $S$ is a K3 and  $K_v(S,H)\simeq K_{\widetilde{v}}(\widehat{S},\widehat{H})$ if $S$ is Abelian.
		This easier statement allows  to prove Theorem \ref{thm:main} with exception of the case where the rank of $v$ is $0$ and $\rho(S)>1$. To handle this remaining case we need the following corollary of Proposition \ref{prop:yoshi1}.}
\end{oss}

\begin{cor}
	\label{cor:yoshi1}
	Let $S$ be a projective K3 or Abelian surface, $k\in\mathbb{N}^{*}$ and $\xi\in NS(S)$ the first Chern class of an effective divisor on $S$ such that $\xi^{2}=2k$.
	For a nonzero  $a\in \mathbb{Z}$, let $v=(0,\xi,a)$ be a Mukai vector on $S$ and let $H$ be a $v-$generic polarization. For $d\in\mathbb{Z}$ set 
	$v_{dH}:=v\cdot ch(\mathcal{O}(dH))=(0,\xi, a+d\xi\cdot H)$.
	\begin{enumerate}
		\item If $S$ is K3, there exists  $d_0\in\mathbb{N}$ such that  $a+d_{0}\xi\cdot H>0$ and a $\widetilde{v_{d_{0}H}}-$generic polarization $H'$ such that  the functor $F_{K3}$ induces isomorphisms $M_{v_{d_{0}H}}(S,H)\simeq M_{\widetilde{v_{d_{0}H}}}(S,H')$ and $M^{s}_{v_{d_{0}H}}(S,H)\simeq M^{s}_{\widetilde{v_{d_{0}H}}}(S,H')$.
		\item If $S$ is Abelian, 
		there exists  $d_0\in\mathbb{N}$ such that  $a+d_{0}\xi\cdot H>0$ and a $\widetilde{v_{d_{0}H}}-$generic polarization $\widehat{H}'$ such that  the functor $F_{Ab}$
		induces isomorphisms $M_{v_{d_{0}H}}(S,H)\simeq M_{\widetilde{v_{d_{0}H}}}(\widehat{S},\widehat{H}')$, $M^{s}_{v_{d_{0}H}}(S,H)\simeq M^{s}_{\widetilde{v_{d_{0}H}}}(\widehat{S},\widehat{H}')$, $K_{v_{d_{0}H}}(S,H)\simeq K_{\widetilde{v_{d_{0}H}}}(\widehat{S},\widehat{H}')$ and $K^{s}_{v_{d_{0}H}}(S,H)\simeq K^{s}_{\widetilde{v_{d_{0}H}}}(\widehat{S},\widehat{H}')$.
	\end{enumerate} 
\end{cor}
\proof We only deal with the case where $S$ is a K3, the Abelian case being very similar. 

We first claim that if $U$ is a small compact neighborhood of $H$ in the ample cone of $S$, then there is $N_U \in \mathbb{N}$ such that $N_U>N_{S,L,\xi}$ for every $L\in U$. To prove this, let $\overline{NE(S)}$ be the closure in $NS(S)\otimes\mathbb{R}$ of the cone of effective curves of $S$. For every $L\in Amp(S)$ we define $F_{L}:\overline{NE(S)}\longrightarrow\mathbb{R}$ as $F_{L}(\alpha):=\alpha\cdot L$ for every $\alpha\in\overline{NE(S)}$, and let $$D_{L}:=\{\alpha\in\overline{NE(S)}\,|\,F_{L}(\alpha)\leq \xi\cdot L\}.$$ As a consequence of Kleiman's ampleness  criterion (see Corollary 1.19 (2) of \cite{KM2}) the locus $D_{L}$ is compact for every $L\in U$. As $U$ is compact, it follows that $\bigcup_{L\in U}D_{L}$ is compact as well, so its subset $Y$ of integral nonzero classes is finite. 

Fix now a nonzero class $C\in\overline{NE(S)}$. The function $$f_{C}:U\longrightarrow\mathbb{R},\,\,\,\,\,\,\,\,\,\,f_{C}(L):=\frac{(C^2+2)(\xi\cdot L)}{2C\cdot L}$$is continuous. As $U$ is compact, the function $f_{C}$ has a maximum and, since $Y$ is finite, there exists $N_U>N_{S,L,\xi}$ for every $L\in U$.

Choose now $d_0\in\mathbb{N}$ such that  $a+d_{0}\xi\cdot H> N_U$ and such that every $H-$semistable sheaf $\mathcal{E}$
with Mukai vector $v_{d_{0}H}$ satisfies WIT(0) with respect to $F_{K3}$ and $F_{K3}(\mathcal{E})$ is locally free (this is possible by Lemma \ref{lem:step1}(2)). 

By point (2) of Lemma \ref{lem:tensogen} the polarization $H$ is also $v_{d_{0}H}-$generic, so it lies in a $v_{d_{0}H}-$chamber $\mathcal{C}$: as the $\widetilde{v_{d_{0}H}}-$walls are locally finite, there is a $\widetilde{v_{d_{0}H}}-$generic polarization $H'\in U\cap\mathcal{C}$. But then $M_{v_{d_{0}H}}(S,H)=M_{v_{d_{0}H}}(S,H')$ (as $H,H'\in\mathcal{C}$), and as $a+d_{0}\xi\cdot H> N_{S,H',\xi}$ (since $H'\in U$), by Proposition \ref{prop:yoshi1} the functor $F_{K3}$ induces an isomorphism between $M_{v_{d_{0}H}}(S,H')$ and 
$M_{\widetilde{v_{d_{0}H}}}(S,H')$, concluding the proof.\endproof

\subsection{The proof of Theorem \ref{thm:main}}

This section is devoted to the proof of Theorem \ref{thm:main}: the goal is to show that if $S$
is a projective K3 (resp. Abelian) surface, $v$ is a Mukai vector on $S$ of type $(m,k)$  and $H$ is a polarization on $S$ which is general with respect to $v$, the locally trivial deformation equivalence class of $M_{v}(S,H)$ only depends on $(m,k)$. Before giving the proof, we provide several results we will need. 

\subsubsection{Changing polarization and first Chern class.}

We first show the following Lemma, which allows us, if the rank of the Mukai vector is strictly positive, to suppose that the first Chern class of the Mukai vector is a multiple of the polarization. As a consequence, this will allow us to suppose the N\'eron-Severi group of $S$ to have rank 1.

\begin{lem}
\label{lem:polgenxi}Let $m,k\in\mathbb{N}$, $m,k>0$, $S$ a projective K3 or Abelian surface, $v$ a Mukai vector on $S$ of type $(m,k)$ and of the form $v=m(r,\xi,a)$ with $r>0$. Let $g:=gcd(r,\xi)$, and suppose that $H$ is a $v-$generic polarization on $S$. Moreover, suppose that $\rho(S)\geq 2$, and let $\mathcal{C}$ be the $v-$chamber such that $H\in\mathcal{C}$. Then there exists a Mukai vector $v'=m(r,\xi',a')$ and a primitive polarization $H'$ in $\mathcal{C}$ such that:
\begin{enumerate}
 \item $v'$ is equivalent to $v$;
 \item $\xi'=gc_{1}(H')$;
\item $(H')^{2}\gg 0$.
\end{enumerate}
In particular $M_{v}(S,H)\simeq M_{v'}(S,H')$ and $M^{s}_{v}(S,H)\simeq M^{s}_{v'}(S',H')$. If $S$ is Abelian, we have $K_{v}(S,H)\simeq K_{v'}(S,H')$ and $K^{s}_{v}(S,H)\simeq K^{s}_{v'}(S,H')$.
\end{lem}

\proof This is a generalization of Lemma II.6 of \cite{OG1}. First, notice that as $g:=gcd(r,\xi)$ there are two coprime integers $s,p\in\mathbb{N}$ and a primitive class $\zeta\in NS(S)$ such that $r=gs$ and $\xi=gp\zeta$.

Replacing $H$ by another polarization inside $\mathcal{C}$ if necessary, 
we may suppose that $\xi\notin\mathbb{R}\cdot c_{1}(H)$. Moreover, since $\mathcal{C}$ is an open cone in $NS(S)$, we may also assume that the sublattice $\Lambda\subseteq NS(S)$ spanned by $c_{1}(H)$ and $\zeta$ is saturated.

Now, let $d\in\mathbb{N}$ and set $v':=v\cdot ch(\mathcal{O}_{S}(dH))$. Then $v'$ is equivalent to $v$, and $$v'=m(r,\xi+rdc_{1}(H),a+d\xi \cdot H+rd^{2}H^{2}/2)=$$ $$=m(gs,gp\zeta+gsdc_{1}(H),a+dgp\zeta \cdot H+gsd^{2}H^{2}/2).$$
Notice now that, if $d\gg0$ then $p\zeta+sdc_{1}(H)\in\mathcal{C }$: if we let $H'$ be an ample divisor such that $c_{1}(H')=p\zeta+sdc_{1}(H)$, we have 
$$\xi':=gp\zeta+gsdc_{1}(H)=g(p\zeta+sdc_{1}(H))=gc_1(H').$$
As $gcd(s,p)=1$, if we choose $d$ such that  $gcd(d,p)=1$, the class $c_1(H')$ is primitive in $\Lambda$ and, since $\Lambda$ is saturated, it is primitive in $NS(S)$.
The isomorphism $M_{v}(S,H)\simeq M_{v'}(S,H')$ follows from (2) of Lemma \ref{lem:tensorlb}.

Finally, as $d\gg 0$, we have $$(H')^{2}=p^2\zeta^{2}+2psd\zeta\cdot H+s^{2}d^{2} H^{2}\gg 0.$$

To conclude the proof, we notice that if $S$ is Abelian, then by point (1) of Lemma \ref{lem:tensorlb} the tensorization with $\mathcal{O}_{S}(dH)$ induces an isomorphism between the fibers of the corresponding Yoshioka fibrations.\endproof  

\subsubsection{Deformation to elliptic surfaces.}

Elliptic surfaces having a section and whose Picard number is 2 prove to be particularly useful, as in this case we have a privileged class of polarizations, called $v-$suitable. Let $Y$ be an elliptic K3 or Abelian surface such that $NS(Y)=\mathbb{Z}\cdot f\oplus\mathbb{Z}\cdot\sigma$, where $f$ is the class of a fiber and $\sigma$ is the class of a section. Let $v$ be a Mukai vector on $Y$ of the form $v=(r,\xi,a)$ with $r>0$, and recall the following definition (see \cite{OG1}):

\begin{defn}
{\rm A polarization $H$ on $Y$ is called $v-$\textit{suitable} if $H$ is in the unique $v-$chamber whose closure contains $f$.}
\end{defn}

We have an easy numerical criterion to guarantee that a polarization on $Y$ is $v-$suitable (see Lemma I.0.3 of \cite{OG1} for K3 surfaces, and point (2) of Lemma 2.24 of \cite{PR} for Abelian surfaces):

\begin{lem}
\label{lem:vsuit}Let $Y$ be aprojective elliptic K3  or an Abelian surface  with $NS(Y)=\mathbb{Z}\cdot\sigma\oplus\mathbb{Z}\cdot f$, where $\sigma$ is the class of a section and $f$ is the class of  a fibre, and let $v=(r,\xi,a)$ be a Mukai vector on $Y$ such that $r>0$. Let $H$ be a polarization, and suppose that $c_{1}(H)=\sigma+tf$ for some $t\in\mathbb{Z}$. 
\begin{enumerate}
 \item If $Y$ is K3, then $H$ is $v-$suitable if $t\geq|v|+1$.
 \item If $Y$ is Abelian, then $H$ is $v-$suitable if $t\geq |v|$.
\end{enumerate}
\end{lem}

In the next Lemma, by deforming a triple $(S,v,H)$ to a triple $(Y,v',H')$ where $Y$ is an elliptic surface and $H'$ is $v'-$suitable, we show that the locally trivial deformation class of $M_{v}(S,H)$ (resp. $K_{v}$) only depends on numerical data associated with $v$: in the particular and important  case where the rank $r$ is strictly positive  and prime to the first Chern class of $v$ it only depends on $r$ and $v^2$. 

\begin{lem}
\label{lem:defor}Let $m,k\in\mathbb{N}$, $m,k>0$, and for $i=1,2$ let $S_{i}$ be a projective K3 (resp. Abelian) surface, $v_{i}$ a Mukai vector on $S_{i}$ of type $(m,k)$ and $H_{i}$ a $v_{i}-$generic polarization on $S_{i}$. Write $v_{i}=m(r_{i},\xi_{i},a_{i})$ for $i=1,2$, and suppose that the following conditions are verified:
\begin{enumerate}
 \item $r_{1}=r_{2}=:r>0$;
 \item $gcd(r,\xi_{1})=gcd(r,\xi_{2})=:g$;
 \item $a_{1}\equiv a_{2}$ mod $g$;
\end{enumerate}
Then $M_{v_{1}}(S_1,H_1)$ and $M_{v_{2}}(S_2,H_2)$ (resp. $K_{v_{1}}(S_1,H_1)$ and $K_{v_{2}}(S_2,H_2)$) are locally trivially deformation equivalent. 
\end{lem}

\proof The argument we present here was first used by O'Grady in \cite{OG1} and by Yoshioka in \cite{Y1} for primitive Mukai vectors, and by the authors in \cite{PR} in the case of $m=2$ and $k=1$. 

First, we may assume $\rho(S_{i})>1$. Indeed, consider a deformation $(\mathscr{X}_{i},\mathscr{L}_{i},\mathscr{H}_{i})$ of the triple $(S_i,v_i,H_i)$ over a smooth connected curve $C_{i}$ inducing a nontrivial deformation of  $S_i$. 
By the Main Theorem of \cite{O}, the locus parametrizing points $t\in C_{i}$ such that $\rho(\mathscr{X}_{i,t})>1$ is dense in the classical topology of $C_{i}$. Since the locus of points $t\in C_{i}$ such that $\mathscr{H}_{i,t}$ is not general with respect to $v_{i,t}$ is finite (see Remark \ref{oss:ap}) and since by Lemma \ref{lem:loctriv}(1) and Lemma \ref{lem:relgenvgen} we may replace a polarization which is general with respect to Mukai vector $v$ by a $v-$generic polarization, we may suppose  $\rho(S_{i})>1$.

Since the class of $a_i$ modulo $g$ does not change replacing $v_i$ by an equivalent Muaki vector, 
by Lemma \ref{lem:polgenxi} we may even suppose that $v_{i}=m(r,g c_{1}(H_{i}),a_{i})$, where $H_{i}$ is ample and $H_{i}^{2}=2d_{i}$ with $d_{i}\gg 0$. 

Let $Y$ be a K3 or an  Abelian surface with an elliptic fibration and such that $NS(Y)=\mathbb{Z}\cdot\sigma\oplus\mathbb{Z}\cdot f$, where $f$ is the class of a fiber, and $\sigma$ is the class of a section. 
For $i=1,2$, because of the connectedness of the moduli spaces of polarized K3 or Abelian surfaces, there is a smooth, connected curve $T_{i}$ and a deformation $(\mathscr{X}_{i},\mathscr{L}_{i},\mathscr{H}_{i})$ over $T_{i}$ of $(S_{i},v_{i},H_{i})$ such that there is $t_i\in T_{i}$ with the property $(\mathscr{X}_{i,t_i},v_{i,t_i},H_{i,t_i})=(Y,v'_{i},H_{i}')$, where  
\begin{enumerate}
 \item $c_{1}(H'_{i})=\sigma+p_{i}f$, 
 \item $v'_{i}=m(r,gc_{1}(H'_{i}),a_{i})$.
\end{enumerate}
By Lemma \ref{lem:loctriv}(1) the varieties $M_{v_{i}}(S_i,H_i)$ and $M_{v'_{i}}(Y,H_{i}')$ are locally trivially deformation equivalent.
Let $\xi'_{i}:=c_{1}(H'_{i})$. Notice that $(v'_{1})^{2}=(v'_{2})^{2}$ and they have the same rank: hence $|v'_{1}|=|v'_{2}|$, so by Lemma \ref{lem:vsuit} a polarization is $v'_{1}-$suitable if and only if it is $v'_{2}-$suitable. 

Notice that
$p_i=d-\sigma^2/2$, i.e. $p_i=d_i+1$ if $Y$ is a K3 and $p_i=d_i$ if $Y$ is Abelian: hence $p_i\ge d\gg0$ in both cases and, by Lemma \ref{lem:vsuit}, we have that $H'_{i}$ is $v'_{i}-$suitable for $i=1,2$, hence $H'_{1}$ and $H'_{2}$ are in the same $v'_{1}-$chamber $\mathcal{C}$. Using Proposition \ref{prop:chamber} we then change to a common generic polarization $H\in\mathcal{C}$, which is $v'_{i}-$generic for $i=1,2$.

The statement for  $M_{v_{1}}(S_1,H_1)$ and $M_{v_{2}}(S_2,H_2)$ follows if we show that $v'_1=m(r,g(\sigma+p_{1}f),a_{1})$ and $v'_2=m(r,g(\sigma+p_{2}f),a_{2})$ are equivalent since, by point (2) of Lemma \ref{lem:tensorlb}, this implies that $M_{v'_{1}}(Y,H)\simeq M_{v'_{2}}(Y,H)$.
 
 As $(v'_{1})^{2}=(v'_{2})^{2}$, we have $$g^2\sigma^2+2g^2p_1-2ra_1=g^2\sigma^2+2g^2p_2-2ra_2,$$ hence $$2r(a_2-a_1)=2g^2(p_2-p_1)$$ and, since $g$ divides $a_2-a_1$, $r$ divides $g(p_2-p_1)$: so there exists $l\in \mathbb{Z}$ such that $g(p_2-p_1)=rl$. Letting $L\in Pic(Y)$ be a line bundle whose first Chern class is $lf$, the last equality implies 
$v'_2$ and $v'_1\cdot ch(L)$ have the same component in $NS(Y)$: since they also have the same rank and the same square we get $v'_2=v'_1\cdot ch(L)$  and
$v'_2$ and $v'_1$ are equivalent.

Replacing Lemma \ref{lem:loctriv}(1) by Lemma \ref{lem:loctriv}(2), as the tensorization with a line bundle preserves the Yoshioka fibration the same argument shows that, if $S_{1}$ and $S_{2}$ are both Abelian surfaces, then $K_{v_{1}}(S_1,H_1)$ and $K_{v_{2}}(S_2,H_2)$ are locally trivially deformation equivalent.\endproof

\subsubsection{A numerical result on equivalent Mukai vectors.} 

The following numerical Lemma, together with Propositions \ref{prop:yoshi2} and \ref{prop:yoshi1}, will allow to show that certain  moduli spaces parametrizing  semistable sheaves with different ranks, on a K3 or Abelian surface whose N\'eron-Severi group has rank 1, are isomorphic.

\begin{lem}
\label{lem:primi}Let $m,k\in\mathbb{N}$, $m,k>0$, $S$ a projective K3 or Abelian surface with $NS(S)=\mathbb{Z}\cdot h$, where $h=c_{1}(H)$ and $H$ an ample line bundle. Let $v=m(r,nh,a)$ be a Mukai vector of type $(m,k)$, where $r>0$. For every $s\in\mathbb{Z}$ let $$v_{sH}:=v\cdot ch(\mathcal{O}_{S}(sH))=m(r,n_{s}h,a_{s}).$$
\begin{enumerate}
 \item For every $N\in\mathbb{N}$ there is $s>N$ such that $n_{s}\gg 0$ and $gcd(n_{s},a_{s})=1$.
 \item If $n=1$ and $a=0$, then for every $N\in\mathbb{N}$ there is $s>N$ such that $n_{s}\gg 0$, $gcd(n_{s},a_{s})=1$ and $a_{s}\in 2k\mathbb{Z}$.
\end{enumerate}
\end{lem}

\proof Write $H^{2}=2l$. A direct computation shows that 
\begin{equation}
	\label{eq:nt}
	n_{s}=n+rs=n_{s-1}+r,
\end{equation}
\begin{equation}
	\label{eq:at}
	a_{s}=a+2lns+rls^{2}.
\end{equation} 
Since $2ln^2-2ra= 2ln_s^2-2ra_s=2k$, if there exists $s$ such that $n_s$ and $a_s$ are not coprime, a common prime divisor $p$ has to divide $k$. 

Let $p$ be a prime factor of $k$ dividing both $n_{s_{0}}$ and $a_{s_{0}}$ for some $s_0\in\mathbb{Z}$.
Since the cup product with the Chern character of a line bundle gives an isometry of the Mukai lattice, and since the Mukai vector $(r,nh,a)$ is primitive, 
the Mukai vector $(r,n_{s_0}h,a_{s_0})$ is primitive as well, so $p$ if divides both $n_{s_0}$ and $a_{s_0}$, then it cannot divide $r$.

Now, if $p$ also divides $n_{s}$ and $a_{s}$, then $p$ has to divide the difference  $n_s-n_{s_0}=r(s-s_0)$, hence $p$ is a factor of $s-s_0$ and  we have $s=s_{0}+pm$ for $m\in\mathbb{Z}$. 

Since the prime divisors of $k$ are finite, by choosing $s$ outside the union of a finite number of arithmetic sequences we obtain that $gcd(n_{s},a_{s})=1$; moreover, we may choose $s\gg 0$ in order to obtain $n_s\gg0$.

If we now suppose $n=1$ and $a=0$ (so that $l=k$), then equations (\ref{eq:nt}) and (\ref{eq:at}) give $n_{s}=1+rs$ and $a_{s}=2ls+rls^{2}=2ks+rks^{2}$.

If $p$ is a prime number dividing both $a_{s}$ and $n_{s}$, then $p$ divides $k$. Choosing $s=2ks'$ for some $s'\in\mathbb{Z}$, then any prime dividing $k$ cannot divide $n_{s}=1+2ks'r$, hence $a_{s}$ and $n_{s}$ are coprime. Moreover, since $s$ is even, $2k$ divides $a_s$ and, finally, $n_s\gg0$ if $s'$ is chosen big enough.\endproof

\subsubsection{Conclusion of the proof of Theorem \ref{thm:main}.}

We now proceed with the proof of:\vspace{3mm}

\noindent\textbf{Theorem \ref{thm:main}.} \textit{Let $m,k\in\mathbb{N}$, $m,k>0$, and for $i=1,2$ let $S_{i}$ be a projective K3 or Abelian surface, $v_{i}$ a Mukai vector on $S_{i}$ of type $(m,k)$ and $H_{i}$ a polarization on $S_{i}$ which is general with respect to $v_{i}$.
\begin{enumerate}
 \item If $S_{1}$ and $S_{2}$ are both K3 surfaces or both Abelian surfaces, then $M_{v_{1}}(S_1,H_1)$ and $M_{v_{2}}(S_2,H_2)$ are deformation equivalent, and the deformation is locally trivial.
 \item If $S_{1}$ and $S_{2}$ are two Abelian surfaces, then $K_{v_{1}}(S_1,H_1)$ and $K_{v_{2}}(S_2,H_2)$ are deformation equivalent, and the deformation is locally trivial.
\end{enumerate}}

\proof We start with the proof of (1). The proof for K3 and for Abelian surfaces is formally the same: to make notation easier we only discuss the case of K3 surfaces and leave to the reader the obvious modifications for the Abelian case.  

For every $\ell\in\mathbb{N}$, $\ell>0$, let $X_{\ell}$ be a projective K3 with $NS(X_{\ell})=\mathbb{Z}\cdot h_{\ell}$, where $h_{\ell}=c_{1}(H_{\ell})$ and $H_{\ell}$ is an ample divisor with $H_{\ell}^{2}=2\ell$. Let $u_{\ell}:=m(0,h_{\ell},0)$, so $u_{\ell}$ is a Mukai vector of type $(m,\ell)$ and $H_{\ell}$ is $u_{\ell}-$generic.

We will show that if $S$ is a K3 surface, $v=m(r,\xi,a)$ is a Mukai vector of type $(m,k)$ on $S$ and $H$ is a polarization on $S$ which is general with respect to $v$, then $M_{v}(S,H)$ is locally trivially  deformation equivalent to $M_{u_{k}}(X_{k},H_{k})$.

We first notice that we can assume $H$ to be $v-$generic polarization and $a\ne0$ if $r=0$: this follows from  Lemma \ref{lem:relgenvgen}, if $r>0$, and from Lemmata \ref{lem:relgenvgen} and \ref{lem:tensorh0}, if $r=0$. Under these assumptions, we may now start the proof, that will be in several steps: in the first, we reduce to $r>0$; in the second, we reduce to $r$ and $\xi$ relatively prime; in the third we reduce to $r\in 2k\mathbb{N}$; the fourth step concludes the proof.

\textit{Step 1: reduction to $r>0$}. If $v=m(0,\xi,a)$, where $\xi$ is effective and $a\neq 0$, by Lemma \ref{lem:tensorlb}(1), for $d\in\mathbb{Z}$ the tensorization by $\mathcal{O}(dH)$ induces an isomorphism between $M_v(S,H)$
and $M_{v_{dH}}(S,H)$. By Corollary \ref{cor:yoshi1} there exists $d_0\in \mathbb{N}$ and a $\widetilde{v_{d_{0}H}}-$generic polarization $H'$ such that 
$F_{K3}$ induces an isomorphism between $M_{v_{d_{0}H}}(S,H)$ and  $M_{\widetilde{v_{d_{0}H}}}(S,H')$ and the rank of  $\widetilde{v_{d_{0}H}}$ is strictly positive.

\textit{Step 2: reduction to $r\gg 0$ and prime with $\xi$}. By Step 1 we may suppose $r>0$. By Lemma \ref{lem:polgenxi} we may even suppose $v=m(r,\xi,a)$ with $\xi=n c_{1}(H)$. We set $\ell:=\frac{H^{2}}{2}$ and consider the Mukai vector $v'=m(r,nh_{\ell},a)$ on $X_{\ell}$. Notice that $v'$ is of type $(m,k)$ and that $H_{\ell}$ is $v'-$generic. 

As the moduli spaces of polarized K3 are connected, the moduli spaces $M_{v}(S,H)$ and $M_{v'}(X_{\ell},H_{\ell})$ are locally trivially deformation equivalent by Lemma \ref{lem:loctriv}. For $s\in\mathbb{Z}$ write $v'_{s}:=v'_{sH_{l}}=m(r,n_{s}h_{l},a_{s})$: by Lemma \ref{lem:tensorlb}(1) we have $M_{v'}(X_\ell,H_\ell)\simeq M_{v'_{s}}(X_\ell,H_\ell)$,
and by Lemma \ref{lem:primi}(1) there is $s\in\mathbb{Z}$ such that $n_{s}\gg 0$ and $gcd(n_{s},a_{s})=1$.

As $n_{s}\gg 0$, by Proposition \ref{prop:yoshi2}, we get $M_{v'_{s}}(X_\ell,H_\ell)\simeq M_{\widetilde{v'_{s}}}(X_\ell,H_\ell)$ and moreover, since $m,k,r$ are fixed, we may also assume $a_s\gg0$. Finally, as $\widetilde{v'_{s}}=m(a_{s},n_{s}h_{l},r)$, $gcd(a_{s},n_{s})=1$ and $a_{s}\gg 0$, we conclude that it is sufficient to prove the theorem for Mukai vectors of the form $m(r,\xi,a)$ with $r\gg 0$ and prime with $\xi$.

\textit{Step 3: reduction to $r\in 2k\mathbb{Z}$, $r\gg 0$ and prime with $\xi$}. Since, for  a nontrivial deformation of a K3, the locus of the base where the rank of the N\'eron-Severi group jumps is dense in the classical topology (see \cite{O}) and,  by Proposition \ref{prop:generalopen}, the locus where the polarization is not general  is closed, deforming $S$ if necessary, by Lemma \ref{lem:loctriv}, 
we may assume that $\rho(S)\ge 2$.

By Step 2 and Lemma \ref{lem:polgenxi}, we can suppose $v=m(r,c_{1}(H),a)$ where $r\gg 0$. On the surface $X_{k}$ we consider the Mukai vector $v''=m(r,h_{k},0)$, which is of type $(m,k)$ and $H_{k}$ is $v''-$generic. By Lemma \ref{lem:defor} we know that $M_{v}(S,H)$ and $M_{v''}(X_{k},H_{k})$ are locally trivially deformation equivalent. 

Let now $s\in\mathbb{Z}$, and $v''_{s}:=v''_{sH_{k}}=m(r,n_{s}h_{k},a_{s})$, where $n_{s}=1+rs$ and $a_{s}=2ks+rks^{2}$. By Lemma \ref{lem:tensorlb}(1) we have $M_{v''}(X_k,H_k)\simeq M_{v''_{s}}(X_k,H_k)$, and by Lemma \ref{lem:primi}(2) we can choose $s$ such that $n_{s}\gg 0$, $a_{s}\in 2k\mathbb{Z}$ and $gcd(n_{s},a_{s})=1$. 

Moreover, as $n_{s}\gg 0$, by Proposition \ref{prop:yoshi2} we have $M_{v''_{s}}(X_k,H_k)\simeq M_{\widetilde{v''_{s}}}(X_k,H_k)$. But $\widetilde{v''_{s}}=m(a_{s},n_{s}h_{k},r)$, and $a_{s}\in 2k\mathbb{Z}$ and, moreover $a_s>n_s\gg0$. In conclusion, we just need to prove the theorem for Mukai vectors of the form $m(r,\xi,a)$ such that $r\gg 0$, $r$ is an even multiple of $k$, and $r$ and $\xi$ are coprime.

\textit{Step 4: conclusion}. By Step 3 we may suppose $v=m(r,\xi,a)$, with $r=2kp$ is prime with $\xi$ and with $p\gg 0$. We show that $M_{v}(S,H)$ is locally trivially deformation equivalent to $M_{u_{k}}(X_{k},H_{k})$.

By Lemma \ref{lem:defor} we know that $M_{v}(S,H)$ is locally trivially deformation equivalent to $M_{v'''}(X_{k},H_{k})$, where $v'''=m(2kp,h_{k},0)$. As $p\gg 0$, by Proposition \ref{prop:yoshi1}, we have $M_{v'''}(X_{k},H_{k})\simeq M_{\widetilde{v'''}}(X_{k},H_{k})$, where $\widetilde{v'''}=m(0,h_{k},2kp)$. But now notice that $\widetilde{v'''}\cdot ch(\mathcal{O}_{X_{k}}(-pH_{k}))=u_{k}$, so by Lemma \ref{lem:tensorlb}(1) we have $M_{\widetilde{v'''}}(X_{k},H_{k})\simeq M_{u_{k}}(X_{k},H_{k})$, concluding the proof of point (1) of the statement. 

We are now in the position to prove point (2) of the statement. To do so, notice that the equivalence in (1) in the Abelian case is obtained using deformations of the moduli spaces induced by deformations of the corresponding triple along smooth, connected varieties, and isomorphism between moduli spaces induced either by tensor products with line bundles or by the Fourier-Mukai transform whose kernel is the Poincar\'e line bundle. Since the Yoshioka fibration is preserved by these isomorphisms and, by Lemmata \ref{lem:redflat} and \ref{lem:loctriv}, a deformation of $(S,v,H)$ along a smooth connected variety also induces a locally trivial deformation of the fibers of the Yoshioka fibrations,  
point (2) of the statement is implied by (1).\endproof

\section{The moduli spaces are irreducible symplectic varieties}

This section is devoted to the proof of Theorem \ref{thm:mainb1}: if $m,k\in\mathbb{N}$, $m,k>0$, $S$ is a projective K3 or Abelian surface, $v$ is a Mukai vector on $S$ of type $(m,k)$ and $H$ is a polarization on $S$ that is general with respect to $v$, then $M_{v}(S,H)$, if $S$ is a K3, and $K_{v}(S,H)$, if $S$ is Abelian, are irreducible symplectic varieties.

To do so, we first show in Section 3.1 that, if $S$ is a projective K3 surface and $H$ is general with respect to the Mukai vector $v$, then $M_{v}(S,H)$ and $M^{s}_{v}(S,H)$ are simply connected. Similarly, if $S$ is an Abelian surface, then $K_{v}(S,H)$ and $K^{s}_{v}(S,H)$ are simply connected (with the exception of $(m,k)=(2,1)$, where $K_{v}(S,H)$ is still simply connected, but the fundamental group of $K_{v}^{s}(S,H)$ is $\mathbb{Z}/2\mathbb{Z}$). 

This will allow us to show that the exterior algebra of reflexive forms on any finite quasi-\'etale cover of $M_{v}(S,H)$ (resp. of $K_{v}(S,H)$) is generated by the reflexive pull-back of a symplectic form on $M_{v}(S,H)$ (resp. on $K_{v}(S,H)$): this will be done in section 3.2, by showing that for a particular choice of $S$, $v$ and $H$, there is a rational dominant map from  a moduli space $M_{u}(S,H)$ (resp. $K_{u}(S,H)$) with primitive Mukai vector, to the moduli space $M_{v}(S,H)$ (resp. $K_{v}(S,H)$).

\subsection{Simple connectedness}

We first show in this section that, for $H$ general with respect to $v$, the moduli spaces $M_{v}(S,H)$ and $M^{s}_{v}(S,H)$ (resp. $K_{v}(S,H)$ and $K^{s}_{v}(S,H)$) are simply connected. We will divide the proof of this in two main parts: the first one is devoted to the case of K3 surfaces; in the second we will consider Abelian surfaces. In both cases, the proof has the same structure: we first show the simple connectedness of the moduli space associated to a particular choice of the surface $S$, of the Mukai vector $v$ and of the polarization, and then use Theorem \ref{thm:main} to conclude.

Before doing this, we state and prove the following result about the codimension of the subset of reducible curves in the linear system of a multiple of the polarization on a K3 surface with Picard number 1.

\begin{lem}
\label{lem:codimr}Let $S$ be a projective K3 surface or an Abelian surface such that $NS(S)=\mathbb{Z}\cdot h$, where $h=c_{1}(H)$ and $H$ is an ample line bundle such that $H^{2}=2k$. Let $m\in\mathbb{N}$, $m>0$, and consider the subset $R\subseteq|mH|$ parameterizing reducible curves. If $(m,k)\neq(2,1)$, then ${\rm codim}_{|mH|}(R)\geq 2$.
\end{lem}

\proof If $C\in|mH|$ is a reducible curve, then if $S$ is K3 (resp. if $S$ is Abelian) there must be $1\leq m_{1},m_{2}\leq m$ such that $m=m_{1}+m_{2}$ and two curves $C_{1}\in|m_{1}H|$ and $C_{2}\in|m_{2}H|$ (resp. $C_{1}\in|m_{1}H+L|$ and $C_{2}\in|m_{2}H-L|$ for some $L\in\widehat{S}$) such that $C=C_{1}+C_{2}$.

For every $1\leq m_{1},m_{1}\leq m$ such that $m_{1}+m_{2}=m$ we let $$P_{m_{1},m_{2}}:=\left\{\begin{array}{ll} |m_{1}H|\times|m_{2}H|, & S\,\,{\rm is}\,\,{\rm K3}\\ \prod_{L\in\widehat{S}}|m_{1}H+L|\times|m_{2}H-L|, & S\,\,{\rm is}\,\,{\rm Abelian}\end{array}\right.$$We then get that $$R=\bigcup_{\substack{1\leq m_{1},m_{2}\leq m,\\ m_{1}+m_{2}=m}}P_{m_{1},m_{2}}.$$

Notice that $$\dim(|pH|)=\left\{\begin{array}{ll} 1+kp^{2}, & S\,\,{\rm is}\,\,{\rm K3}\\ kp^{2}-1, & S\,\,{\rm is}\,\,{\rm Abelian}\end{array}\right.$$and that if $S$ is Abelian and $L\in\widehat{S}$, then $\dim(|pH|)=\dim(|pH|\pm L)$. It follows that $$\dim(P_{m_{1},m_{2}})=\left\{\begin{array}{ll} 2+k(m_{1}^{2}+m_{2}^{2}), & S\,\,{\rm is}\,\,{\rm K3}\\ k(m_{1}^{2}+m_{2}^{2}), & S\,\,{\rm is}\,\,{\rm Abelian}\end{array}\right.$$and the codimension of $P_{m_{1},m_{2}}$ in $|mH|$ is $2km_{1}m_{2}-1$.

Hence, in order for $R$ to have codimension 1 in $|mH|$, there must be $1\leq m_{1},m_{2}\leq m$ such that $m_{1}+m_{2}=m$, and such that $2m_{1}m_{2}k-1=1$. Hence $m_{1},m_{2},k=1$, so that $m=2$ and $k=1$. Thus, if $(m,k)\neq(2,1)$ we get co$\dim_{|mH|}(R)\geq 2$.\endproof

Under the hypothesis of Lemma \ref{lem:codimr}, we see that the only case where $R$ is a divisor in $|mH|$ is when $(m,k)=(2,1)$.

\subsubsection{The case of K3 surfaces}

Let $X$ be a projective K3 surface with $Pic(X)=\mathbb{Z}\cdot\mathcal{O}_{X}(H)$, where $H$ is an ample divisor such that $H^{2}=2k$. We let $h:=c_{1}(H)$, and we choose $m\in\mathbb{N}$, $m>0$.

We let $V$ be the open subset of $|mH|$ of smooth curves, and $U$ the open subset of $|mH|$ of integral curves. For $u=m(0,h,0)$ we will consider the morphism $p_{u}:M_{u}(X,H)\longrightarrow|mH|$ mapping a sheaf to its Fitting subscheme (see Corollary 20.5 of \cite{EH}, and \cite{LP}).

Let $\mathcal{J}_{V}:=p_{u}^{-1}(V)$ and $\mathcal{J}_{U}:=p_{u}^{-1}(U)$, which are two open subsets of $M_{u}(X,H)$. Notice that if $C\in V$, then $\mathcal{F}\in p_{u}^{-1}(C)$ if and only if there is $L\in Pic(C)$ of degree $m^{2}k$ such that $\mathcal{F}=j_{*}L$, where $j:C\longrightarrow X$ is the inclusion. In particular, we have an isomorphism\footnote{Here and in what follows, if $C$ is a smooth projective curve and $d\in\mathbb{Z}$, we let $Pic^{d}(C)$ be the set of line bundles of degree $d$ on $C$.} $$p_{u}^{-1}(C)\longrightarrow Pic^{m^{2}k}(C)$$obtained by mapping $\mathcal{F}=j_{*}L$ to $L$.

Moreover, if $C\in U$, then $\mathcal{F}\in p_{u}^{-1}(C)$ if and only if $\mathcal{F}=j_{*}L$, where $j:C\longrightarrow X$ is the inclusion, and $L$ is a rank one torsion-free sheaf on $C$ of degree $m^{2}k$, i. e. such that $\chi(L)=0$. We notice that all these sheaves are $H-$stable of Mukai vector $u$, hence we have $$\mathcal{J}_{V}\subseteq\mathcal{J}_{U}\subseteq M^{s}_{u}(X,H)\subseteq M_{u}(X,H).$$

We start by showing the following (the proof is a generalization of the argument proposed in section 4 of \cite{OG2}).

\begin{prop}
\label{prop:simplu}The moduli spaces $M_{u}(X,H)$ and $M^{s}_{u}(X,H)$ are simply connected.
\end{prop}

\proof Notice that $u$ is a Mukai vector of type $(m,k)$ and that $H$ is $u-$generic. If $m=1$, then $M_{u}(X,H)=M^{s}_{u}(X,H)$: this is an irreducible symplectic manifold, and we are done. For $(m,k)=(2,1)$, see section 4 of \cite{OG2}.

For $m\geq 2$, we have that $M_{u}(X,H)$ is a normal, irreducible projective variety (by Theorem 4.4 of \cite{KLS}). Since for normal quasi-projective varieties the inclusion of an open subvariety induces a surjection on the fundamental groups (see Proposition 2.10 of \cite{Koll}) the chain of inclusions $$\mathcal{J}_{V}\stackrel{j}\hookrightarrow\mathcal{J}_{U}\hookrightarrow M^{s}_{u}(X,H)\hookrightarrow M_{u}(X,H)$$ of smooth open subvarieties of $M_{u}(X,H)$ given before induces a chain of surjections $$\pi_{1}(\mathcal{J}_{V})\stackrel{\pi_{1}(j)}\longrightarrow\pi_{1}(\mathcal{J}_{U})\longrightarrow\pi_{1}(M^{s}_{u}(X,H))\longrightarrow\pi_{1}(M_{u}(X,H)).$$We then just need to show that $\pi_{1}(j)$ is the trivial map.

To show this, notice that the homotopy exact sequence of the fibration $p_{u|\mathcal{J}_{V}}:\mathcal{J}_{V}\longrightarrow V$ gives the exact sequence $$\pi_{1}(p_{u}^{-1}(C))\longrightarrow\pi_{1}(\mathcal{J}_{V})\longrightarrow\pi_{1}(V)\longrightarrow\{1\},$$where $C\in V$. As remarked above we have $p^{-1}(C)\simeq Pic^{m^{2}k}(C)$, hence the exact sequence is 
\begin{equation}
\label{eq:expi1}
\pi_{1}(Pic^{m^{2}k}(C))\stackrel{j_{C}}\longrightarrow\pi_{1}(\mathcal{J}_{V})\longrightarrow\pi_{1}(V)\longrightarrow\{1\}.
\end{equation}

We start by proving the following:

\begin{lem}
\label{lem:jc}The morphism $\pi_{1}(j)\circ j_{C}:\pi_{1}(Pic^{m^{2}k}(C))\longrightarrow\pi_{1}(\mathcal{J}_{U})$ is trivial.
\end{lem}

\proof Let $\ell\subseteq|mH|$ be a generic line, and suppose it is generated by two smooth curves intersecting transversally. By Lemma \ref{lem:codimr}, we can suppose that all the curves in $\ell$ are reduced and irreducible. 

If $b:\widetilde{X}\longrightarrow X$ is the blow-up of $X$ along the base locus $Bs(\ell)$ of $\ell$, then $\widetilde{X}$ is the total space of $\ell$: this means that for every $s\in\widetilde{X}$ there is a unique curve $C_{s}$ of $\ell$ such that $s\in\widetilde{C}_{s}$, where $\widetilde{C}_{s}$ is the proper transform of $C_{s}$. We have a natural fibration $p_{\ell}:\widetilde{X}\longrightarrow\ell$ mapping $s\in\widetilde{X}$ to $C_{s}$.

We now define an embedding $g:\widetilde{X}\longrightarrow\mathcal{J}_{U}$ of fibrations over $\ell$. First, fix $p\in Bs(\ell)$, and let $d:=1+m^{2}k$. Let $s\in\widetilde{X}$: then $b(s)\in C_{s}$, and consider the rank one torsion-free sheaf $L_{s}:=\mathcal{I}_{b(s)}\otimes\mathcal{O}_{C_{s}}(dp)$, whose degree on $C_{s}$ is $m^{2}k$. If $j_{s}:C_{s}\longrightarrow X$ is the inclusion, then $j_{s*}(L_{s})\in\mathcal{J}_{U}$, and we let $g(s):=j_{*}L_{s}$. The inclusion $g$ then fits in a commutative diagram (where $i$ is the inclusion)
$$\begin{CD}
\widetilde{X} @>g>> \mathcal{J}_{U}\\
@V\text{$p_{\ell}$}VV @VV\text{$p_{u}$}V\\
\ell @>>\text{$i$}> U
\end{CD}$$

Notice that if $t\in\ell$ is a generic point and $C$ is the corresponding curve in $\ell$, then $p_{\ell}^{-1}(t)=\widetilde{C}$, the proper transform of $C$ under $b$, while $p_{u}^{-1}(t)\simeq Pic^{m^{2}k}(C)$. The restriction $g_{t}:\widetilde{C}\longrightarrow Pic^{m^{2}k}(C)$ of $g$ to $p_{\ell}^{-1}(t)$ can be identified with the Abel-Jacobi map from $C$ to its Jacobian. It then induces a surjective morphism $\pi_{1}(g_{t}):\pi_{1}(\widetilde{C})\longrightarrow\pi_{1}(Pic^{m^{2}k}(C))$.

Now, let $C\in\ell$ be a smooth curve. We have a commutative diagram 
$$\begin{CD} 
\widetilde{C}  @>\text{$g_{t}$}>> Pic^{m^{2}k}(C)\\ 
@ViVV @VV\text{$\widetilde{i}$}V\\ 
\widetilde{X} @>>g> \mathcal{J}_{U}
\end{CD},$$inducing a commutative diagram 
$$\begin{CD} \pi_{1}(\widetilde{C}) @>\text{$\pi_{1}(g_{t})$}>> \pi_{1}(Pic^{m^{2}k}(C))\\ 
@V\text{$\pi_{1}(i)$}VV @VV\text{$\pi_{1}(j)\circ j_{C}$}V\\ 
\pi_{1}(\widetilde{X}) @>>\text{$\pi_{1}(g)$}> \pi_{1}(\mathcal{J}_{U})
\end{CD}.$$

As $\pi_{1}(\widetilde{X})=\{1\}$ and the morphism $\pi_{1}(g_{t})$ is surjective, it follows that $\pi_{1}(j)\circ j_{C}$ is trivial, thus concluding the proof.\endproof

An immediate consequence of Lemma \ref{lem:jc} is that the surjective morphism $\pi_{1}(j)$ factors through a surjective morphism $$\overline{\pi_{1}(j)}:\pi_{1}(\mathcal{J}_{V})/{\rm im}(j_{C})\longrightarrow\pi_{1}(\mathcal{J}_{U}).$$The exact sequence (\ref{eq:expi1}) gives an isomorphism between $\pi_{1}(V)$ and $\pi_{1}(\mathcal{J}_{V})/{\rm im}(j_{C})$, hence we get a surjective map $\iota:\pi_{1}(V)\longrightarrow\pi_{1}(\mathcal{J}_{U})$, which is then trivial if and only if $\pi_{1}(j)$ is trivial: we then just need to show that $\iota$ is trivial. 

To do so, consider the generic line $\ell\subseteq|mH|$ of the proof of Lemma \ref{lem:jc}: all the curves parametrized by $\ell$ are reduced and irreducible, and we can suppose that it is transversal to $W:=U\setminus V$, where $\ell\cap W:=\{x_{1},...,x_{p}\}$ is given by smooth points of $W$.
 
As $\ell$ is generic, by the Zariski Main Theorem (see Theorem 3.22 of \cite{Vo}), the inclusion of $\ell\setminus W$ in $V$ gives a surjection $\pi_{1}(\ell\setminus W)\longrightarrow\pi_{1}(V)$, hence we finally get a surjective morphism $\iota_{\ell}:\pi_{1}(\ell\setminus W)\longrightarrow\pi_{1}(\mathcal{J}_{U})$, and we just need to show that $\iota_{\ell}$ is trivial. More precisely, if $\gamma_{1},\cdots,\gamma_{p}$ are the generators of $\pi_{1}(\ell\setminus W)$, we need to show that $\iota_{\ell}(\gamma_{i})$ is trivial.

Now, notice that the fibration $p_{\ell}:\widetilde{X}\longrightarrow\ell$ has a section $\sigma_{\ell}$ (fixing $p\in Bs(\ell)$, we let $\sigma_{\ell}(t):=\pi^{-1}(p)\cap p_{\ell}^{-1}(t)$). Hence every $\gamma_{i}$ has a lifting $\widetilde{\gamma}_{i}$ in $\pi_{1}(\widetilde{X})$, and by construction its image in $\pi_{1}(\mathcal{J}_{U})$ under $\pi_{1}(g)$ is $\iota_{\ell}(\gamma_{i})$. But as $\pi_{1}(\widetilde{X})=\pi_{1}(X)$ (since $\pi:\widetilde{X}\longrightarrow X$ is a blow-up) and as $\pi_{1}(X)$ is trivial (since $X$ is K3), it follows that $\iota_{\ell}(\gamma_{i})=0$.\endproof

The main consequence of Proposition \ref{prop:simplu} is that the moduli spaces of (semi)stable sheaves associated to $(m,k)-$triples are simply connected:

\begin{thm}
\label{thm:simcon}Let $m,k\in\mathbb{N}$, $m,k>0$, $S$ a projective K3 surface, $v$ a Mukai vector on $S$ of type $(m,k)$ and $H$ a polarization that is general with respect to $v$. Then $M_{v}(S,H)$ and $M^{s}_{v}(S,H)$ are simply connected.
\end{thm}

\proof Suppose first that $m=1$, that $v$ is a Mukai vector such that if $v=(0,v_{1},v_{2})$ then $v_{2}\neq 0$, and that $H$ is $v-$generic. In this case we then have $M_{v}(S,H)=M_{v}^{s}(S,H)$: by \cite{OG1} and \cite{Y1} we know that $M_{v}(S,H)$ is an irreducible symplectic manifold, and we are done. 

If $v$ is any Mukai vector of type $(1,k)$ and $H$ is general with respect to $v$, the result follows from the case we considered above and by Lemmata \ref{lem:relgenvgen}, \ref{lem:tensorh0} and \ref{lem:polchangeab}.

Fix now $m\geq 2$ and $k\geq 1$. By point (1) of Theorem \ref{thm:main}, the moduli spaces arising from K3 surfaces, Mukai vectors of type $(m,k)$ and polarization which are general with respect to them are all deformation equivalent. As this deformation equivalence is obtained using only isomorphism of moduli spaces (coming from Fourier-Mukai transforms) and deformations of the moduli spaces induced by deformations of triples, by point (1) of Lemma \ref{lem:toploctriv} these deformation equivalent moduli spaces are also homeomorphic and the same holds for their stable (i. e. smooth) loci. 

It is then enough to prove that $M_{v}(S,H)$ and $M^{s}_{v}(S,H)$ are simply connected for one particular choice of $S$, $v$ and $H$ where $v$ is of type $(m,k)$ and $H$ is a general with respect to $v$: hence the result follows from  Proposition \ref{prop:simplu}.\endproof

\subsubsection{The case of Abelian surfaces}

Let $A$ be an Abelian surface with $NS(A)=\mathbb{Z}\cdot h$, where $h=c_{1}(H)$ and $H$ is an ample divisor such that $H^{2}=2k$. We let $h:=c_{1}(H)$, $m\in\mathbb{N}$ and $u:=m(0,h,0)$.

Let $Y_{mH}$ be the Hilbert scheme of curves on $A$ which are deformation of curves in $|mH|$, and let $p_{u}:M_{u}(A,H)\longrightarrow Y_{mH}$ be the morphism mapping a sheaf to its Fitting subscheme. We moreover let $p_{u}^{K}:K_{u}(A,H)\longrightarrow|mH|$ be the restriction of $p_{u}$ to $K_{u}(A,H)$.

We first prove  simple connectedness in the particular case of $K_{v}(A,H)$ and $K^{s}_{v}(A,H)$.

\begin{prop}
\label{prop:simpluab}If $(m,k)\neq(2,1)$, then $K_{u}(A,H)$ and $K^{s}_{u}(A,H)$ are simply connected.
\end{prop}

\proof Notice that $u$ is a Mukai vector of type $(m,k)$, and that $H$ is $u-$generic. If $m=1$, then $K_{u}(A,H)=K^{s}_{u}(A,H)$, and this is a point (if $k=1$) or an irreducible symplectic manifold (if $k>1$), and we are done.

For $m\geq 2$, we have that $K_{u}(A,H)$ is a normal, irreducible projective variety (see Remark A.1 of \cite{PR2}). As a consequence we have a surjective map $\pi_{1}(K_{u}^{s}(A,H))\longrightarrow\pi_{1}(K_{u}(A,H))$ (see Proposition 2.10 of \cite{Koll}), and it will be sufficient to prove that $K^{s}_{u}(A,H)$ is simply connected.

To show this, let $p^{K}_{u|K^{s}_{u}(A,H)}:K^{s}_{u}(A,H)\longrightarrow|mH|$ be the restriction of $p_{u}^{K}$ to $K^{s}_{u}$(A,H). By the Theorem in section 1.1, Part II of \cite{GM}, the fundamental group of a smooth connected variety admitting a dominant mapping to $\mathbb{P}^{N}$ (for some $N$) is generated by the fundamental group of the inverse image of a generic line in $\mathbb{P}^{N}$. As a consequence, if $\ell\subseteq|mH|$ is a generic line and $K^{0}:=(p^{K}_{u|K^{s}_{u}(A,H)})^{-1}(\ell)\subseteq K^{s}_{u}(A,H)$, we have a surjective morphism $\pi_{1}(K^{0})\longrightarrow\pi_{1}(K^{s}_{u}(A,H))$. It is then enough to show that $K^{0}$ is simply connected.

As $\ell$ is generic in $|mH|$, by Bertini's Theorem we know that $K^{0}$ is smooth. Moreover, by Lemma \ref{lem:codimr} all the curves parameterized by $\ell$ are reduced and irreducible. It then follows that $K^{0}=(p^{K}_{u})^{-1}(\ell)$.

To show that $K^{0}$ is simply connected, we show that $K^{0}$ is a fiber of an isotrivial fibration, and then use the homotopy exact sequence of this fibration to conclude. The domain of this isotrivial fibration will be $M^{0}:=p_{u}^{-1}(\ell)$ (which is a subset of $M_{u}(A,H)$), that will be identified with the relative compactified Jacobian of $\ell$. By construction, there is an inclusion $f:K^{0}\longrightarrow M^{0}$ fitting in a commutative diagram
$$\begin{CD}
K^{0} @>f>> M^{0}\\
@V\text{$p^{0}_{K}$}VV @VV\text{$p^{0}$}V\\
\ell @>>\text{$id_{\ell}$}> \ell
\end{CD}$$
where $p^{0}_{K}$ is the restriction of $p^{K}_{u}$ to $K^{0}$, and $p^{0}$ is the restriction of $p_{u}$ to $M^{0}$.

We now let $\sigma:M^{0}\longrightarrow A$ be the restriction to $M^{0}$ of the map $\beta:M_{u}(A,H)\longrightarrow A$ defined in section 2.2, mapping a sheaf $\mathcal{F}$ to the Albanese image of $\textbf{c}_{2}(\mathcal{F})$. As the determinant of $\mathcal{F}\in M^{0}$ is represented by the Fitting subscheme of $\mathcal{F}$, which is a divisor in $|mH|$, by Lemma \ref{lem:avbv} we have $$K^{0}=M^{0}\cap K_{u}(A,H)=M^{0}\cap b_{u}^{-1}(0_{A},\mathcal{O}_{A})=\sigma^{-1}(0_{A}),$$where $b_{u}:M_{u}(A,H)\longrightarrow A\times\widehat{A}$ is the O'Grady fibration of $M_{u}(A,H)$ defined in section 2.2.

Next, we claim that $\sigma:M^{0}\longrightarrow A$ is an isotrivial fibration. Indeed, if $L\in Pic^{0}(A)$ is represented by a divisor $D$, and $\delta$ is a $0-$cycle of degree 0 on $A$ representing $mH\cdot D$ in the Chow ring of $A$, then the tensorization with $L$ induces an automorphism of $M_{u}(A,H)$ mapping $K^{0}$ to $\sigma^{-1}(\delta)$. It follows that the connected algebraic group $Pic^{0}(A)$ acts transitively on the fibers of the projective morphism $\sigma$: this implies that $\sigma$ is an isotrivial fibration.

Finally, notice that $K^{0}$ is connected. Indeed, it is the inverse image, under the dominant map $p^{K}_{u}:K_{u}(A,H)\longrightarrow |mH|$, of a linear space of the projective space $|mH|$: by Theorem 1.1 of \cite{FL}, it follows that $K^{0}$ is connected. 

To resume, we have an isotrivial fibration $\sigma:M^{0}\longrightarrow A$, and $K^{0}$ is one of the fibers. The homotopy exact sequence associated to this fibration gives then $$\pi_{2}(A)\longrightarrow\pi_{1}(K^{0})\stackrel{\pi_{1}(f)}\longrightarrow\pi_{1}(M^{0})\stackrel{\pi_{1}(\sigma)}\longrightarrow\pi_{1}(A)\longrightarrow\{1\},$$where the last term comes from the fact that $K^{0}$ is connected. As $A$ is an Abelian surface, we have $\pi_{2}(A)=\{1\}$, hence in order to show that $K^{0}$ is simply connected, we just need to prove that the morphism $\pi_{1}(\sigma):\pi_{1}(M^{0})\longrightarrow\pi_{1}(A)$ is injective.

To do so, suppose that $\ell$ is generated by two smooth curves intersecting transversally at a finite number of points. Let $Bs(\ell)$ be the base locus of $\ell$, and $\pi:\widetilde{A}\longrightarrow A$ the blow-up of $A$ along $Bs(\ell)$.

The surface $\widetilde{A}$ is the total space of $\ell$: for every $a\in\widetilde{A}$ there is a unique curve $C_{a}\in\ell$ such that $a\in\widetilde{C}_{a}$, where $\widetilde{C}_{a}$ is the proper transform of $C_{a}$ under $\pi$. We then have a fibration $p_{\ell}:\widetilde{A}\longrightarrow\ell$, mapping $a\in\widetilde{A}$ to the point of $\ell$ corresponding to $C_{a}$.

There is a natural morphism $g:\widetilde{A}\longrightarrow M^{0}$ of fibrations over $\ell$ obtained as follows: first, choose $p\in Bs(\ell)$, and let $d:=m^{2}k+1$. For every $a\in\widetilde{A}$, the rank 1 torsion-free sheaf $\mathcal{I}_{g(a)}\otimes\mathcal{O}_{C_{a}}(dp)$ has degree $m^{2}k$. We then let $g(a):=\mathcal{I}_{g(a)}\otimes\mathcal{O}_{C_{a}}(dp)$, so to have a commutative diagram
$$\begin{CD}
\widetilde{A} @>\text{$g$}>> M^{0}\\
@V\text{$p_{\ell}$}VV @VV\text{$p^{0}$}V\\
\ell @>>\text{$id_{\ell}$}> \ell
\end{CD}$$

If $t\in\ell$ is a generic point, the curve $C$ corresponding to $t$ is smooth, $p_{\ell}^{-1}(t)=\widetilde{C}$ and $(p^{0})^{-1}(t)\simeq Pic^{m^{2}k}(C)$. Let $p_{1},\cdots,p_{n}\in\ell$ be the points corresponding to singular curves. The fundamental group of $M^{0}$ is generated by $\pi_{1}(Pic^{m^{2}k}(C))$ and by liftings $\widetilde{\gamma}_{1},\cdots,\widetilde{\gamma}_{n}$ of the generators $\gamma_{1},\cdots,\gamma_{n}$ of $\pi_{1}(\ell\setminus\{p_{1},\cdots,p_{n}\})$.

Moreover, the morphism $g_{t}:\widetilde{C}\longrightarrow Pic^{m^{2}k}(C)$ given by the restriction of $g$ to $p_{\ell}^{-1}(t)$ can be identified to the Abel-Jacobi map from $C$ to its Jacobian: it then induces a surjective map $\pi_{1}(\widetilde{C})\longrightarrow\pi_{1}(Pic^{m^{2}k}(C))$. 

As $\widetilde{C}\subseteq\widetilde{A}$, it follows that $\pi_{1}(M^{0})$ is generated by $\pi_{1}(\widetilde{A})$ and by the $\widetilde{\gamma}_{1},\cdots,\widetilde{\gamma}_{n}$. Now, notice that the fibration $p_{\ell}:\widetilde{A}\longrightarrow\ell$ has a section: fixing $p\in Bs(\ell)$, this section is obtained by mapping $t\in\ell$ to the unique intersection point of $\pi^{-1}(p)$ and $p_{\ell}^{-1}(t)$.  

We can then choose the liftings $\widetilde{\gamma}_{1},\cdots,\widetilde{\gamma}_{n}$ to be in the image of $\pi_{1}(\widetilde{A})$ in $\pi_{1}(M^{0})$, hence $g$ induces a surjection  $\pi_1(g):\pi_1(\widetilde{A})\rightarrow\pi_1(M^0)$. As $\sigma\circ g:\widetilde{A}\longrightarrow A$ induces an isomorphism between $\pi_{1}(\widetilde{A})$ and $\pi_{1}(A)$,  the morphism $\pi_1(\sigma)$ is injective and this concludes the proof.\endproof

Theorem \ref{thm:main} allows us to extend Proposition \ref{prop:simpluab} to all Abelian surfaces (provided $(m,k)\neq(2,1)$).

\begin{thm}
\label{thm:simconab}Let $m,k\in\mathbb{N}$, $m,k>0$, $S$ an Abelian surface, $v$ a Mukai vector on $S$ of type $(m,k)$ and $H$ a polarization which is general with respect to $v$.
\begin{enumerate}
 \item If $(m,k)\neq(2,1)$, then $K_{v}(S,H)$ and $K^{s}_{v}(S,H)$ are simply connected.
 \item If $(m,k)=(2,1)$, then $K_{v}(S,H)$ is simply connected, and $\pi_{1}(K^{s}_{v}(S,H))=\mathbb{Z}/2\mathbb{Z}$.
\end{enumerate}
\end{thm}

\proof Suppose first that $m=1$, that $v$ is a Mukai vector such that if $v=(0,v_{1},v_{2})$ then $v_{2}\neq 0$, and that $H$ is $v-$generic. Then $K_{v}(S,H)=K^{s}_{v}(S,H)$, and this is a point (if $k=1$) or an irreducible symplectic manifold (if $k\geq 2$): the statement is then clear in this case.

If $v$ is any Mukai vector of type $(1,k)$ and $H$ is general with respect to $v$, the result follows from the case we considered above and by Lemmata \ref{lem:relgenvgen}, \ref{lem:tensorh0} and \ref{lem:polchangeab}.

Fix now $m\geq 2$ and $k\geq 1$, and suppose that $(m,k)\neq(2,1)$. By point (2) of Theorem \ref{thm:main}, the fibers of the Yoshioka fibration of a moduli spaces arising from Abelian surfaces, Mukai vectors of type $(m,k)$ and polarizations which are general with respect to them are all deformation equivalent. As this deformation equivalence is obtained using only isomorphism of moduli spaces (coming from Fourier-Mukai transforms) and deformations of the moduli spaces induced by deformations of triples, by point (2) of Lemma \ref{lem:toploctriv} the homeomorphism type of $K_{mw}(S,H)$ and $K^s_{mw}(S,H)$ only depends on $m$ and $k=w^{2}/2$. 

It is then enough to show that $K_{v}(S,H)$ and $K^{s}_{v}(S,H)$ are simply connected for a particular choice of $S$, $v$ of type $(m,k)$ and $H$ which is general with respect to $v$: the result follows then from Proposition \ref{prop:simpluab}.

If $(m,k)=(2,1)$, by Lemmata \ref{lem:relgenvgen}, \ref{lem:tensorh0} and \ref{lem:polchangeab} we may suppose that $v$ is a Mukai vector such that if $v=(0,v_{1},v_{2})$ then $v_{2}\neq 0$, and that $H$ is $v-$generic. 

In this case we know that $K_{v}(S,H)$ admits a symplectic resolution $\widetilde{K}_{v}(S,H)$, which is an irreducible symplectic manifold by point (2) of Theorem 1.6 in \cite{PR}. As $K_{v}(S,H)$ has canonical singularities, by \cite{T} we have $\pi_{1}(K_{v}(S,H))=\pi_{1}(\widetilde{K}_{v}(S,H))$: it follows that $K_{v}(S,H)$ is simply connected. 

By Theorem 4.2 and Proposition 5.3 of \cite{MRS}, we know that $K^{s}_{v}(S,H)$ has an \'etale cover of degree 2 from an open subset $U$ of an irreducible symplectic manifold $Y$ which is deformation equivalent to a Hilbert scheme of 3 points on a K3 surface. This open subset $U$ is obtained by removing from $Y$ 256 copies of $\mathbb{P}^{3}$ and one copy of a desingularization of the singular locus of $K_{v}(S,H)$. It follows that the complement of $U$ has codimension at least 2 in $Y$, so that $\pi_{1}(U)=\pi_{1}(Y)=\{1\}$. It then follows that the fundamental group of $K^{s}_{v}(S,H)$ is $\mathbb{Z}/2\mathbb{Z}$.\endproof

Recall that if $X$ is a normal projective variety having at most rational singularities, it is possible to define the Albanese variety $Alb(X)$ as the Albanese variety of any desingularization $\widetilde{X}$ of $X$, and  construct the Albanese morphism $alb:X\longrightarrow Alb(X)$ by descending the usual Albanese morphism of $\widetilde{X}$ (see Proposition 2.3 of \cite{Re}, and Lemma 8.1 of \cite{Kaw}). 

As a consequence of Theorem \ref{thm:simconab}, we show in the next result that the Yoshioka fibration is the Albanese morphism of the moduli space $M_{v}(S,H)$.

\begin{cor}
\label{cor:albanese}
Let $m,k\in\mathbb{N}$, $m,k>0$, $S$ an Abelian surface, $v$ a Mukai vector on $S$ of type $(m,k)$ and $H$ a polarization that is general with respect to $v$. The morphism $a_{v}:M_{v}(S,H)\longrightarrow S\times\widehat{S}$ is the Albanese morphism of $M_{v}(S,H)$.
\end{cor}

\proof By Lemmata \ref{lem:vgengen}, \ref{lem:relgenvgen}, \ref{lem:tensorh0} and \ref{lem:tensorlb} we may assume that $v$ is not of the form $(0,v_1,0)$ and that $H$ is a $v-$generic polarization. Under this assumption, for $m=1$ the map $a_{v}$ is the Albanese map by point (1) of Theorem 0.1 in \cite{Y2}. We then suppose $m\geq 2$.

For  $(m,k)=(2,1)$, we know by Th\'eor\`eme 1.1 of \cite{LS} that $M_{v}(S,H)$ admits a symplectic resolution of the singularities $\pi:\widetilde{M}_{v}(S,H)\longrightarrow M_{v}(S,H)$, which is obtained by blowing up the singular locus $\Sigma$ with reduced structure.

Now, for every $(p,L)\in S\times\widehat{S}$ the fiber $K_{p,L}:=a_{v}^{-1}(p,L)$ is a singular symplectic variety whose singular locus is $\Sigma_{p,L}:=\Sigma\cap K_{p,L}$, and $\widetilde{K}_{p,L}:=\pi^{-1}(K_{p,L})$ is the symplectic resolution of $K_{p,L}$, which is an irreducible symplectic manifold by point (2) of Theorem 1.6 of \cite{PR}. It follows that $$a_{v}\circ\pi:\widetilde{M}_{v}(S,H)\longrightarrow S\times\widehat{S}$$is the Albanese morphism of $\widetilde{M}$, so that $a_{v}:M_{v}(S,H)\longrightarrow S\times\widehat{S}$ is the Albanese morphism of $M_{v}(S,H)$.

Let us now finally consider the case $(m,k)\neq(2,1)$, and $m\geq 2$. Let $\pi:\widetilde{M}\longrightarrow M_{v}(S,H)$ be a desingularization of $M_{v}(S,H)$, where $\widetilde{M}$ is a smooth projective variety. The inclusion $j:M^{s}_{v}(S,H)\longrightarrow\widetilde{M}$ induces a surjective morphism $\pi_{1}(j):\pi_{1}(M^{s}_{v}(S,H))\longrightarrow\pi_{1}(\widetilde{M})$. If we now let $a^{s}_{v}:M^{s}_{v}(S,H)\longrightarrow S\times\widehat{S}$ be the restriction of $a_{v}$ to $M^{s}_{v}(S,H)$, then $a^{s}_{v}$ is an isotrivial fibration whose fibers are all isomorphic to $K^{s}_{v}(S,H)$.

As the fiber  $K^{s}_{v}(S,H)$ is simply connected (since $(m,k)\neq(2,1)$), the isotrivial fibration $a^{s}_{v}$ induces an isomorphism $\pi_{1}(a^{s}_{v}):\pi_{1}(M^{s}_{v}(S,H))\longrightarrow\pi_{1}(S\times\widehat{S})$. 

Now, notice that $a_{v}\circ\pi\circ j=a^{s}_{v}$, hence $\pi_{1}(j)$ is injective, and hence an isomorphism. But this implies that $\pi_{1}(a_{v}\circ\pi)$ is an isomorphism, so that $a_{v}\circ\pi:\widetilde{M}\longrightarrow S\times\widehat{S}$ is the Albanese morphism for $\widetilde{M}$. It then follows that $a_{v}:M_{v}(S,H)\longrightarrow S\times\widehat{S}$ is the Albanese morphism of $M_{v}(S,H)$, concluding the proof.\endproof

\subsection{The proof of Theorem \ref{thm:mainb1}}

We are finally in the position to prove Theorem \ref{thm:mainb1}, i. e. that $M_{v}(S,H)$ and $K_{v}(S,H)$ are irreducible symplectic varieties. Before doing this, we calculate the dimension of the space of reflexive $p-$forms for a particular choice of the surface, of the Mukai vector and of the polarization.

\begin{lem}
\label{lem:hppart}Let $m,k\in\mathbb{N}$, $m,k>0$, and $X$ a projective K3 or Abelian surface such that $NS(X)=\mathbb{Z}\cdot h_k$ where $h_k^{2}=2k$. Let $H_k$ be a polarization on $S$ such that $c_{1}(H_k)=h_k$, and let $u_k=m(0,h_k,0)$.
\begin{enumerate}
 \item If $X$ is K3 and $p\in\mathbb{N}$ is such that $0\leq p\leq\dim(M_{u_k}(X,H_k))$, then $$h^{0}(M_{u_k}(X,H_k),\Omega^{[p]}_{M_{u_k}(X,H_k)})=\left\{\begin{array}{ll} 1 & p\,\,is\,\,even\\0 & p\,\,is\,\,odd\end{array}\right. .$$
 \item If $X$ is Abelian and $p\in\mathbb{N}$ is such that $0\leq p\leq\dim(K_{u_k}(X,H_k))$, then $$h^{0}(K_{u_k}(X,H_k),\Omega^{[p]}_{K_{u_k}(X,H_k)})=\left\{\begin{array}{ll} 1 & p\,\,is\,\,even\\0 & p\,\,is\,\,odd\end{array}\right. .$$
\end{enumerate}
\end{lem}

\proof Suppose first that $X$ is K3. We let $u:=(0,mh_k,1-m^{2}k)$, which is a primitive Mukai vector on $X$.

If $C\in|mH_k|$ is an integral curve and $j:C\longrightarrow X$ is the inclusion, for every $L\in Pic^{1}(C)$ the sheaf $j_{*}L$ is $H_k-$stable of Mukai vector $u$. The sheaves of this type form an open subset $U$ of $M_{u}(X,H_k)$. 

Moreover, if $L\in Pic^{1}(C)$ then $L^{\otimes m^{2}k}\in Pic^{m^{2}k}(C)$, hence $j_{*}(L^{\otimes m^{2}k})$ is an $H_k-$stable sheaf of Mukai vector $v$. We then have a rational map $$g:M_{u}(X,H_k)\dashrightarrow M_{u_k}(X,H_k),\,\,\,\,\,\,\,\,\,g(j_{*}L):=j_{*}L^{\otimes m^{2}k}.$$

We first show that $g$ is dominant. To do so, consider the two fibrations $p_{u}:M_{u}(X,H_k)\longrightarrow|mH_k|$ and $p_{u_k}:M_{u_k}(X,H_k)\longrightarrow|mH_k|$ mapping a sheaf to its Fitting subscheme. If $C\in|mH_k|$ is smooth, we have $p_{u}^{-1}(C)\simeq Pic^{1}(C)$ and $p_{u_k}^{-1}(C)\simeq Pic^{m^{2}k}(C)$, hence $p_{u}^{-1}(C)\simeq p_{u_k}^{-1}(C)\simeq Pic^{0}(C)$ and the restriction of $g$ to $p_{u}^{-1}(C)$ can be identified with the multiplication by $m^{2}k$ on $Pic^{0}(C)$: hence it is surjective. This shows that, if $V\subseteq|mH_k|$ is the open subset of smooth curves, then $g$ maps $p_{u}^{-1}(V)$ surjectively to $p_{u_k}^{-1}(V)$. As $M_{u}(X,H_k)$ and $M_{u_k}(X,H_k)$ are two projective varieties which are both irreducible and of the same dimension, it follows that $g$ is dominant.

Since $M_{u_k}(X,H_k)$ has canonical singularities, letting  $\widetilde{M}_{u_k}$ be a resolution of the singularities, by Theorem 1.4 of \cite{GKKP} we have the equality $h^{0}(M_{u_k}(X,H_k),\Omega_{M_{u_k}(X,H_k)}^{[p]})=h^{0}(\widetilde{M}_{u_k},\Omega_{\widetilde{M}_{u_k}}^{p})$  for every $p\in\mathbb{N}$. As $g$ is dominant, we also have  $h^{0}(\widetilde{M}_{u_k},\Omega^{p}_{\widetilde{M}_{u_k}})\leq h^{0}(M_{u},\Omega^{p}_{M_{u}(X,H_k)})$ for every $p\in\mathbb{N}$. Since $u$ is primitive, $M_{u}(X,H_k)$ is an irreducible symplectic manifold and we conclude that $$h^{0}(M_{u_k}(X,H_k),\Omega^{[p]}_{M_{u_k}(X,H_k)})\le\left\{\begin{array}{ll} 1 & p\,\,is\,\,even\\0 & p\,\,is\,\,odd\end{array}\right. .$$
Since $M_{u_k}(X,H_k)$ is a symplectic variety  $h^{0}(M_{u_k}(X,H_k),\Omega^{[p]}_{M_{u_k}(X,H_k)})\geq 1$ if $p$ is even and we are done.

If $S$ is Abelian, the same proof works replacing $M_{u_k}(X,H_k)$ with $K_{u_k}(X,H_k)$ and $M_{u}(X,H_k)$ with $K_{u}(X,H_k)$.\endproof

We now prove the following:\vspace{3mm}

\noindent\textbf{Theorem \ref{thm:mainb1}.} \textit{Let $m,k\in\mathbb{N}$, $m,k>0$, $S$ a projective K3 surface, $v$ a Mukai vector on $S$ of type $(m,k)$ and $H$ a polarization that is general with respect to $v$. 
\begin{enumerate}
 \item If $S$ is K3, then $M_{v}(S,H)$ is an irreducible symplectic variety.
 \item If $S$ is Abelian and $(m,k)\neq(1,1)$, then $K_{v}(S,H)$ is an irreducible symplectic variety.
\end{enumerate}}

\proof We first consider $S$ to be a K3 surface. If $m=1$, then $M_{v}(S,H)$ is an irreducible symplectic manifold by Theorem 0.1 of \cite{Y1} and by Lemmata
\ref{lem:relgenvgen}, \ref{lem:tensorh0} and \ref{lem:polchangeab} if $H$ is not $v-$generic or $v=(0,v_1,0)$.

If $m\geq 2$, then $M_{v}(S,H)$ is a symplectic variety, and let $\sigma$ be a symplectic form on it. We have to show that if $f:Y\longrightarrow M_{v}(S,H)$ is a finite quasi-\'etale morphism, then the exterior algebra of reflexive forms on on the normal variety $Y$ is spanned by $f^{[*]}\sigma$.

Let then $f:Y\longrightarrow M_{v}(S,H)$ be a finite quasi-\'etale cover: it induces a finite quasi-\'etale cover of $M^{s}_{v}(S,H)$. But a finite quasi-\'etale morphism of a smooth variety is \'etale, and $M^{s}_{v}(S,H)$ is simply connected by Theorem \ref{thm:simcon}: hence $f$ is an isomorphism.

We then just need to show that the exterior algebra of reflexive forms on $M_{v}(S,H)$ if spanned by $\sigma$. This follows if we show that $h^{0}(M_{v}(S,H),\Omega^{[p]}_{M_{v}(S,H)})=1$ if $p$ is even, and $h^{0}(M_{v}(S,H),\Omega^{[p]}_{M_{v}(S,H)})=0$ if $p$ is odd.

For this, let $X$ be a projective K3 surface with $Pic(X)=\mathbb{Z}\cdot H_k$, where $H_k$ is an ample line bundle with $(H_k)^{2}=2k$, and let $u_k:=m(0,h_k,0)$, where $h_k=c_{1}(H_k)$.
By Lemma \ref{lem:hppart} we have $$h^{0}(M_{u_k}(X,H_k),\Omega^{[p]}_{M_{u_k}(X,H_k)})=\left\{\begin{array}{ll} 1 & p\,\,is\,\,even\\0 & p\,\,is\,\,odd\end{array}\right. . $$

By Theorem \ref{thm:main} $M_{v}(S,H)$ and $M_{u_k}(X,H_k)$ are deformation equivalent, and the deformation is locally trivial. Hence they have resolutions $\widetilde{M}_{v}$ and $\widetilde{M}_{u_k}$ of the singularities which are deformation equivalent as smooth varieties, so their Hodge numbers are equal. By Theorem 1.4 of \cite{GKKP} we then have $$h^{0}(M_{v}(S,H),\Omega^{[2]}_{M_{v}(S,H)})=h^{0}(\widetilde{M}_{v},\Omega^{2}_{\widetilde{M}_{v}})=$$$$h^{0}(\widetilde{M}_{u_k},\Omega^{2}_{\widetilde{M}_{u_k}})=h^{0}(M_{u_k}(X,H_k),\Omega^{[2]}_{M_{u_k}(X,H_k)}),$$and we are done.

If $S$ is an Abelian surface, the proof is identical if $(m,k)\neq(2,1)$, replacing moduli spaces of sheaves by the corresponding Albanese fibers and using Theorem \ref{thm:simconab}(1) instead of Theorem \ref{thm:simcon}. 
 
The case $(m,k)=(2,1)$ has to be treated differently: by Theorem \ref{thm:simconab}(2) we have $\pi_{1}(K_{v}^{s}(S,H))=\mathbb{Z}/2\mathbb{Z}$, hence $K_{v}(S,H)$ has a unique (up to isomorphism) nontrivial connected finite quasi-\'etale cover $Y_{v}$. 
We need to show that the exterior algebra of reflexive forms on $K_{v}(S,H)$ and $Y_{v}$ are spanned by the reflexive pull-back of a symplectic form on $K_{v}(S,H)$. To do so, it will be enough to show that both $K_{v}(S,H)$ and $Y_{v}$ are birational to  irreducible symplectic manifolds.

For $K_{v}(S,H)$, by point (2) of Theorem 1.6 in \cite{PR} we know that $K_{v}(S,H)$ has a symplectic resolution which is an irreducible symplectic manifold (in the deformation class $OG_{6}$). For $Y_{v}$, by Proposition 5.3 of \cite{MRS} we know that it is birational to an irreducible symplectic manifold deformation equivalent to $Hilb^{3}(K3)$. This concludes the proof.\endproof

\end{document}